\documentclass[11pt]{amsart}
\usepackage{amsmath,amsthm,color,amssymb}
\usepackage{fullpage}
\usepackage{rotating}
\makeindex
\usepackage{wrapfig,epsfig}

\input xy
\xyoption{all}

\theoremstyle{plain}
\newtheorem{thm}{Theorem}[section]
\newtheorem{lem}[thm]{Lemma}
\newtheorem{cor}[thm]{Corollary}
\newtheorem{prop}[thm]{Proposition}
\newtheorem{rem}[thm]{Remark}

\newcommand{\eh}{\mathbb{A}}

\newcommand{\eff}{\mathbb{F}}
\newcommand{\gee}{\mathbb{G}}
\newcommand{\en}{\mathbb{N}}
\newcommand{\pe}{\mathbb{P}}
\newcommand{\cue}{\mathbb{Q}}
\newcommand{\arr}{\mathbb{R}}
\newcommand{\zed}{\mathbb{Z}}

\newcommand{\oh}{\mathcal{O}}

\newcommand{\Ann}{\operatorname{Ann}}
\newcommand{\Aut}{\operatorname{Aut}}

\newcommand{\bs}{\backslash}

\newcommand{\coker}{\operatorname{coker}}
\newcommand{\Comm}{\operatorname{Comm}}

\newcommand{\diag}{\operatorname{diag}}
\newcommand{\dirlim}{\underrightarrow{\lim}\,}
\newcommand{\End}{\operatorname{End}}
\newcommand{\fp}{\mathbb{F}_p}
\newcommand{\fpbar}{\overline{\mathbb{F}}_p}
\newcommand{\fptwo}{\mathbb{F}_{p^2}}
\newcommand{\fq}{\mathbb{F}_q}
\newcommand{\Fr}{\operatorname{Fr}}
\newcommand{\Frob}{\operatorname{Frob}}
\newcommand{\Gal}{\operatorname{Gal}}
\newcommand{\GL}{\operatorname{GL}}
\newcommand{\GSp}{\operatorname{GSp}}
\newcommand{\GU}{\operatorname{GU}}
\renewcommand{\H}{\operatorname{H}}
\newcommand{\Hom}{\operatorname{Hom}}
\newcommand{\im}{\operatorname{Im}}

\newcommand{\into}{\hookrightarrow}

\newcommand{\invlim}{\underleftarrow{\lim}\,}

\newcommand{\ord}{\operatorname{ord}}
\newcommand{\Pic}{\operatorname{Pic}}

\newcommand{\R}{\operatorname{R}}

\newcommand{\Res}{\operatorname{Res}}

\newcommand{\SL}{\operatorname{SL}}
\newcommand{\Sp}{\operatorname{Sp}}
\newcommand{\Spec}{\operatorname{Spec}}

\newcommand{\std}{\operatorname{std}}
\newcommand{\Sym}{\operatorname{Sym}}

\newcommand{\U}{\operatorname{U}}

\begin{document}
\title[Siegel modular forms]{Siegel modular forms (mod $p$) and algebraic modular forms}
\author[A. Ghitza]{Alexandru Ghitza}
\date{}
\thanks{This is the author's PhD thesis, submitted to the Department of Mathematics of the Massachusetts Institute of Technology in June 2003.}
\address{Department of Mathematics\\
Massachusetts Institute of Technology\\
Cambridge, MA 02139}
\email{aghitza@math.mit.edu}
\begin{abstract}
In his letter~\cite{serre1}, J.-P. Serre proves that the systems of Hecke eigenvalues given by modular forms (mod $p$) are the same as the ones given by locally constant functions $\mathbb{A}_B^\times/B^\times\to\bar{\mathbb{F}}_p$, where $B$ is the endomorphism algebra of a supersingular elliptic curve.  After giving a detailed exposition of Serre's result, we prove that the systems of Hecke eigenvalues given by Siegel modular forms (mod $p$) of genus $g$ are the same as the ones given by algebraic modular forms (mod $p$) on the group $\text{GU}_g(B)$, as defined in~\cite{gross2} and~\cite{gross4}.  The correspondence is obtained by restricting to the superspecial locus of the moduli space of abelian varieties.
\end{abstract}

\maketitle
\pagestyle{headings}
\tableofcontents

\section{Introduction}
The study of arithmetic properties of modular forms goes back to Ramanujan, who observed that the coefficients $\tau(n)$ of the modular form
$$
\Delta=q\prod_{n=1}^\infty (1-q^n)^{24}=\sum_{n=1}^\infty\tau(n)q^n
$$
satisfy interesting congruence relations modulo various primes.  After these congruence relations had been proved by classical methods, Swinnerton-Dyer used results of Serre and Deligne in order to conduct a systematic study of the problem.  His main tool is the description of the algebra of modular forms mod $p$ (see~\cite{swinnerton-dyer1},~\cite{swinnerton-dyer2}). 

Consider a modular form $f$ of weight $k$ and level $N$
$$
f=\sum_{n=0}^\infty a_nq^n
$$
whose coefficients $a_n$ are rational and $p$-integral (that is, $p$ does not divide the denominator of $a_n$).  Write $\tilde{f}$ for the reduction of $f$ modulo $p$, i.e. the element of $\fp[[q]]$ obtained by reducing each coefficient of $f$ modulo $p$.  We denote the set of all such series $M_k(N)$.  Set $M(N)=\sum_{k\in\zed} M_k(N)$, known as the algebra of level $N$ modular forms (mod $p$).  For $\ell$ prime to $pN$, one defines the Hecke operator $T_{\ell}$ by 
$$
T_{\ell} f=\sum_{n=0}^\infty a_{\ell n}q^n+\ell^{k-1}\sum_{n=0}^\infty a_nq^{\ell n}.
$$

An eigenform is an element $f$ of $M_k(N)$ which is a simultaneous eigenvector for all $T_{\ell}$, $\ell\nmid pN$.  Given an eigenform, one obtains a system of eigenvalues $(b_{\ell})_{\ell\nmid pN}$ defined simply by $T_{\ell} f=b_{\ell} f$. 

A considerable amount of effort has gone into the study of this Hecke action on modular forms.  One of the approaches, due to Serre and inspired by the philosophy of the Langlands program, is to interpret these modular forms in terms of ad\`eles on quaternion algebras.  More precisely, let $B$ be the quaternion algebra ramified at $p$ and $\infty$ and let $\eh_B^\times$ denote the group of ad\`eles on $B$.  Serre proved the following
\begin{thm}[\cite{serre1}]\label{thm:serre}
The systems of eigenvalues $(b_{\ell})$ (with $b_{\ell}\in\fpbar$) given by the Hecke operators on modular forms (mod $p$) of any level $N$ are the same as those obtained from locally constant functions $f:\eh_B^\times/B^\times\to\fpbar$.
\end{thm}

Given an eigenform with eigenvalues $(b_{\ell})$, Deligne shows in~\cite{deligne3} (see also~\cite{gross5}) how to construct a mod $p$ Galois representation of degree $2$
$$
\rho:\Gal(\overline{\cue}/\cue)\to\GL_2(\fpbar)
$$
such that $\rho(\Frob_{\ell})$ has trace $b_{\ell}$ and determinant $\ell^{k-1}$.  Serre shows therefore that the systems of traces of Frobenii for a modular representation are the same as the eigenvalues of locally constant functions $\eh_B^\times/B^\times\to\fpbar$.  We remark that he was lead to this result by working on his famous conjectures, which basically say that any Galois representation of degree~$2$ that satisfies certain properties comes from a modular form in the way we mentioned (see~\cite{serre5}).

Inspired by this work, Gross has introduced (\cite{gross2}) the notion of algebraic modular forms (mod $p$) on certain reductive algebraic groups $G$ over $\cue$.  The hope is to establish a relation between these forms and higher-dimensional Galois representations
$$
\rho:\Gal(\bar{\cue}/\cue)\to \hat{G}(\fpbar),
$$
where $\hat{G}$ is the $L$-group of $G$ (see~\cite{gross4} for a precise statement).  Serre's result gives one direction of this relation in the case $G=B^\times$, $\hat{G}=\GL_2$, while the opposite direction is the object of the Serre conjectures.

The aim of this thesis is to generalize Serre's work to the case $G=\GU_g(B)=$ the group of unitary similitudes of the quaternionic hermitian space $B^g$, for $g>1$.  We start by giving a detailed and explicit treatment of Serre's theorem in \S~\ref{chap:elliptic}.  After recalling some basic facts about modular forms, we reduce Theorem~\ref{thm:serre} to the existence of a Hecke-invariant bijection between a finite subset of the modular curve $X(N)$ and a finite double coset space.  The reduction is based on restricting modular forms to the supersingular locus, and on the use of the Hasse invariant to prove that this restriction is compatible with the Hecke action.  The construction of the bijection of finite sets is explicit but somewhat lengthy and takes the remainder of the section.

We attack the generalization to higher dimensions in \S~\ref{chap:abelian}.  After recalling the definition of Siegel modular forms, we briefly discuss properties of superspecial abelian varieties and of the algebraic group $\GU_g(B)$.  The main result is the following
\begin{thm}
Fix a dimension $g>1$, a level $N\geq 3$ and a prime $p$ not dividing $N$.  The systems of Hecke eigenvalues coming from Siegel modular forms (mod $p$) of dimension $g$, level $N$ and any weight $\rho$, are the same as the systems of Hecke eigenvalues coming from algebraic modular forms (mod $p$) of level $N$ and any weight $\rho_S$ on the group $\GU_g(B)$.
\end{thm}

The proof consists of two main steps: first we show that the restriction to the locus of superspecial principally polarized abelian varieties preserves the systems of eigenvalues, and that systems of eigenvalues occuring on the superspecial locus can be lifted to the entire moduli space; we then construct a Hecke bijection between the superspecial locus and a double coset space built from $\GU_g(B)$.

\section{Elliptic modular forms}\label{chap:elliptic}

We give a detailed exposition of Serre's letter to Tate (see~\cite{serre1}) linking elliptic modular forms (mod $p$) to quaternion algebras.

\subsection{The geometric theory of modular forms (mod $p$)}
We review some of the definitions and results from Chapter~1 of~\cite{katz1}.  The reader unfamiliar with the geometric definition of modular forms is encouraged to consult Katz' article for details.

An \emph{elliptic curve}\index{elliptic curve} over a scheme $S$ is a proper smooth morphism
$\pi:E\to S$, whose geometric fibers are connected curves of genus
one, together with a section $0:S\to E$:
$$
\UseTips
\xymatrix{
E \ar[d]_{\pi}\\
S \ar@/_.8pc/[u]_{0}
}
$$
We define $\omega_{E/S}:=\pi_*(\Omega_{E/S}^1)$.  This is an invertible
sheaf on $S$; by Serre duality it is canonically dual to
$\R^1\pi_*(\oh_E)$.

For each integer $N\geq 1$, $E[N]:=\ker([N]:E\to E)$ is a finite flat
commutative group scheme of rank $N^2$ over $S$.  It is \'etale over
$S$ if and only if $S$ is a scheme over $\zed[\frac{1}{N}]$.  In this
case, we define a \emph{level $N$ structure}\index{level structure} on $E$ to be an isomorphism
$$
\alpha:E[N]\xrightarrow{\sim}(\zed/N\zed)^2_S.
$$
If such an isomorphism exists and $S$ is connected, then the set of
all level $N$ structures is principal homogeneous under
$\Aut((\zed/N\zed)^2_S)=\GL_2(\zed/N\zed)$.

For $N\geq 3$, the functor ``isomorphism classes of elliptic curves
with level $N$ structure'' is representable by a scheme $Y(N)$ which
is an affine smooth curve over $\zed[\frac{1}{N}]$, finite and flat of
degree equal to $\#(\GL_2(\zed/N\zed)/\{\pm 1\})$ over the affine
$j$-line $\zed[\frac{1}{N},j]$, and \'etale over the open set of the
affine $j$-line where $j$ and $j-1728$ are invertible.  The
normalization of the projective $j$-line $\pe^1_{\zed[\frac{1}{N}]}$
in $Y(N)$ is a proper and smooth curve $X(N)$ over
$\zed[\frac{1}{N}]$, the global sections of whose structure sheaf are
$\zed[\frac{1}{N},\zeta_N]$.  Let $(\mathcal{E}/Y(N),\alpha)$ be the
universal elliptic curve with level $N$ structure.  There is a unique
invertible sheaf $\omega$ on $X(N)$ whose restriction to
$Y(N)$ is $\omega_{\mathcal{E}/Y(N)}$ and whose sections
over the completion $\zed[\frac{1}{N},\zeta_N][[q]]$ at each cusp are
precisely the $\zed[\frac{1}{N},\zeta_N][[q]]$ multiples of the
canonical differential of the Tate curve.  The \emph{Kodaira-Spencer isomorphism}\index{Kodaira-Spencer map}
$$
\left(\omega_{\mathcal{E}/Y(N)}\right)^2 \cong
\Omega^1_{Y(N)/\zed[\frac{1}{N}]}
$$
extends to an isomorphism
$$
\omega^2 \cong \Omega^1_{X(N)/\zed[\frac{1}{N}]}(\log(X(N)-Y(N))). 
$$

A \emph{modular form (mod $p$) of level $N$ and weight $k$}\index{modular form} is a global section of $\omega^k$ on $X(N)\otimes\fpbar$ or equivalently a global section of the quasi-coherent sheaf $\omega^k\otimes\fpbar$ on $X(N)$.

An important fact in the mod $p$ theory is the existence of a modular form $A$ of level $1$ and weight $p-1$, called the \emph{Hasse invariant}\index{Hasse invariant}.  For its properties, see \S{}IV.4 in~\cite{hartshorne1}, \S{}2.0 and \S{}2.1 in~\cite{katz1}.

\subsection{Main result and reductions}\label{sect:main_results}
Fix a prime $p$ and an integer $N\geq 3$ prime to $p$.  Set
$X:=X(N)\otimes\fpbar$.

We work with Katz' definition of modular forms (mod $p$) of level $N$
and weight $k$: 
$$
M_k(N)=\H^0(X,\omega^k).
$$

Multiplication by the Hasse invariant $A$ gives a natural Hecke embedding
$$
M_{k-(p-1)}(N)\into M_k(N)
$$
and we're interested in the structure of the quotient
$W_k(N)=M_k(N)/M_{k-(p-1)}(N)$ as a module over the Hecke operators
$T_{\ell}$, $\ell\nmid pN$.  On the level of sheaves, we define
$\mathcal{S}_k$ by the following 
$$
0\to\omega^{k-(p-1)}\xrightarrow{\times A}\omega^k\to\mathcal{S}_k\to 0.
$$
The Hasse invariant $A$ vanishes to order $1$ when evaluated on the supersingular elliptic curves (and does not vanish anywhere else), so $\mathcal{S}_k$ is $0$ over the ordinary locus and one-dimensional over the supersingular points.

We take global sections and get the exact sequence
$$
0\to W_k(N)\to S_k(N)\to \H^1(X,\omega^{k-(p-1)})\to \H^1(X,\omega^k)\to 0.
$$
Actually by Theorem 1.7.1 in~\cite{katz1}, $\H^1(X,\omega^{k-(p-1)})=0$ if $k\geq p+1$, so $W_k(N)=S_k(N)$ when $k\geq p+1$.

An \emph{$\eff_q$-structure}\index{Fq-structure@$\eff_q$-structure} on an elliptic curve $E/\fpbar$ is an elliptic curve $E'/\eff_q$ such that $E$ is isomorphic to $E'\otimes\fpbar$.
\begin{lem}\label{lem:canonical_supersing}
Let $E$ be a supersingular elliptic curve over $\fpbar$.  Then $E$ has
a canonical $\fptwo$-structure $E'$, namely the one whose geometric
Frobenius is $[-p]$.  The correspondence $E\mapsto E'$ is
functorial.
\end{lem}
\begin{proof}
It is well-known that $E$ is supersingular if and only if $[-p]:E\to
E$ is purely inseparable (see for instance Theorem~V.3.1
of~\cite{silverman1}).  Multiplication by $-p$ is an isogeny of degree $p^2$, so by Corollary~II.2.12 of~\cite{silverman1} there is a commutative triangle
$$
\UseTips
\xymatrix{
E\ar[r]^{[-p]}\ar[dr]_{\Fr^2} & E\\
& E^{(p^2)}\ar[u]_{\sim},
}
$$
where $E^{(p^2)}$ is obtained by raising to the power $p^2$ the
coefficients of a Weierstrass equation for $E$ and $\Fr^2:E\to
E^{(p^2)}$ is the morphism $[x:y:z]\mapsto [x^{p^2}:y^{p^2}:z^{p^2}]$.
Since $E\cong E^{(p^2)}$ we have
$j(E)=j\left(E^{(p^2)}\right)=j(E)^{p^2}$, i.e. $j(E)\in\fptwo$.
Therefore there exists a supersingular elliptic curve $E''/\fptwo$
such that $j(E'')=j(E)$.  In other words, $E$ has an
$\fptwo$-structure.

The curve $E''$ is supersingular so the morphism $[-p]:E''\to E''$ is purely
inseparable.  As before we have a commutative diagram
$$
\UseTips
\xymatrix{
E''\ar[r]^{[-p]}\ar[dr]_{\pi''} & E''\\
& E''\ar[u]^{\sim}_{\lambda},
}
$$
where $\pi'':E''\to E''$ is the geometric Frobenius
$[x:y:z]\mapsto[x^{p^2}:y^{p^2}:z^{p^2}]$, and $\lambda\in\Aut(E'')$.
Therefore
$$
\pi''=\lambda^{-1}\circ [-p]=[-p]\circ\lambda^{-1}.
$$

As usual, we write $G_{\fptwo}$ to denote $\Gal(\fpbar/\fptwo)$.  Let
$\sigma\in G_{\fptwo}$ be the Frobenius element, defined by $x\mapsto
x^{p^2}$; it is a topological generator of $G_{\fptwo}$.  Consider the
continuous $1$-cocycle $\xi:G_{\fptwo}\to\Aut(E'')$ defined by
$\sigma\mapsto\lambda^{-1}$.  Since $\lambda^{-1}$ is a morphism
defined over $\fptwo$, $\sigma$ acts trivially on $\lambda^{-1}$ and
$\xi$ is a group homomorphism.  Using a combination of Theorem~X.2.2
and Proposition~X.5.3 in~\cite{silverman1}, we know that given
$\xi\in\H^1(G_{\fptwo},\Aut(E''))$ one can construct an elliptic curve
$E'/\fptwo$ and an isomorphism $\phi:E'\otimes\fpbar\to
E''\otimes\fpbar$ satisfying
$$
\xi(\tau)=\phi^\tau\circ\phi^{-1},\quad\text{ for all }\tau\in
G_{\fptwo}.
$$
We recall that the action of $G_{\fptwo}$ on $\fpbar$-morphisms is
given by
$$
\phi^\tau(P)=\left(\phi\left(P^{\tau^{-1}}\right)\right)^\tau.
$$

We'll now work with quasi-isogenies so that we can make sense of
things like $(\pi')^{-1}$.  If we work on points the previous
relation can be written
$$
\phi^\sigma(P)=\pi''\circ\phi\circ(\pi')^{-1}(P).
$$
Therefore
$$
\lambda^{-1}(P)=\xi(\sigma)(P)=\pi''\circ\phi\circ(\pi')^{-1}\circ\phi^{-1}(P).
$$
Now
$$
\pi''(P)=[-p]\circ\lambda^{-1}(P)=
[-p]\circ\pi''\circ\phi\circ(\pi')^{-1}\circ\phi^{-1}(P)\quad\Rightarrow\quad 
\pi'(P)=[-p](P),
$$
where we've made extensive use of the fact that $[-p]$ commutes with
everything.  We conclude that $\pi'=[-p]$ as maps $E'(\fpbar)\to
E'(\fpbar)$, from which it follows that they are equal as morphisms.

It remains to prove the functoriality.  Suppose we have a morphism
$f:E_1\to E_2$ of supersingular elliptic curves over $\fpbar$.  Let
$E_1'$, $E_2'$ be the respective canonical $\fptwo$-structures, then
$f$ induces a morphism $f':E_1'\otimes\fpbar\to E_2'\otimes\fpbar$.
On points, we have as before
$$
(f')^\sigma(P)=\pi_2'\circ f'\circ(\pi_1')^{-1}(P)=[-p]\circ f'\circ
[-p]^{-1}(P)=f'(P).
$$
This equality holds on points, so it also holds as morphisms.  Since
$f'$ is fixed by the generator of $G_{\fptwo}$, it is fixed by the
whole group therefore it is defined over $\fptwo$.
\end{proof}

Since $E$ has a canonical $\fptwo$-structure, so does its
cotangent space $\omega(E)$.  Therefore 
the vector space
$\omega^{p^2-1}(E)$ has a canonical $\fptwo$-basis and we can
identify $\omega^k(E)$ with $\omega^{k+p^2-1}(E)$.

Let $\Sigma(N)$ denote the finite set of
$\fpbar$-isomorphism classes of triples $(E,\alpha,\eta)$,
where $E$ is a supersingular elliptic curve over $\fpbar$,
$\alpha:E[N]\to(\zed/N\zed)^2$ is a level $N$ structure and $\eta\neq
0$ is an invariant differential defined over $\fptwo$.  To save
brackets, we'll write the isomorphism class of the triple
$(E,\alpha,\eta)$ as $[E,\alpha,\eta]$.

On the other hand, fix a supersingular elliptic curve $E_0$ over
$\fpbar$ and let $\oh=\End(E_0)$, $B=\End^0(E_0)=\oh\otimes\cue$.  For
any prime $\ell\neq p$ let $\oh_{\ell}^\times(N)$ denote those elements of
$\oh_{\ell}^\times$ which are congruent to $1$ modulo $\ell^n$, where
$\ell^n\|N$.  Let $\oh_p^\times(1)$ be the kernel of the map
$\oh_p^\times\to\fptwo^\times$ given by reduction modulo the
uniformizer $\pi$ of $\oh_p$.  Finally let 
$$
U:=B_\infty^\times\times\oh_p^\times(1)\times\prod_{\ell\neq p}
\oh_{\ell}^\times(N)
$$
and
$$
\Omega(N):=U\bs\eh_B^\times/B^\times.
$$

For $\ell\nmid pN$ we have Hecke operators $T_{\ell}$ acting on both $\Sigma(N)$ and $\Omega(N)$ (see \S{}\ref{sect:hecke_elliptic} for the definition).  The main technical result of this section is
\begin{thm}\label{thm:bijection}
There exists a bijection $\Sigma(N)\cong\Omega(N)$ which is compatible with 
\begin{itemize}
\item the action of the Hecke operators $T_{\ell}$, $\ell\nmid pN$,
\item the action of $\GL_2(\zed/N\zed)$,
\item raising the level $N$.
\end{itemize}
\end{thm}

We now return to the topic of modular forms.  We can identify $S_k(N)$ with the functions $f:\Sigma(N)\to\fpbar$ such that 
$$
f([E,\alpha,\lambda\eta])=\lambda^{-k}
f([E,\alpha,\eta])\quad\quad\forall \lambda\in\fptwo^\times.
$$
The \emph{modular forms $M_k(\Omega(N))$ of weight $k$ on $\Omega(N)$}\index{modular form!on $\Omega(N)$} are
functions $f:\Omega(N)\to\fpbar$ such that
$$
f(\lambda[x])=\lambda^{-k} f([x])\quad\quad\forall
\lambda\in\oh_p^\times/\oh_p^\times(1)=\fptwo^\times.
$$

The Hecke operators on $\Sigma(N)$, respectively $\Omega(N)$, induce Hecke operators acting on the corresponding modular forms.  It follows immediately from Theorem~\ref{thm:bijection} that there is a bijection 
$$
S_k(N)\xrightarrow{\sim} M_k(\Omega(N))
$$
which is compatible with the action of the Hecke operators.

\begin{lem} \label{lem:sk_periodic}
As $k$ varies, $S_k(N)$ is periodic of period $p^2-1$.
\end{lem}
\begin{proof}
By Lemma~\ref{lem:canonical_supersing} we know that any supersingular
curve $E$ over $\fpbar$ has a canonical $\fptwo$-structure.  Therefore
the cotangent space $\omega(E)$ has a canonical $\fptwo$-structure and
$\omega^{p^2-1}(E)$ has a canonical $\fptwo$-basis.  Hence for any $k$
we can identify $\omega^k(E)$ with $\omega^{k+p^2-1}(E)$, and $S_k(N)$
with $S_{k+p^2-1}(N)$.
\end{proof}

\begin{lem}
Every sequence $(a_{\ell})$ occurring as a system of eigenvalues of the
Hecke operators acting on some $M_k(N)$ also occurs in some $S_{k'}(N)$, and 
conversely.
\end{lem}
\begin{proof}
Suppose $(a_{\ell})$ is given by $f\in M_k(N)$, i.e. $T_{\ell}f=a_{\ell}f$ for all $\ell\nmid pN$.  Factor $f=A^mg$, where $m$ is a nonnegative integer and $A$ does not divide $g$.  Let $k'=k-m(p-1)$, then $g$ has weight $k'$.  Since multiplication by $A$ is a Hecke map we conclude that $g$ is a Hecke eigenform with eigenvalues $a_\ell$.

Conversely, if $(a_{\ell})$ is given by $g\in S_k(N)$ we may assume by Lemma~\ref{lem:sk_periodic} that $k\geq p+1$, in which case the restriction map
$$
M_k(N)\to S_k(N)
$$
 is surjective.  We now use Proposition 1.2.2 of~\cite{ash1} to conclude that $(a_{\ell})$ is given by some $f\in M_k(N)$.
\end{proof}

Let $k$ vary.  We conclude that there is a bijection between the
systems of eigenvalues coming from 
$$
M_*(N):=\bigoplus_k M_k(N)
$$
and those coming from
$$
M_*(\Omega(N)):=\bigoplus_{k\text{ mod }(p^2-1)} M_k(\Omega(N)).
$$

If we let $N$ vary, we get
$$
\bigoplus_N M_*(N)=\left\{\text{functions }\bigcup_N\Omega(N)\to\fpbar\right\}.
$$

The right hand side is the space of $\oh_p^\times(1)$-invariant locally constant functions from $\eh_B^\times/B^\times$ to $\fpbar$.  The following lemma is the last step of the proof of Theorem~\ref{thm:serre}:
\begin{lem}
Let $G$ be a pro-$p$-group with a continuous action on an
$\fpbar$-vector space $V$, and let $T_{\ell}$ be endomorphisms of $V$ that
commute with the action of $G$.  Let $(a_{\ell})$ be a system of
eigenvalues of the $T_{\ell}$ given by a common eigenvector $v\in V$.  Then
we can choose $v$ in such a way that it is $G$-invariant.
\end{lem}
\begin{proof}
We are given a continuous representation
$$
\rho:G\to\GL(V),
$$
where $G=\invlim G_j$ has the profinite topology and $\GL(V)$ has the
discrete topology.  The kernel of $\rho$ is the inverse image of the
open set $\{1\}$ so it is open, but a basis of open neighborhoods of
the identity in $G$ is given by the kernels of the maps $G\to G_j$.
Therefore $\rho$ factors as in the following diagram:
$$
\UseTips
\xymatrix{
G \ar[rr]^{\rho} \ar[dr]_{\pi} & & \GL(V)\\
& \bar{G}, \ar[ur]_{\bar{\rho}}
}
$$
where $\bar{G}$ is a finite group.  Define
$$
w:=\sum_{g\in\bar{G}} \bar{\rho}(g)v.
$$
We have
$$
T_{\ell}w=\sum_{g\in\bar{G}} T_{\ell}\bar{\rho}(g)v=\sum_{g\in\bar{G}}
\bar{\rho}(g)T_{\ell}v=\sum_{g\in\bar{G}}
\bar{\rho}(g)a_{\ell}v=a_{\ell}\sum_{g\in\bar{G}} \bar{\rho}(g)v=a_{\ell}w,
$$
so $w$ is a common eigenvector of the $T_{\ell}$, with the same eigenvalues
as $v$.  Finally for $h\in G$ we have
$$
\rho(h)w=\bar{\rho}(\pi(h))w=
\sum_{g\in\bar{G}}\bar{\rho}(\pi(h))\bar{\rho}(g)v=
\sum_{g\in\bar{G}}\bar{\rho}(\pi(h)g)v=
\sum_{g'\in\bar{G}}\bar{\rho}(g')v=w,
$$
so $w$ is $G$-invariant.
\end{proof}

\subsection{Preliminary results}

\subsubsection{$p$-divisible groups and Dieudonn\'e modules}\label{subsect:dieudonne_elliptic}
We start by recalling the basic terminology and results of
contravariant Dieudonn\'e theory, following~\cite{fontaine1}.

A \emph{$p$-divisible group of height $h$}\index{p-divisible group@$p$-divisible group} is a system
$$
G_0\xrightarrow{i_0} G_1\xrightarrow{i_1}\ldots\to
G_n\xrightarrow{i_n}\ldots,
$$
where for all $n$, $G_n$ is a finite commutative group scheme of rank
$p^{hn}$, $i_n$ is a group homomorphism and the following sequence is
exact:
$$
0\to G_n\xrightarrow{i_n} G_{n+1}\xrightarrow{p^n} G_{n+1}.
$$
This definition was tailored specifically so that
$A[p^\infty]=(A[p^n])$ is a $p$-divisible group of height $2\dim A$.

If $G=(G_n,i_n)$ and $H=(H_n,j_n)$ are $p$-divisible groups, a
homomorphism $f:G\to H$ is a system of group scheme homomorphisms
$f_n:G_n\to H_n$ such that the following diagram commutes for all
$n\geq 1$:
$$
\UseTips
\xymatrix{
G_n\ar[r]^{i_n}\ar[d]_{f_n} & G_{n+1}\ar[d]^{f_{n+1}}\\
H_n\ar[r]^{j_n} & H_{n+1}.
}
$$

Let $k$ be a perfect field of characteristic $p>0$.  Let $W:=W(k)$ be
the ring of infinite Witt vectors over $k$, i.e. the ring of integers
of the absolutely unramified complete extension $L$ of $\cue_p$ with
residue field $k$.  Let $\sigma:W\to W$ be the automorphism induced by
the Frobenius $x\mapsto x^p$ of $k$.  Let $\mathcal{A}:=W[F,V]$, where
the Frobenius $F$ and the Verschiebung $V$ are variables satisfying
$$
FV=VF=p,\quad F\lambda=\sigma(\lambda)F,\quad
V\lambda=\sigma^{-1}(\lambda) V\quad\forall \lambda\in W.
$$
Let $\mathbb{W}_n$ denote the $n$-th Witt group scheme.  If $G$ is a
commutative $k$-group scheme, we define its \emph{Dieudonn\'e module}\index{Dieudonn\'e module} by
$$
M(G):=\dirlim \Hom_{k\text{-gp}}(G,\mathbb{W}_n)\oplus
\left(W(\bar{k})\otimes_\zed
\Hom_{\bar{k}\text{-gp}}(G_{\bar{k}},(\mathbb{G}_m)_{\bar{k}})
\right)^{\Gal(\bar{k}/k)}.
$$

\begin{thm}[Dieudonn\'e-Cartier-Barsotti-Oda]
The functor $M$ gives an anti-equi\-valence from the category of finite commutative
$k$-group schemes of $p$-power rank to the category of left
$\mathcal{A}$-modules of finite $W$-length, taking a group of rank
$p^n$ to a module of length $n$.  It is compatible with perfect base
extension, i.e. if $K/k$ is perfect then
$$
M(G_K)\cong W(K)\otimes_W M(G).
$$
\end{thm}
The proof of this is on page 69 of \cite{demazure1}.  The main
ingredient is to be found in \S{}V.1.4 of~\cite{demazure2}.

\begin{thm}[Oda]
If $G=(G_n)$ is a $p$-divisible group of height $h$, $M(G):=\invlim
M(G_n)$ is a left $\mathcal{A}$-module which is $W$-free of rank $h$.
This gives an equivalence of these two categories which is compatible
with perfect base extension.
\end{thm}

One can define the \emph{dual}\index{Dieudonn\'e module!dual} of a Dieudonn\'e module $M$ by setting
$M^*:=\Hom_W(M,W)$ and
$$
(F^*u)(m):=\sigma(u(Vm)),\quad (V^*u)(m):=\sigma^{-1}(u(Fm))
$$
for all $u\in M^*$, $m\in M$.
\begin{lem}\label{lem:finite-W-length}
If $W$ is a local ring with residue field $k$, then any $W$-module of
length $1$ is isomorphic to $k$.  If $W$ is a DVR with uniformizer
$\pi$, then any $W$-module of length $n$ is isomorphic to one of the
form 
\[\bigoplus_{i=1}^l W/(\pi^{e_i}W),\quad e_i>0,\sum e_i=n.\]
\end{lem}
\begin{proof}
Let $M$ be a $W$-module of length $1$.  Let $x\in M$ nonzero and consider $Wx$.  This is a $W$-module and $0\subsetneq Wx\subset M$ so because $M$ has length $1$ we have $Wx=M$.  From this we conclude that the $W$-module homomorphism $W\to M$ given by $w\to wx$ is surjective with kernel $\Ann M=\{w\in W:wM=0\}$, so $M\cong W/\Ann M$.  Since $\pi W$ is the unique maximal ideal of $W$ we have a canonical surjective homomorphism of $W$-modules
$$
\phi:W/\Ann M\to W/\pi W\cong k.
$$  
But $M$ has length $1$ so $\ker\phi$ is either $0$ or $M$.  The latter is impossible since $\phi$ is surjective, so $M\cong k$.

We now prove the second assertion by induction on the length $n$.  The
base case $n=1$ is what we just proved above.  Now suppose the
assertion is true for all $W$-modules of length strictly less than $n$
and let $M$ be a $W$-module of length $n$.  We have two possibilities:
\begin{itemize}
\item Suppose $M=Wx$.  Then as before $M\cong W/\Ann M$.  But $\Ann M$ is
a proper ideal of $W$, hence $\Ann M=\pi^i W$ for some $i$, since $W$
is a DVR.  The quotient $W/\pi^i W$ has length $i$, on the other hand $M$ has
length $n$ so $i=n$ and $M\cong W/\pi^n W$.
\item Suppose $M$ is not generated by a single element; say
$\{x_1,\ldots,x_j\}$ is a minimal set of generators for $M$.  I claim
that $M$ is actually the direct sum of the $Wx_i$.  Suppose $Wx_1\cap
Wx_2\neq 0$.  Then there exist nonzero $w_1,w_2\in W$ such that
$w_1x_1=w_2x_2$.  But $W$ is a DVR, so $w_1=u_1\pi^{n_1}$,
$w_2=u_2\pi^{n_2}$ where $u_1,u_2\in W^\times$ and $n_1,n_2\in\en$.
Suppose without loss of generality that $n_2\geq n_1$.  Then
$x_1=u_1^{-1}u_2\pi^{n_2-n_1}x_2$, so $Wx_1\subset Wx_2$,
contradicting the minimality of the chosen set of generators.  So
$$
M=Wx_1+Wx_2+\ldots Wx_j
$$
and we can apply the induction hypothesis to each of the $Wx_i$ since
they have length smaller than $n$.
\end{itemize}
\end{proof}

Let $E'$ be a supersingular elliptic curve over $\fptwo$, whose
Frobenius satisfies $\pi_{E'}=-p$.  The goal of this section is to
describe the structure of the Dieudonn\'e module of the
$p$-divisible group $E'[p^\infty]$, and to identify its endomorphism
ring.  The idea is to understand $M(E'[p])$ first, and then lift it to
get $M(E'[p^2])$, and so on all the way to $M(E'[p^\infty])$.

Within this section we'll use $k$ to denote the field $\fptwo$.  If $\tau$ is an automorphism of $k$, $T:k\to k$ is a $\tau$-linear operator and $A$ is the matrix of $T$ with respect to some basis, we shall write
$$
T=A\tau.
$$
This helps eliminate the confusion that arises when composing several
such operators.

\paragraph{The structure of $M(E'[p])$}
We're looking for a left $\mathcal{A}$-module of $W$-length $2$ which
is killed by $p$ and such that 
$$F^2=-p\sigma^2=0.$$
By Lemma~\ref{lem:finite-W-length}, we have that $M\cong
(W/pW)^2=k^2$.  Alternatively, we know that
$$
M(E'[p^n])=M(E'[p^\infty])/(p^n)
$$
and $M(E'[p^\infty])$ is free of rank $2$ because $E'[p^\infty]$ has
height $2$.  Therefore $M(E'[p^n])\cong W_n^2$.

It remains to find the semi-linear operators $F$ and
$V$.  A priori there are several cases according to the $k$-dimension
of $\ker F$.  But the following result, applied with $X=Y=E'$,
$f=f^t=$ multiplication by $p$ says that $E'[p]$ is self-dual:
\begin{thm}[See Theorem III.19.1 in~\cite{oort2}]
Let $f:X\to Y$ be an isogeny of abelian varieties with kernel $K$, and
let $\tilde{K}$ be the kernel of the dual isogeny $f^t:Y^t\to X^t$.
Then there is a canonical isomorphism between the dual $K^t$ of $K$
and $\tilde{K}$.
\end{thm}
Therefore $M=M(E'[p])$ is self-dual, which in particular means that
the kernels of $F$ and $V$ have the same dimension.  First we rule out
the case $\dim\ker F=\dim\ker V=0$: $F$ cannot be bijective because
$F^2=0$.

We also rule out the case $\dim\ker F=\dim\ker V=0$, i.e. $F=V=0$.
This can be done in two different ways:
\begin{itemize}
\item If $F=0$, then the Frobenius of $M(E'[p^2])$ can be expressed as
  $pA$ for some $A\in M_2(k)$, and its square will be $p^2A^2=0\neq
  -p$, contradiction.
\item If $F=V=0$, then $M=M(\alpha_p)\oplus M(\alpha_p)$, where
  $M(\alpha_p)=(k,F=0,V=0)$ is the Dieudonn\'e module of the finite
  group scheme $\alpha_p$.  So we want to show that $E'[p]$ is not
  isomorphic to $\alpha_p^2$.  The dimension of the tangent space to
  $E'[p]$ at the origin is the same as the dimension of the tangent
  space to $E'$ at the origin, which is $1$ since $E'$ is a smooth
  curve.  The dimension of the tangent space to $\alpha_p^2$ at the
  origin is twice the dimension of the tangent space to $\alpha_p$ at
  the origin, i.e. $2\times 1=2$.  Therefore we get a contradiction.
\end{itemize}

The only remaining possibility is $\dim\ker F=1$.  Say $\ker F=km_1$
and pick $m_2\in M$ such that $\{m_1, m_2\}$ is a basis of $M$, then
$F=\left(\begin{smallmatrix}0&a_1\\0&a_2\end{smallmatrix}\right)\sigma$.  We impose
$F^2=0$:
$$
F^2=\left(\begin{smallmatrix}0&a_1a_2^p\\0&a_2^{p+1}\end{smallmatrix}\right)\sigma^2=0\Rightarrow
a_2=0.
$$
We can also change the basis in such a way that
$a_1$ becomes $1$: let $m_1'=m_1$, $m_2'=a_1^{-1/p}m_2$ (can take
$p$-th roots because $\eff_{p^2}$ is perfect).  With respect to the
basis $\{m_1', m_2'\}$ we have
$$
F=\left(\begin{smallmatrix}0&1\\0&0\end{smallmatrix}\right)\sigma.
$$
The Verschiebung is given a priori by some matrix
$$
V=\left(\begin{smallmatrix}b_{11}&b_{12}\\b_{21}&b_{22}\end{smallmatrix}\right)\sigma^{-1}.
$$
We impose $FV=VF=0$:
\begin{eqnarray*}
FV&=&\left(\begin{smallmatrix}b_{21}^p&b_{22}^p\\0&0\end{smallmatrix}\right)=0\Rightarrow
b_{21}=b_{22}=0,\\
VF&=&\left(\begin{smallmatrix}0&b_{11}\\0&0\end{smallmatrix}\right)=0\Rightarrow b_{11}=0.
\end{eqnarray*}
Therefore the matrix of $V$ is of the form
$$
V=\left(\begin{smallmatrix}0&\lambda\\0&0\end{smallmatrix}\right)\sigma^{-1},
$$
where $\lambda\neq 0$ since $\dim\ker V=1$.  Any change of basis that
fixes the matrix of $F$ also fixes the matrix of $V$, so it looks like
we're stuck with the parameter $\lambda$.  But we can find the actual
value of $\lambda$ by investigating $M_2=M(E'[p^2])$.

\paragraph{The structure of $M(E'[p^n])$}
We use Nakayama's lemma to lift the basis of $k^2$ to a basis of
$W_2^2$.  With respect to this basis the operators $F$ and $V$ look
like
$$
F=\left(\begin{smallmatrix}0&1\\0&0\end{smallmatrix}\right)\sigma+pA\sigma,\quad
V=\left(\begin{smallmatrix}0&\lambda\\0&0\end{smallmatrix}\right)\sigma^{-1}+pB\sigma^{-1}.
$$
We impose the conditions $F^2=-p\sigma^2$, $FV=VF=p$:
\begin{eqnarray*}
F^2&=&p\left(\begin{smallmatrix}a_{21}^p&a_{11}+a_{22}^p\\0&a_{21}\end{smallmatrix}\right)\sigma^2
=-p\sigma^2 \Rightarrow a_{21}=-1, a_{11}=-a_{22}^p,\\
FV&=&p\left(\begin{smallmatrix}b_{21}^p&b_{22}^p-\lambda^pa_{22}^p\\0&-\lambda^p\end{smallmatrix}\right)=p
\Rightarrow b_{21}=1,\lambda=-1,a_{22}=-b_{22},\\
VF&=&p\left(\begin{smallmatrix}1&b_{22}^{1/p}+b_{11}\\0&1\end{smallmatrix}\right)=p
\Rightarrow b_{22}=-b_{11}^p.
\end{eqnarray*}
This implies that $F$ and $V$ are of the form
$$
F=\left(\begin{smallmatrix}0&1\\0&0\end{smallmatrix}\right)\sigma+
p\left(\begin{smallmatrix}-b_{11}&a_{12}\\-1&b_{11}^p\end{smallmatrix}\right)\sigma,\quad 
V=\left(\begin{smallmatrix}0&-1\\0&0\end{smallmatrix}\right)\sigma^{-1}+
p\left(\begin{smallmatrix}b_{11}&b_{12}\\1&-b_{11}^p\end{smallmatrix}\right)\sigma^{-1}.
$$
In particular,
$$
M(E'[p])=\left(k^2, F=\left(\begin{smallmatrix}0&1\\0&0\end{smallmatrix}\right)\sigma,
V=\left(\begin{smallmatrix}0&-1\\0&0\end{smallmatrix}\right)\sigma^{-1}\right).
$$
We're only a few matrix multiplications away from having the structure
of $M(E'[p^2])$, so we might as well do it.  We want to change the
basis of $W_2^2$ in such a way that the matrices for $F$ and $V$ look
nicer.  To keep this simple, we will change the basis by a matrix of
the form
$$
I+pC,\quad\text{ where }
C=\left(\begin{smallmatrix}c_{11}&c_{12}\\c_{21}&c_{22}\end{smallmatrix}\right)\in M_2(k).
$$

The new matrix for the Frobenius will be
\begin{multline*}
\left(I-p\left(\begin{smallmatrix}c_{11}&c_{12}\\c_{21}&c_{22}\end{smallmatrix}\right)\right)
\left(\left(\begin{smallmatrix}0&1\\0&0\end{smallmatrix}\right)+
p\left(\begin{smallmatrix}-b_{11}&a_{12}\\-1&b_{11}^p\end{smallmatrix}\right)\right)
\sigma\left(I+p\left(\begin{smallmatrix}c_{11}&c_{12}\\c_{21}&c_{22}\end{smallmatrix}\right)\right)=\\
=\left(\begin{smallmatrix}0&1\\0&0\end{smallmatrix}\right)\sigma+p\left(\begin{smallmatrix}-b_{11}+c_{21}^p&a_{12}-c_{11}+c_{22}^p\\-1&b_{11}^p-c_{21}\end{smallmatrix}\right)\sigma.
\end{multline*}
So we can set $c_{21}=b_{11}^p$, $c_{11}=c_{22}^p+a_{12}$ and get a
very simple expression for $F$.  What about $V$?  The new matrix for
it is
\begin{multline*}
\left(I-p\left(\begin{smallmatrix}a_{12}+c_{22}^p&c_{12}\\b_{11}^p&c_{22}\end{smallmatrix}\right)\right)
\left(\left(\begin{smallmatrix}0&-1\\0&0\end{smallmatrix}\right)+
p\left(\begin{smallmatrix}b_{11}&b_{12}\\1&-b_{11}^p\end{smallmatrix}\right)\right)\sigma^{-1}
\left(I+p\left(\begin{smallmatrix}a_{12}+c_{22}^p&c_{12}\\b_{11}^p&c_{22}\end{smallmatrix}\right)\right)=\\
=\left(\begin{smallmatrix}0&-1\\0&0\end{smallmatrix}\right)\sigma^{-1}+
p\left(\begin{smallmatrix}0&a_{12}+b_{12}\\1&0\end{smallmatrix}\right)\sigma^{-1}.
\end{multline*}
All we can say so far is that there exists a basis of $W_2^2$ with
respect to which
\begin{eqnarray*}
F&=&\left(\begin{smallmatrix}0&1\\0&0\end{smallmatrix}\right)\sigma+
  p\left(\begin{smallmatrix}0&0\\-1&0\end{smallmatrix}\right)\sigma,\\
V&=&\left(\begin{smallmatrix}0&-1\\0&0\end{smallmatrix}\right)\sigma^{-1}+
  p\left(\begin{smallmatrix}0&\lambda\\1&0\end{smallmatrix}\right)\sigma^{-1},
\end{eqnarray*}
where $\lambda\in\fptwo$.  It turns out that we can pin down
$\lambda$, once again by going a step higher.  Lift this basis to one
of $W_3^2$; we have
\begin{eqnarray*}
F&=&\left(\begin{smallmatrix}0&1\\0&0\end{smallmatrix}\right)\sigma+
  p\left(\begin{smallmatrix}0&0\\-1&0\end{smallmatrix}\right)\sigma+p^2A\sigma,\\ 
V&=&\left(\begin{smallmatrix}0&-1\\0&0\end{smallmatrix}\right)\sigma^{-1}+
  p\left(\begin{smallmatrix}0&\lambda\\1&0\end{smallmatrix}\right)\sigma^{-1}+p^2B\sigma^{-1}, 
\end{eqnarray*}
We impose the conditions $F^2=-p\sigma^2$, $FV=p$:
\begin{eqnarray*}
F^2&=&p\left(\begin{smallmatrix}-1&0\\0&-1\end{smallmatrix}\right)\sigma^2+
p^2\left(\begin{smallmatrix}a_{21}^p&a_{22}^p+a_{11}\\0&a_{21}\end{smallmatrix}\right)\sigma^2=
-p\sigma^2\Rightarrow a_{21}=0,a_{11}=-a_{22}^p,\\
FV&=&p\left(\begin{smallmatrix}1&0\\0&1\end{smallmatrix}\right)+
p^2\left(\begin{smallmatrix}b_{21}^p&b_{22}^p+a_{22}^p\\0&\lambda^p\end{smallmatrix}\right)=p
\Rightarrow b_{21}=0,a_{22}=-b_{22},\lambda=0.
\end{eqnarray*}
Therefore
$$
M(E'[p^2])=\left(W_2^2,F=\left(\begin{smallmatrix}0&1\\-p&0\end{smallmatrix}\right)\sigma,
V=\left(\begin{smallmatrix}0&-1\\p&0\end{smallmatrix}\right)\sigma^{-1}\right).
$$

\begin{prop}\label{prop:struct_p^n}
For any integer $n\geq 2$ we have
$$
M(E'[p^n])=\left(W_n^2,F=\left(\begin{smallmatrix}0&1\\-p&0\end{smallmatrix}\right)\sigma,
V=\left(\begin{smallmatrix}0&-1\\p&0\end{smallmatrix}\right)\sigma^{-1}\right).
$$
\end{prop}
\begin{proof}
We proceed by induction on $n$, the base case $n=2$ having been
done above.  Suppose the statement is true for $n\geq 2$.  Use
Nakayama's lemma to lift the basis of $W_n^2$ to one of $W_{n+1}^2$.
With respect to the latter we have
\begin{eqnarray*}
F&=&\left(\begin{smallmatrix}0&1\\0&0\end{smallmatrix}\right)\sigma+
  p\left(\begin{smallmatrix}0&0\\-1&0\end{smallmatrix}\right)\sigma+p^nA\sigma,\\ 
V&=&\left(\begin{smallmatrix}0&-1\\0&0\end{smallmatrix}\right)\sigma^{-1}+
  p\left(\begin{smallmatrix}0&0\\1&0\end{smallmatrix}\right)\sigma^{-1}+p^nB\sigma^{-1}. 
\end{eqnarray*}
We impose the conditions $F^2=-p$, $FV=VF=p$:
\begin{eqnarray*}
F^2&=&-p\sigma^2+p^n\left(\begin{smallmatrix}a_{21}^p&a_{22}^p+a_{11}\\0&a_{21}\end{smallmatrix}\right)\sigma^2
=-p\sigma^2\Rightarrow a_{21}=0,a_{11}=-a_{22}^p,\\
FV&=&p+p^n\left(\begin{smallmatrix}b_{21}^p&b_{22}^p+a_{22}^p\\0&0\end{smallmatrix}\right)=p
\Rightarrow a_{22}=-b_{22},a_{11}=b_{22}^p,b_{21}=0,\\
VF&=&p+p^n\left(\begin{smallmatrix}0&b_{22}^{1/p}+b_{11}\\0&0\end{smallmatrix}\right)=p
\Rightarrow b_{22}=-b_{11}^p,
\end{eqnarray*}
and conclude that
$$
A=\left(\begin{smallmatrix}-b_{11}&a_{12}\\0&b_{11}^p\end{smallmatrix}\right),\quad
B=\left(\begin{smallmatrix}b_{11}&b_{12}\\0&-b_{11}^p\end{smallmatrix}\right).
$$

We'll apply a change of basis with matrix $I+p^nC$.  The new matrix of
the Frobenius is
\begin{multline*}
\left(I-p^n\left(\begin{smallmatrix}c_{11}&c_{12}\\c_{21}&c_{22}\end{smallmatrix}\right)\right)
\left(\left(\begin{smallmatrix}0&1\\0&0\end{smallmatrix}\right)+p\left(\begin{smallmatrix}0&0\\-1&0\end{smallmatrix}\right)+
  p^n\left(\begin{smallmatrix}-b_{11}&a_{12}\\0&b_{11}^p\end{smallmatrix}\right)\right)\sigma
\left(I+p^n\left(\begin{smallmatrix}c_{11}&c_{12}\\c_{21}&c_{22}\end{smallmatrix}\right)\right)=\\
=\left(\begin{smallmatrix}0&1\\0&0\end{smallmatrix}\right)\sigma+
p\left(\begin{smallmatrix}0&0\\-1&0\end{smallmatrix}\right)\sigma+
p^n\left(\begin{smallmatrix}-b_{11}+c_{21}^p&a_{12}-c_{11}+c_{22}^p\\0&b_{11}^p-c_{21}\end{smallmatrix}\right)\sigma.
\end{multline*}
So if we let $c_{11}=a_{12}+c_{22}^p$ and $c_{21}=b_{11}^p$ then $F$
is of the desired form.

The new matrix of the Verschiebung is
\begin{multline*}
\left(I-p^n\left(\begin{smallmatrix}a_{12}+c_{22}^p&c_{12}\\b_{11}^p&c_{22}\end{smallmatrix}\right)\right)
\left(\left(\begin{smallmatrix}0&-1\\0&0\end{smallmatrix}\right)+p\left(\begin{smallmatrix}0&0\\1&0\end{smallmatrix}\right)+
p^n\left(\begin{smallmatrix}b_{11}&b_{12}\\0&-b_{11}^p\end{smallmatrix}\right)\right)\sigma^{-1}
\left(I+p^n\left(\begin{smallmatrix}a_{12}+c_{22}^p&c_{12}\\b_{11}^p&c_{22}\end{smallmatrix}\right)\right)=\\
=\left(\begin{smallmatrix}0&-1\\0&0\end{smallmatrix}\right)\sigma^{-1}+
p\left(\begin{smallmatrix}0&0\\1&0\end{smallmatrix}\right)\sigma^{-1}+
p^n\left(\begin{smallmatrix}0&\lambda\\0&0\end{smallmatrix}\right)\sigma^{-1},
\end{multline*}
where $\lambda\in\fptwo$.  We want to show that $\lambda=0$; for this
we need to go to $M(E'[p^{n+2}])$.  Lift the new basis of $W_{n+1}^2$
to $W_{n+2}^2$.  Our operators are of the form
\begin{eqnarray*}
F&=&\left(\begin{smallmatrix}0&1\\0&0\end{smallmatrix}\right)\sigma+
  p\left(\begin{smallmatrix}0&0\\-1&0\end{smallmatrix}\right)\sigma+ p^{n+1}A\sigma,\\
V&=&\left(\begin{smallmatrix}0&-1\\0&0\end{smallmatrix}\right)\sigma^{-1}+
  p\left(\begin{smallmatrix}0&0\\1&0\end{smallmatrix}\right)\sigma^{-1}+
  p^n\left(\begin{smallmatrix}0&\lambda\\0&0\end{smallmatrix}\right)\sigma^{-1}+p^{n+1}B\sigma^{-1}.
\end{eqnarray*}
We impose the conditions $F^2=-p\sigma^2$ and $FV=p$:
\begin{eqnarray*}
F^2&=&-p\sigma^2+
p^{n+1}\left(\begin{smallmatrix}a_{21}^p&a_{22}^p+a_{11}\\0&a_{21}\end{smallmatrix}\right)\sigma^2=
-p\sigma^2\Rightarrow a_{21}=0,a_{11}=-a_{22}^p,\\
FV&=&p+p^{n+1}\left(\begin{smallmatrix}b_{21}^p&b_{22}^p+a_{22}^p\\0&-\lambda^p\end{smallmatrix}\right)=
  p\Rightarrow b_{21}=0,a_{22}=-b_{22},\lambda=0.
\end{eqnarray*}
This concludes the proof of the proposition.
\end{proof}

A direct consequence of this is that
$$
M(E'[p^\infty])=\left(W^2,
F=\left(\begin{smallmatrix}0&1\\-p&0\end{smallmatrix}\right)\sigma,
V=\left(\begin{smallmatrix}0&-1\\p&0\end{smallmatrix}\right)\sigma^{-1}\right).
$$

\paragraph{The endomorphism ring of $M(E'[p^\infty])$}
\begin{cor}\label{cor:isom_p}
Let $M=M(E'[p^\infty])$.  Then $\End(M)=\oh_p=\oh\otimes\zed_p$, where
$\oh=\End(E')$.  Moreover, $\oh_p^\times(1)$ can be identified with
the group of automorphisms of $M$ which lift the identity map on
$M/FM$.
\end{cor}
\begin{proof}
Let $g\in\End(M)$; it is a $W$-linear map that commutes with $F$ and
$V$.  Suppose $g$ is given by a matrix $(g_{ij})\in M_2(W)$.  We have
\begin{eqnarray*}
F\circ g &=& \left(\begin{smallmatrix}0&1\\-p&0\end{smallmatrix}\right)\sigma
\left(\begin{smallmatrix}g_{11}&g_{12}\\g_{21}&g_{22}\end{smallmatrix}\right)= 
\left(\begin{smallmatrix}g_{21}^p&g_{22}^p\\-pg_{11}^p&-pg_{12}^p\end{smallmatrix}\right)\sigma,\\
g\circ F &=& \left(\begin{smallmatrix}g_{11}&g_{12}\\g_{21}&g_{22}\end{smallmatrix}\right)
\left(\begin{smallmatrix}0&1\\-p&0\end{smallmatrix}\right)\sigma=
\left(\begin{smallmatrix}-pg_{12}&g_{11}\\-pg_{22}&g_{21}\end{smallmatrix}\right)\sigma.
\end{eqnarray*}
These should be equal so we get $g_{21}^p=-pg_{12}$,
$g_{11}=g_{22}^p$.  We also impose the condition $V\circ g=g\circ V$,
but this doesn't give anything new.  Therefore
$$
\End(M)=\left\{\left(\begin{smallmatrix}x&y\\-py^p&x^p\end{smallmatrix}\right):x,y\in
W(\fptwo)\right\}=\left\{\left(\begin{smallmatrix}x&0\\0&x^p\end{smallmatrix}\right)+
F\left(\begin{smallmatrix}y&0\\0&y^p\end{smallmatrix}\right):x,y\in W(\fptwo)\right\}.
$$
But $W(\fptwo)$ is the ring of integers of the unique unramified
quadratic extension $L$ of $\cue_p$.  Let $\pi$ be a solution of
$X^2+p=0$ in $\bar{L}$.  The map $\sigma:x\mapsto x^p$ is the unique
nontrivial automorphism of $L$.  It is now easy to see that the map
\begin{eqnarray*}
\varphi:\End(M) &\to& B_p=\{L,-p\}=B\otimes\cue_p\\
\left(\begin{smallmatrix}x&0\\0&x^p\end{smallmatrix}\right)+
  F\left(\begin{smallmatrix}y&0\\0&y^p\end{smallmatrix}\right)&\mapsto& x+\pi y 
\end{eqnarray*}
is an injective ring homomorphism.  It identifies $\End(M)$ with
$\oh_p=\{x+\pi y:x,y\in\oh_L\}$, the unique maximal order of $B_p$.

It remains to prove the last statement.  Let
$g=\left(\begin{smallmatrix}x&y\\-py^p&x^p\end{smallmatrix}\right)\in\End(M)^\times=\oh_p^\times$.
Note that 
$$
M/FM=\left\{\left(\begin{smallmatrix}0\\a\end{smallmatrix}\right)+FM:a\in\fptwo\right\}.
$$
Let $\bar{x}$ be the reduction of $x$ modulo $\pi$, then $g$ restricts to multiplication by $\bar{x}^p$ on $M/FM$.

Therefore $g$ restricts to the identity if and only if $\bar{x}=1$, which means that the group of such automorphisms is identified with the kernel of the reduction modulo $\pi$, i.e. with $\oh_p^\times(1)$.
\end{proof}

\begin{cor}\label{cor:struct_p}
Let $E_1$, $E_2$ be two supersingular elliptic curves over $\fpbar$
and let $E_1'$, $E_2'$ denote their canonical $\fptwo$-structures.
Put $M_1=M(E_1'[p^\infty])$, $M_2=M(E_2'[p^\infty])$.  Then $M_1\cong
M_2$ and any isomorphism $M_1/FM_1\cong M_2/FM_2$ lifts to an
isomorphism $M_1\cong M_2$.
\end{cor}
\begin{proof}
We already proved the first part by describing the structure of
$M=M(E'[p^\infty])$ for such curves.  As for the second part, it
suffices to show that any automorphism of $M/FM$ lifts to an
automorphism of $M$.  From the description of $M/FM$ in the proof of
the previous corollary we know that the automorphisms are given by
multiplication by some $\lambda\in\fptwo^\times$.  But then the matrix
$$
\left(\begin{smallmatrix}\lambda^p&0\\0&\lambda\end{smallmatrix}\right)
$$
represents an automorphism of $M$ which restricts to multiplication by
$\lambda$ on $M/FM$, which is what we wanted to show.
\end{proof}

\subsubsection{Local study of isogenies}
\begin{lem}\label{lem:ker_isom_coker}
Let $\phi:A\to B$ be an isogeny of abelian varieties.  For \emph{any}
prime $\ell$ and any $n\geq\ord_{\ell}\deg\phi$, there exists a canonical
isomorphism between the kernel and the cokernel of the restriction of
$\phi$ to the $\ell^n$-torsion.
\end{lem}
\begin{proof}
First suppose $\ell\neq p$.  Let $\ell^k\|\deg\phi$ and $d:=\ell^{-k}\deg\phi$.
For any $n>0$ let $K_{\ell,n}$ and $C_{\ell,n}$ denote the kernel and the
cokernel of the restriction of $\phi$ to $A[\ell^n]$: 
$$
0\to K_{\ell,n}\to A[\ell^n]\xrightarrow{\phi} B[\ell^n]\to C_{\ell,n}\to 0.
$$
This is an exact sequence of finite groups so
$$
1=\frac{\#(K_{\ell,n})\cdot\#(B[\ell^n])}{\#(A[\ell^n])\cdot\#(C_{\ell,n})},
$$
but $\#(A[\ell^n])=\ell^{2gn}=\#(B[\ell^n])$ (where $g=\dim A=\dim B$) so $\#(K_{\ell,n})=\#(C_{\ell,n})$.

Since $(d,\ell)=1$, multiplication by $d$ is an injective endomorphism of
$A[\ell^n]$.  But by the order counting that we did above the cokernel is
also trivial, so $d$ gives an automorphism of $A[\ell^n]$.  Therefore it
makes sense to define a map $f:B[\ell^n]\to A[\ell^n]$ as the composition:
$$
\UseTips
\xymatrix{
B[\ell^n]\ar[r]^{\phi^\vee}\ar@(dr,dl)[rrr]_f & A[\ell^n]\ar[r]^{\ell^{n-k}} &
A[\ell^n]\ar[r]^{d^{-1}} & A[\ell^n].
}
$$
Note that 
$$
\phi\circ f= \phi\circ(d^{-1}\ell^{n-k}\phi^\vee)=
(\deg\phi)d^{-1}\ell^{n-k}= \ell^n
$$
and $\ell^n:B[\ell^n]\to B[\ell^n]$ is just the zero map, so we actually have
$f:B[\ell^n]\to\ker\phi=K_{\ell,n}$.  Next we see that
$$
f\circ\phi= (d^{-1}\ell^{n-k}\phi^\vee)\circ\phi= d^{-1}\ell^{n-k}(\deg\phi)
=\ell^n,
$$
so $f$ factors as
$$
\UseTips
\xymatrix{
A[\ell^n]\ar[r]^{\phi} & B[\ell^n]\ar[d]_f\ar[r] & C_{\ell,n}\ar[dl]\ar[r] &
0.\\
& K_{\ell,n}
}
$$
To save letters we'll denote this map $C_{\ell,n}\to K_{\ell,n}$ by $f$.  
Now we define $g:K_{\ell,n}\to C_{\ell,n}$ such that $g$ is the inverse of
$f$.  Let $a\in K_{\ell,n}$ and pick $a'\in A$ such that $\ell^na'=a$.  Note
that $\ell^n\phi(a')=\phi(\ell^na')=\phi(a)=0$ since $a\in\ker\phi$.  So 
$\phi(a')\in B[\ell^n]$; let $g(a)$ be the image of $\phi(a')$ in
$C_{\ell,n}$.  This is easily seen to be well-defined: if $a''\in A$ is
such that $\ell^na''=a$, then $\ell^n(a''-a')=a-a=0$ so $a''-a'\in A[\ell^n]$
and we get
$$
\phi(a'')=\phi(a'+(a''-a'))=\phi(a')+\phi(a''-a'),
$$
and $\phi(a''-a')$ gets mapped to $0$ in $C_{\ell,n}$.  It is clear that
$f$ and $g$ are inverses:
$$
f(g(a))= f(\phi(\ell^{-n}a))= d^{-1}\ell^{n-k}\phi^\vee(\phi(\ell^{-n}a))=
d^{-1}\ell^{n-k}(\deg\phi)\ell^{-n}a= a,
$$
and similarly $g(f(b))=b$.

When $\ell=p$ we have an exact sequence of finite commutative group
schemes
$$
0\to K_{p,n}\to A[p^n]\xrightarrow{\phi} B[p^n]\to C_{p,n}\to 0,
$$
where the map $\phi$ is obtained by the usual argument from
$$
\UseTips
\xymatrix{
0\ar[r] & A[p^n]\ar[r] & A\ar[r]^{p^n}\ar[d]_{\phi} & A\ar[r]\ar[d]_{\phi}
& 0\\
0\ar[r] & B[p^n]\ar[r] & B\ar[r]^{p^n} & B\ar[r] & 0.
}
$$
Applying the Dieudonn\'e functor $M$ we get an exact sequence of left
$\mathcal{A}$-modules of finite length
$$
0\to M(C_{p,n})\to M(B[p^n])\xrightarrow{M(\phi)} M(A[p^n])\to
M(K_{p,n})\to 0.
$$
A trivial modification of the proof given above for $\ell\neq p$ will do
the trick.  The advantage of working with the Dieudonn\'e modules
instead of the group schemes is that for the construction of the
inverse we need to work with actual elements, and they are not
available in group-scheme-land.
\end{proof}

\begin{prop}
Let $\phi:A\to B$ be an isogeny of abelian varieties.  For any prime $\ell\neq p$, $\phi$ induces an injective $\zed_{\ell}$-linear map $T_{\ell}\phi:T_{\ell} A\to T_{\ell} B$ whose cokernel is canonically isomorphic to $(\ker\phi)_{\ell}$.\\
When $\ell=p$, $\phi$ induces an injective $W$-linear map $M(\phi):M(B[p^\infty])\to M(A[p^\infty])$ whose cokernel is canonically isomorphic to $M((\ker\phi)_p)$.
\end{prop}
\begin{proof}
We start with $\ell\neq p$.  Let $(a_n)\in T_{\ell} A$, i.e. $a_n\in A[\ell^n]$
and $a_{n-1}=\ell a_n$ for all $n$.  Set $T_{\ell}\phi((a_n)):=(\phi(a_n))$.
We have $\ell^n\phi(a_n)=\phi(\ell^na_n)=\phi(0)=0$ and
$\phi(a_{n-1})=\phi(\ell a_n)=\ell\phi(a_n)$ so indeed $(\phi(a_n))\in T_{\ell}
B$.

Suppose $T_{\ell}\phi((a_n))=0$, i.e. $\phi(a_n)=0$ for all $n$.  But
$\phi$ is an isogeny so its kernel is finite, therefore the sequence
$(a_n)$ has an infinite constant subsequence.  Fix an integer $N>0$,
then there exist $n>m>N$ such that $a_n=a_m$.  But $a_n=\ell^{n-m}a_m$,
so we get $(\ell^{n-m}-1)a_m=0$, therefore $a_m=0$.  In particular
$a_j=0$ for all $j\leq N$.  This works for arbitrarily large $N$ so
$(a_n)=0$ and $T_{\ell}\phi$ is injective.

By Lemma~\ref{lem:ker_isom_coker} we know that if
$n\geq\ord_{\ell}\deg\phi$, the kernel $K_{\ell,n}$ of the restriction of
$\phi$ to $A[\ell^n]$ is isomorphic to the cokernel $C_{\ell,n}$ of the same
map.  But for $n\geq\ord_{\ell}\deg\phi$ we have $K_{\ell,n}=(\ker\phi)_{\ell}$, so
in particular $C_{\ell,n}$ stabilizes.  Therefore for
$n\geq\ord_{\ell}\deg\phi$ we have 
$$
\UseTips
\xymatrix{
{\coker T_{\ell}\phi} = C_{\ell,n} \ar[r]^f_{\sim} & K_{\ell,n} = (\ker\phi)_{\ell},
}
$$
which is what we wanted to show.

Now suppose $\ell=p$.  The map $M(\phi):M(A[p^\infty])\to
M(B[p^\infty])$ is simply
$$
\UseTips
\xymatrix{
\ldots\ar[r]^p & M(B[p^3])\ar[r]^p\ar[d]_{M(\phi)} &
M(B[p^2])\ar[r]^p\ar[d]_{M(\phi)} & M(B[p])\ar[d]_{M(\phi)}\\
\ldots\ar[r]^p & M(A[p^3])\ar[r]^p & M(A[p^2])\ar[r]^p & M(A[p])
}
$$
To show that $M(\phi)$ is injective, first note that $K_{p,n}$ is a
subgroup scheme of $\ker\phi$ for all $n$.  Since $\phi$ is an
isogeny, $\ker\phi$ is a finite group scheme so the $K_{p,n}$
stabilize to $(\ker\phi)_p$.  Applying the Dieudonn\'e functor gives
that the cokernel of $M(\phi)$ is isomorphic to $M((\ker\phi)_p)$.  

Suppose $M(\phi)((a_n))=0$, then $a_n\in M(C_{p,n})$ for all $n$.  But
the $M(C_{p,n})$ stabilize so there exists $n_0$ such that $a_n\in
M(B[p^{n_0}])\cong W_{n_0}^2$ for all $n$.  Since $p^{n_0}=0$ in
$W_{n_0}^2$ we get that $a_n=p^na_0=0$ for all $n\geq n_0$.  Therefore
$(a_n)=0$ and $M(\phi)$ is injective.
\end{proof}

\subsubsection{Differentials defined over $\fptwo$}
\begin{lem}\label{lem:ec_diff}
Let $E$ be a supersingular elliptic curve over $\fpbar$.  Then a
non-zero invariant differential\index{invariant differential!on supersingular elliptic curve} on $E$ defined over $\fptwo$ is
equivalent to a choice of nonzero element of $M/FM$, where
$M=M(E'[p^\infty])$ and $E'$ is the canonical $\fptwo$-structure of
$E$.
\end{lem}
\begin{proof}
Differentials of $E$ defined over $\fptwo$ can easily be identified
with the differentials of $E'$, i.e. with the cotangent space
$\omega(E')$.  Consider the exact sequence
$$
0\to E'[p]\to E'\xrightarrow{p} E'\to 0.
$$
The contravariant functor $\omega$ is exact so we get an exact
sequence
$$
0\to \omega(E')\xrightarrow{p} \omega(E')\to \omega(E'[p])\to 0.
$$
But $E'$ is supersingular so $p$ is purely inseparable, i.e. it
induces the zero map on (co)tangent spaces.  We conclude that
$\omega(E')\cong\omega(E'[p])$.  By Proposition~III.4.3
in~\cite{fontaine1}, $\omega(E'[p])$ is isomorphic to $M/FM$, so
$\omega(E')\cong M/FM$ and a choice of nonzero invariant differential
of $E$ defined over $\fptwo$ amounts to a choice of nonzero element in
$M/FM$.
\end{proof}

\subsection{Construction of the bijection}\label{sect:elliptic_bijection}
Fix a supersingular elliptic curve $E_0$ over $\fpbar$ and let
$\oh=\End(E_0)$, $B=\End^0(E_0)=\oh\otimes\cue$,
$M_0=M(E_0'[p^\infty])$, where $E_0'$ is the canonical
$\fptwo$-structure given by Lemma~\ref{lem:canonical_supersing}.  We
know that $B$ is isomorphic to the unique 
quaternion algebra ramified at $p$ and $\infty$.  Also fix a level $N$
structure $\alpha_0:E_0[N]\xrightarrow{\sim}(\zed/N\zed)^2$ and a
nonzero invariant differential
$\eta_0\in M_0/FM_0$ defined over $\fptwo$.

We consider triples of the form $(E_1,\alpha_1,\eta_1)$, where $E_1$
is a supersingular elliptic curve over $\fpbar$,
$\alpha_1:E_1[N]\xrightarrow{\sim} (\zed/N\zed)^2$ is a full level $N$
structure and $0\neq \eta_1\in M_1/FM_1$ is a nonzero
invariant differential form on $E_1$ defined over $\fptwo$.  

Recall from \S{}\ref{sect:main_results} that we defined $\Sigma(N)$ to be the finite set of $\fpbar$-isomorphism classes
of triples $(E_1,\alpha_1,\eta_1)$ as above.  By Honda-Tate theory
(\cite{tate3}) any two supersingular elliptic curves are isogenous; also any curve isogenous to a supersingular one is supersingular.  Therefore after the choice of $(E_0,\alpha_0,\eta_0)$ the set
$\Sigma(N)$ can be identified with the set $\Sigma^0(N)$ of triples 
$$
\left(\phi_1:E_0\to E_1, \alpha_1:E_1[N]\xrightarrow{\sim}(\zed/N\zed)^2,
0\neq\eta_1\in M_1/FM_1\right),
$$ 
where $\phi_1:E_0\to E_1$ is an isogeny.  Two such triples
$(\phi_1,\alpha_1,\eta_1)$ and $(\phi_2,\alpha_2,\eta_2)$ are considered the
same if there exists an isomorphism $f:E_1\to E_2$ such that
$M(f')(\eta_2)=\eta_1$ and the following diagram commutes: 
\begin{equation}\label{diag:elliptic_iso}
\UseTips
\xymatrix{
E_1[N] \ar[r]^f \ar[d]_{\alpha_1}^{\sim} & E_2[N]
\ar[d]_{\alpha_2}^{\sim} \\
(\zed/N\zed)^2 \ar@{=}[r] & (\zed/N\zed)^2.
}
\end{equation}

For any prime $\ell\neq p$ let $\oh_{\ell}^\times(N)$ denote those elements of
$\oh_{\ell}^\times$ which are congruent to $1$ modulo $\ell^n$, where
$\ell^n\|N$.  Let $\oh_p^\times(1)$ be the kernel of the map
$\oh_p^\times\to\fptwo^\times$ given by reduction modulo the
uniformizer $\pi$ of $\oh_p$.  Finally let 
$$
U=B_\infty^\times\times\oh_p^\times(1)\times\prod_{\ell\neq p}
\oh_{\ell}^\times(N)
$$
and
$$
\Omega(N)=U(1,N)\bs\eh_B^\times/B^\times.
$$

Since the level $N$ is fixed throughout this section, we'll drop the reference to $N$ and write simply $\Sigma$, $\Sigma^0$ and $\Omega$ for our finite sets.
The purpose of the section is to exhibit a canonical bijection
between $\Sigma^0$ and $\Omega$.

Let $[\phi_1,\alpha_1,\eta_1]\in\Sigma^0$ and pick a representative $(\phi_1,\alpha_1,\eta_1)$ of this isomorphism class.  

By~\cite{tate2} we know that $\End(T_{\ell} E_0)=\oh_{\ell}:=\oh\otimes\zed_{\ell}$.  Choose a $\zed_\ell$-linear isomorphism $k_{\ell,1}:T_{\ell} E_0\xrightarrow{\sim} T_{\ell} E_1$ whose restriction gives a commutative diagram
$$
\UseTips
\xymatrix{
E_0[\ell^n] \ar^{k_{\ell,1}}_{\sim}[r] \ar_{\alpha_0}^{\sim}[d] & E_1[\ell^n]
\ar_{\alpha_1}^{\sim}[d]\\ 
(\zed/\ell^n\zed)^2 \ar@{=}[r] & (\zed/\ell^n\zed)^2.
}
$$

Set $x_{\ell}:=k_{\ell,1}^{-1}\circ T_{\ell}\phi_1$, then $x_{\ell}\in\End(T_{\ell} E_0)=\oh_{\ell}$.
Note that $x_{\ell}$ is not necessarily invertible in $\oh_{\ell}$! 
But it is nonzero, so $x_{\ell}\in B_{\ell}^\times$.  If $\ell\nmid N\cdot\deg\phi_1$
then $T_{\ell}\phi_1$ is actually an isomorphism extending
$\alpha_1^{-1}\circ\alpha_0$, so we can take $x_{\ell}=1$ in this case.  

How does this depend on the choice of $k_{\ell,1}$?  Let
$\tilde{k}_{\ell,1}:T_{\ell} E_0\to T_{\ell} E_1$ be an isomorphism extending
$\alpha_1^{-1}\circ\alpha_0:E_0[\ell^n]\to 
E_1[\ell^n]$.  Let $u=(\tilde{k}_{\ell,1})^{-1}\circ k_{\ell,1}\in\End(T_{\ell}
E_0)^\times=\oh_{\ell}^\times$.  Observe that $u$ restricts to the identity 
on $E_0[\ell^n]=(T_{\ell} E_0)/\ell^n(T_{\ell} E_0)$, so actually
$u\in\oh_{\ell}^\times(N)$.  Conversely, if $u\in\oh_{\ell}^\times(N)$, then
$k_{\ell,1}\circ u^{-1}:T_{\ell} E_0\to T_{\ell} E_1$ is an isomorphism extending
$\alpha_1^{-1}\circ\alpha_0$.  Therefore an isogeny $\phi_1:E_0\to E_1$
defines for each prime $\ell\neq p$ an element
$[x_{\ell}]\in\oh_{\ell}^\times(N)\bs B_{\ell}^\times$.  Moreover, $[x_{\ell}]=1$ for all
but finitely many $\ell$.

We need to carry out the same program at the prime $p$.  The isogeny
$\phi_1:E_0\to E_1$ induces an injective $W$-linear map
$M(\phi_1'):M_1\to M_0$ (note the contravariance!), whose cokernel is
isomorphic to $(\ker\phi_1')_p$.  By Corollary~\ref{cor:struct_p}
there exists an isomorphism of Dieudonn\'e modules $k_{p,1}:M_1\xrightarrow{\sim}M_0$ such
that the induced isomorphism $M_1/FM_1\to M_0/FM_0$ maps $\eta_1$ to
$\eta_0$; by Corollary~\ref{cor:isom_p} $k_{p,1}$ is well-defined up
to multiplication by $\oh_p^\times(1)$.  Set $x_p:=M(\phi_1')\circ
k_{p,1}^{-1}\in\End(M_0)=\oh_p$, then we get
$[x_p]\in\oh_p^\times(1)\bs B_p^\times$:
$$
\UseTips
\xymatrix{
M_1 \ar[r]^{M(\phi_1')} \ar[dr]_{k_{p,1}}^{\sim} & M_0\\
& M_0 \ar[u]_{x_p}.
}
$$

So far, our construction associates to a triple $(\phi_1,\alpha_1,\eta_1)$
an element $[x]\in U(1,N)\bs\eh_B^\times$.  This depends on the choice
of isogeny $E_0\to E_1$:
\begin{lem}\label{lem:isogenies}
Any isogeny $\tilde{\phi}_1:E_0\to E_1$ is of the form
$\tilde{\phi}_1=\phi_1\circ u$, where 
$$
u\in (\End(E_0)\otimes\cue)^\times=B^\times.
$$
\end{lem}
\begin{proof}
We treat $\phi_1$, $\tilde{\phi}_1$ as quasi-isogenies, i.e. elements of
$\Hom(E_0,E_1)\otimes\cue$.  Let $n=\deg\phi_1$, then we have that as
quasi-isogenies:
$$
(\hat{\phi_1}\otimes\frac{1}{n})\circ\phi_1 =
n\otimes\frac{1}{n}=1=\phi_1\circ(\hat{\phi_1}\otimes\frac{1}{n}).
$$
We can therefore write $\phi_1^{-1}=\hat{\phi_1}\otimes\frac{1}{n}$ and
we've shown that any isogeny has an inverse quasi-isogeny -- actually
a trivial modification of the argument shows that any quasi-isogeny is
invertible.  Now all we need to do is to set
$$
u=\phi_1^{-1}\circ\tilde{\phi}_1\in(\End(E_0)\otimes\cue)^\times=B^\times.
$$
\end{proof}

Therefore we cannot hope to get anything better than a map
$\Sigma^0\to U\bs\eh_B^\times/B^\times$.  We now show that
this map is well-defined, i.e. that it only depends on the isomorphism
class $[\phi_1,\alpha_1,\eta_1]$.  Let $f:E_1\to E_2$ be an isomorphism such
that the diagrams~(\ref{diag:elliptic_iso}) commute.  By
Lemma~\ref{lem:isogenies} it does not matter which isogeny
$\phi_2:E_0\to E_2$ we work with, so we might as well take
$\phi_2=f\circ\phi_1$.  We get the following diagrams
$$
\UseTips
\xymatrix{
T_{\ell} E_0 \ar[d]_{x_{\ell}} \ar[r]^{T_{\ell}\phi_1} \ar@(ur,ul)[rr]^{T_{\ell}\phi_2}
& T_{\ell} E_1 \ar[r]^{T_{\ell} f}_{\sim} & T_{\ell} E_2\\
T_{\ell} E_0 \ar@{=}[r] & T_{\ell} E_0 \ar[u]_{\sim}^{k_{{\ell},1}} \ar@{=}[r] & T_{\ell}
E_0 \ar[u]_{\sim}^{k_{{\ell},2}} 
}
\quad\quad
\xymatrix{
E_0[{\ell}^n] \ar[d]_{\alpha_0}^{\sim} \ar[r]^{k_{{\ell},1}}_{\sim}
\ar@(ur,ul)[rr]^{k_{{\ell},2}}_{\sim} & E_1[{\ell}^n] \ar[r]^{T_{\ell} f}_{\sim}
\ar[d]_{\alpha_1}^{\sim} & E_2[{\ell}^n] \ar[d]_{\alpha_2}^{\sim}\\
(\zed/{\ell}^n\zed)^2 \ar@{=}[r] & (\zed/{\ell}^n\zed)^2 \ar@{=}[r] &
(\zed/{\ell}^n\zed)^2  
},
$$
where $k_{{\ell},2}:=T_{\ell} f\circ k_{{\ell},1}$.  It is now clear that we end up
with the same $x_{\ell}\in\oh_{\ell}^\times(N)\bs B_{\ell}^\times$ as the one
obtained from $E_1$.  The exact same thing happens at the prime $p$.

We obtain a map 
$$
\gamma:\Sigma^0\to\Omega.
$$  

We need to construct an inverse.\label{page:inverse}  Let $[x]\in\Omega$ and choose a
representative $x=(x_v)\in\eh_B^\times$.  Let ${\ell}\neq p$.  We have
$x_{\ell}\in B_{\ell}^\times=\End(V_{\ell} E_0)^\times$, i.e. $x_{\ell}:V_{\ell} E_0\to V_{\ell}
E_0$ is a linear isomorphism.  But $V_{\ell} E_0=T_{\ell} E_0\otimes_{\zed_{\ell}}
\cue_{\ell}$ and $\cue_{\ell}=\zed_{\ell}[\frac{1}{{\ell}}]$.  Choose an isomorphism $T_{\ell} 
E_0\cong\zed_{\ell}^2$ extending $\alpha_0:E_0[{\ell}^n]\to(\zed/{\ell}^n\zed)^2$.
This gives us a matrix $A_{\ell}\in\GL_2(\cue_{\ell})$ which represents $x_{\ell}$.
Let $n_{\ell}\in\zed$ be the smallest integer such that ${\ell}^{n_{\ell}}A_{\ell}$ has
coefficients in $\zed_{\ell}$. 
Via our isomorphism this matrix corresponds to an endomorphism
$y_{\ell}\in\End(T_{\ell} E_0)$ which is injective with finite cokernel $C_{\ell}$.
Let ${\ell}^k$ be the order of $C_{\ell}$.  Let $K_{\ell}$ be the kernel of the map
induced by $y_{\ell}$ on $E_0[{\ell}^k]$: 
$$
0\to K_{\ell}\to E_0[{\ell}^k]\xrightarrow{y_{\ell}} E_0[{\ell}^k]\to C_{\ell}\to 0.
$$

For ${\ell}=p$ we have $x_p\in B_p^\times=(\End(M_0)\otimes\cue_p)^\times$.
Write $x_p=a+\pi b$, where $a,b\in L_p$ and $\pi^2=-p$.  We have
\begin{eqnarray*}
a=\sum_i a_i\otimes\frac{1}{p^i},\\
b=\sum_j b_j\otimes\frac{1}{p^j}.
\end{eqnarray*}
with $a_i,b_j\in\End(M_0)$.  Let $n_p\in\zed$ be the smallest integer
such that 
$$
p^{n_p}x_p=(a'\otimes 1)+\pi(b'\otimes 1)
$$ 
and set $y_p=a'+\pi b'\in\End(M_0)$.  This $y_p$ is an endomorphism of the Dieudonn\'e module
$M_0$ which induces an automorphism of $M_0\otimes\cue_p$, therefore
this endomorphism must be injective with finite cokernel $C_p$.  
Let $p^k$ be the order of $C_p$, then $y_p$ induces a
map  
$$
M(E_0[p^k])\xrightarrow{y_p} M(E_0[p^k])\to C_p\to 0.
$$
Then $C_p$ is the Dieudonn\'e module of a subgroup scheme $K_p$ of
$E_0$ of rank $p^k$.

Since $x\in\eh_B^\times$, $n_{\ell}=0$ for all but finitely many ${\ell}$.
Therefore it makes sense to set $q:=\prod {\ell}^{n_{\ell}}\in\cue^\times$ and
$y:=xq$; the ${\ell}$-th component of $y$ is precisely the $y_{\ell}$ above, and
clearly $[x]=[y]$.  Now set $K:=\bigoplus K_{\ell}$, then $K$ is a finite
subgroup of $E_0$.  So to the given $(x_v)\in\eh_B^\times$ we can
associate the quotient isogeny $\phi:E_0\to E_0/K$.  For ${\ell}\neq p$,
our construction gives for any positive integer $m$
$$
\UseTips
\xymatrix{
0\ar[r] & \ker \ar[r] \ar@{=}[d] & E_0[{\ell}^m] \ar[r]^{\phi} \ar@{=}[d] &
(E_0/K)[{\ell}^m]\\
0\ar[r] & \ker \ar[r] & E_0[{\ell}^m] \ar[r]^{y_{\ell}} & E_0[{\ell}^m].
}
$$
Due to the structure of ${\ell}^m$-torsion, it is not hard to see that one
can construct an isomorphism (actually, there exist many of them)
$E_0[{\ell}^m]\cong (E_0/K)[{\ell}^m]$ which makes the above diagram commute.
On the level of Tate modules, we get
$$
\UseTips\label{diag:hl}
\xymatrix{
0\ar[r] & T_{\ell} E_0 \ar[r]^>>>>>{T_{\ell}\phi}\ar@{=}[d] & T_{\ell}(E_0/K)\\
0\ar[r] & T_{\ell} E_0 \ar[r]^{y_{\ell}} & T_{\ell} E_0\ar[u]^{k_{\ell}}_{\sim}.
}
$$
In particular, we can set $\alpha:=\alpha_0\circ k_{\ell}^{-1}$, then the
isomorphisms $\alpha:(E_0/K)[{\ell}^n]\xrightarrow{\sim}(\zed/{\ell}^n\zed)^2$
for ${\ell}|N$ piece together to give a level $N$ structure on $E_0/K$.

Let $M:=M\left((E_0/K)'[p^\infty]\right)$.  For ${\ell}=p$ we have
similarly 
$$
\UseTips
\xymatrix{
0 \ar[r] & M \ar[d]^{\sim}_{k_p} \ar[r]^{M(\phi)} & M_0 \ar@{=}[d]
\ar[r] & \coker M(\phi) \ar[d]^{\sim}\ar[r] & 0\\ 
0 \ar[r] & M_0 \ar[r]^{y_p} & M_0 \ar[r] & C_p \ar[r] & 0,
}
$$
and $\eta:=k_p^{-1}(\eta_0)\in M/FM$ gives a nonzero invariant
differential on $E_0/K$.

The next result tells us that we have indeed constructed a map
$$
\delta:\Omega\to\Sigma^0.
$$
\begin{prop}
$\delta$ is well-defined.
\end{prop}
\begin{proof}
First suppose that $\bar{x}=xu$, where $u\in\End(E_0)$ is not
divisible by any rational prime.  Let ${\ell}\neq p$, then
$\bar{x}_{\ell}=x_{\ell}u$, so $\bar{y}_{\ell}=y_{\ell}u$:
$$
\UseTips
\xymatrix{
0\ar[r] & T_{\ell} E_0\ar[r]^{y_{\ell}} & T_{\ell} E_0\ar[r]\ar@{=}[d] & C_{\ell}\ar[r] &
0\\
0\ar[r] & T_{\ell} E_0\ar[r]^{\bar{y}_{\ell}}\ar[u]^{u} & T_{\ell} E_0\ar[r] &
\bar{C}_{\ell}\ar[r]\ar[u]^{v_{\ell}} & 0.
}
$$
The snake lemma gives
$$
\coker v_{\ell}=0,\quad\ker v_{\ell}\cong\coker u.
$$
Let ${\ell}^k$ be the order of $\bar{C}_{\ell}$, then we can restrict the above
diagram to the ${\ell}^k$-torsion and get
$$
\UseTips
\xymatrix{
0\ar[r] & K_{\ell}\ar[r] & E_0[{\ell}^k]\ar[r]^{y_{\ell}} & E_0[{\ell}^k]\ar@{=}[d]\ar[r]
& C_{\ell}\ar[r] & 0\\
0\ar[r] & \bar{K}_{\ell}\ar[u]^{g_{\ell}}\ar[r] &
E_0[{\ell}^k]\ar[u]^{u_{\ell}}\ar[r]^{\bar{y}_{\ell}} & 
E_0[{\ell}^k]\ar[r] & \bar{C}_{\ell}\ar[r]\ar[u]^{v_{\ell}} & 0,
}
$$
where $u_{\ell}$ is the restriction of $u$ to $E_0[{\ell}^k]$ and $g_{\ell}$ is the
restriction of $u$ to $\bar{K}_{\ell}$.  Note that $\coker(u_{\ell}:T_{\ell} E_0\to T_{\ell}
E_0)=\coker(u:E_0[{\ell}^k]\to E_0[{\ell}^k])$.  Since there's no snake lemma
for diagrams of long exact sequences, we split the above diagram in
two:
$$
\UseTips
\xymatrix @C=17pt{
0\ar[r] & K_{\ell}\ar[r] & E_0[{\ell}^k]\ar[r] & E_0[{\ell}^k]/\ker y_{\ell}\ar[r] & 0\\ 
0\ar[r] & \bar{K}_{\ell}\ar[r]\ar[u]^{g_{\ell}} & E_0[{\ell}^k]\ar[r]\ar[u]^{u_{\ell}} &
E_0[{\ell}^k]/\ker \bar{y}_{\ell}\ar[r]\ar[u]^{h_{\ell}} & 0
}
\quad\quad
\xymatrix @C=17pt{
0\ar[r] & \im y_{\ell}\ar[r] & E_0[{\ell}^k]\ar[r]\ar@{=}[d] & C_{\ell}\ar[r] & 0\\
0\ar[r] & \im \bar{y}_{\ell}\ar[r]\ar[u]^{h_{\ell}} & E_0[{\ell}^k]\ar[r] &
\bar{C}_{\ell}\ar[r]\ar[u]^{v_{\ell}} & 0,
}
$$
where I have taken the liberty of using the same label $h_{\ell}$ for two
maps which are canonically isomorphic.  We first apply the snake lemma
to the diagram on the right and get
$$
\ker h_{\ell}=0,\quad\coker h_{\ell}\cong\ker v_{\ell}.
$$
Using this information together with the snake lemma in the diagram on
the left gives
$$
\ker g_{\ell}\cong\ker u_{\ell},\quad 0\to\coker g_{\ell}\to\coker u_{\ell}\to\coker
h_{\ell}\to 0.
$$
But we already have $\coker u_{\ell}=\coker u\cong\ker v_{\ell}\cong\coker h_{\ell}$
so the short exact sequence above becomes $0\to\coker g_{\ell}\to 0$, i.e.
$\coker g_{\ell}=0$.

Let $g=\bigoplus g_{\ell}:\bar{K}\to K$ and let $f:E_0/\bar{K}\to E_0/K$ be
defined by the diagram 
$$
\UseTips
\xymatrix{
0\ar[r] & K\ar[r] & E_0\ar[r]^>>>>>{\phi} & E_0/K\ar[r] & 0\\
0\ar[r] & \bar{K}\ar[u]^g\ar[r] & E_0\ar[u]^u\ar[r]^>>>>{\bar{\phi}} &
E_0/\bar{K}\ar[u]^f\ar[r] & 0.
}
$$
We apply the snake lemma and get an exact sequence
$$
0\to\ker g\to\ker u\to\ker f\to\coker g=0\to\coker u=0\to\coker f\to 0.
$$
But the map $\ker g\to\ker u$ is the sum of the isomorphisms $\ker
g_{\ell}\cong\ker u_{\ell}$, so $\ker u\to\ker f$ is the zero map; therefore
$\ker f=0$.  Clearly $\coker f=0$, so $f$ is an isomorphism.  

We check that this isomorphism preserves level $N$ structures.  We
have a diagram 
$$
\UseTips
\xymatrix{
T_{\ell} E_0\ar[rr]^<<<<<<<<<<<<<<{T_{\ell}\phi}\ar[dr]_{y_{\ell}} & & T_{\ell}(E_0/K)\\
& T_{\ell} E_0\ar[ur]^{k_{\ell}}_{\sim}\ar[dr]_{\bar{k}_{\ell}}^{\sim}\\
T_{\ell} E_0\ar[rr]^<<<<<<<<<<<<<<{T_{\ell}\bar{\phi}}\ar[uu]^{T_{\ell}
u=u_{\ell}}\ar[ur]^{\bar{y}_{\ell}} & & 
T_{\ell}(E_0/\bar{K})\ar[uu]^{T_{\ell} f}_{\sim},
}
$$
where we know that the outer square commutes, and that the triangles
situated over, to the left, and under the central $T_{\ell} E_0$ commute.
Therefore the triangle to the right of the central $T_{\ell} E_0$ also
commutes, i.e. $k_{\ell}=T_{\ell} f\circ \bar{k}_{\ell}$.  The level $N$ structures on
$E_0/K$ and $E_0/\bar{K}$ are defined in such a way that the inner squares
in the following diagram commute:
$$
\UseTips
\xymatrix{
(E_0/K)[{\ell}^n]\ar[r]^>>>>>>{k_{\ell}^{-1}}_>>>>>>{\sim}\ar[d]_{\alpha}^{\sim}
\ar@(ur,ul)[rr]^{f}_{\sim} &
E_0[{\ell}^n]\ar[r]^<<<<<<<{\bar{k}_{\ell}}_<<<<<<<{\sim}\ar[d]_{\alpha_0}^{\sim} &
(E_0/\bar{K})[{\ell}^n]\ar[d]_{\bar{\alpha}}^{\sim}\\
(\zed/{\ell}^n\zed)^2\ar@{=}[r] & (\zed/{\ell}^n\zed)^2\ar@{=}[r] &
(\zed/{\ell}^n\zed)^2,
}
$$
therefore the outer rectangle also commutes, i.e. $f$ preserves the
level $N$ structures.

The same argument with reversed arrows shows that $f$ preserves
differentials.

Now suppose $\bar{x}=x{\ell}$, ${\ell}\neq p$ (the case ${\ell}=p$ is analogous, even
easier).  If ${\ell}'\nmid {\ell}p$, then $\bar{x}_{{\ell}'}=x_{{\ell}'}{\ell}$ and
$\bar{y}_{{\ell}'}=y_{{\ell}'}{\ell}$.  Multiplication by ${\ell}$ is an isomorphism of 
$T_{{\ell}'}E_0$, so it induces an isomorphism $\bar{K}_{{\ell}'}\cong K_{{\ell}'}$ by
applying the same argument as before on the diagram:
$$
\UseTips
\xymatrix{
0\ar[r] & K_{{\ell}'}\ar[r] & E_0[{\ell}'^k]\ar[r]^{y_{{\ell}'}} &
E_0[{\ell}'^k]\ar@{=}[d]\ar[r] & C_{{\ell}'}\ar[r] & 0\\
0\ar[r] & \bar{K}_{{\ell}'}\ar[u]^{g_{{\ell}'}}_{\sim}\ar[r] &
E_0[{\ell}'^k]\ar[u]^{{\ell}}_{\sim}\ar[r]^{\bar{y}_{{\ell}'}} & E_0[{\ell}'^k]\ar[r] &
\bar{C}_{{\ell}'}\ar[r]\ar[u]^{v_{{\ell}'}}_{\sim} & 0.
}
$$
Something similar occurs at $p$.  If ${\ell}'={\ell}$, we get $\bar{x}_{\ell}=x_{\ell}{\ell}$ and
$\bar{y}_{\ell}=y_{\ell}$ so $\bar{K}_{\ell}=K_{\ell}$.  We have an isomorphism
$\bar{K}\cong K$ so 
$E_0/\bar{K}\cong E_0/K$.  We need to check that this isomorphism is
compatible with the level structures and the differentials.  Let
${\ell}'\nmid {\ell}p$, then we have a diagram
$$
\UseTips
\xymatrix{
T_{\ell'}(E_0/K) \ar[d]_{\alpha} \ar@(ur,ul)[rr]^{\ell}_{\sim} & T_{\ell'} E_0
\ar[l]_<<<<<<{k_{{\ell}'}} \ar[d]_{\alpha_0} \ar[r]^>>>>>>{\bar{k}_{{\ell}'}} &
T_{{\ell}'}(E_0/\bar{K}) 
\ar[d]_{\bar{\alpha}}\\
(\zed/{\ell}'^n\zed)^2 \ar@{=}[r] & (\zed/{\ell}'^n\zed)^2 \ar@{=}[r] &
(\zed/{\ell}'^n\zed)^2.
}
$$
Since the top ``triangle'' commutes, we see that the level structures
commute with the isomorphism.  The same thing happens at $p$.  When
${\ell}'={\ell}$, then $\bar{K}_{\ell}=K_{\ell}$ so we get the same diagram as above, except
that the top isomorphism is actually the identity map.

It remains to check the local choices.  $C_{\ell}$ (therefore $K_{\ell}$) depends
on the chosen isomorphism $T_{\ell} E_0\cong\zed_{\ell}^2$, and this can change
$y_{\ell}$ by multiplication by an element of $\oh_{\ell}^\times(N)$.  Suppose
we have another such candidate $\bar{y}_{\ell}=u_{\ell}y_{\ell}$, then we would get a
commutative diagram 
$$
\UseTips
\xymatrix{
0\ar[r] & T_{\ell} E_0\ar@{=}[d]\ar[r]^{y_{\ell}} & T_{\ell} E_0\ar[r] & C_{\ell}\ar[r] & 0\\
0\ar[r] & T_{\ell} E_0\ar[r]^{\bar{y}_{\ell}} & T_{\ell} E_0\ar[r]\ar[u]^{u_{\ell}}_{\sim} &
\bar{C}_{\ell}\ar[u]^{v_{\ell}}_{\sim}\ar[r] & 0,
}
$$
from which we conclude as before that $\bar{K}_{\ell}\cong K_{\ell}$ and
$E_0/\bar{K}\cong E_0/K$.  For the level $N$ structure, we have the
diagram 
$$
\UseTips
\xymatrix{
(E_0/K)[{\ell}^n] \ar@(ur,ul)[rr]^{\alpha} \ar[r]^<<<<<<{k_{\ell}^{-1}} \ar[d]^{\sim}
& E_0[{\ell}^n] \ar@{=}[d] \ar[r]^<<<<<{\alpha_0} & (\zed/{\ell}^n\zed)^2
\ar@{=}[d]\\ 
(E_0/\bar{K})[{\ell}^n] \ar@(dr,dl)[rr]_{\bar{\alpha}}
\ar[r]^<<<<<{(\bar{k_{\ell}})^{-1}} 
& E_0[{\ell}^n] 
\ar[r]^<<<<<{\alpha_0} & (\zed/{\ell}^n\zed)^2
}
$$
and a similar argument holds for the $\eta$ and $\bar{\eta}$.
\end{proof}

\begin{lem}\label{lem:elliptic_bijection}
The map $\gamma$ is bijective with inverse $\delta$.
\end{lem}
\begin{proof}
Suppose we started with $[x]\in\Omega$ and got $[\phi:E_0\to 
E_0/K,\alpha,\eta]$.  For ${\ell}\neq p$ we get the exact sequence
(of $\zed_{\ell}$-modules)
$$
0\to T_{\ell} E_0\xrightarrow{T_{\ell}\phi} T_{\ell}(E_0/K)\to\coker T_{\ell}\phi\to 0.
$$
We see from diagram~(\ref{diag:hl}) that $y_{\ell}=k_{\ell}^{-1}\circ T_{\ell}\phi$,
where $k_{\ell}$ is an isomorphism that restricts to
$\alpha^{-1}\circ\alpha_0$.  Therefore $[y_{\ell}]$ is exactly the local
element that's obtained in the computation of
$\gamma([\phi,\alpha,\eta])$.  The same thing happens at $p$, so
indeed $\gamma\circ\delta=1$.

Conversely, suppose we start with a triple $(\phi:E_0\to
E,\alpha,\eta)$.  We get local elements $x_{\ell}$ forming an ad\`ele
$x$.  We have $\ker\phi=\prod_{\ell}\coker x_{\ell}$.  Now when we apply $\delta$
we already have $x_{\ell}\in\oh_{\ell}$ so $y_{\ell}=x_{\ell}$ and $K=\bigoplus\coker
x_{\ell}=\ker\phi$.  We get an isogeny $E_0\to E_0/K$ which has the same
kernel as $\phi$, therefore $E=E_0/K$.  It is clear from the
construction of $\delta$ that the level $N$ structure and the
invariant differential will stay the same.  
\end{proof}

We have just proved
\begin{thm}
There is a canonical bijection $\Sigma^0\to\Omega$.
\end{thm}

Before moving on we record the following consequence of our arguments.
\begin{cor}
Every class $[\phi,\alpha,\eta]\in\Sigma^0$ has a representative
satisfying $p\nmid\deg\phi$.
\end{cor}
\begin{proof}
Let $[x]=\gamma([\phi,\alpha,\eta])$.  Let $\tilde{x}$ be a
representative of $[x]$.  Let $B=\{L,\pi\}$ be a description of our
quaternion algebra, then $\pi$ is a uniformizer of $B_p$.  Therefore
$\tilde{x}_p=u\pi^k$, with $k\in\zed$, $u\in\oh_p^\times$.  But then
$\tilde{x}_p\pi^{-k}\in\oh_p^\times$.  Since $\pi\in\oh_p^\times(1)$,
we can use the representative $x=\tilde{x}\pi^{-k}$ for computing
$\delta$, in which case $C_p=1$ so the degree of the resulting isogeny
is prime to $p$.
\end{proof}

\begin{rem}
We note here that one can say more: the isogeny $\phi:E_0\to E$ can be
chosen to have square-free degree.  This is not easy to show.
\end{rem}

\subsection{Compatibilities}
It remains to prove the following result:
\begin{thm}
The canonical bijection $\gamma:\Sigma^0(N)\to\Omega(N)$ is compatible with the action of the Hecke algebra, with the action of $\GL_2(\zed/N\zed)$, and with the operation of raising the level.
\end{thm}

\subsubsection{Hecke action}\label{sect:hecke_elliptic}
In this section $\ell$ will denote a fixed prime not dividing $pN$.  For generalities on Hecke algebras, see \S{}3.1.2 in~\cite{andrianov2}.

If $H_1, H_2$ are subgroups of a group $G$, we say that $H_1$ is
\emph{commensurable} with $H_2$ (write $H_1\sim H_2$) if $H_1\cap H_2$ has
finite index in both $H_1$ and $H_2$.  If $H$ is a subgroup of $G$, we
define its \emph{commensurator} by
$$
\Comm(H):=\{g\in G:g^{-1}Hg\sim H\}.
$$
Finally, we say that $(G,H)$ is a \emph{Hecke pair}\index{Hecke pair} if $\Comm(H)=G$.  For example, one can find in \S{}3.2.1 of~\cite{andrianov2} a proof of the fact that $(\GL_2(\cue_\ell),\GL_2(\zed_\ell))$ is a Hecke pair.

If $(G,H)$ is a Hecke pair then any double coset $HgH$ has a finite
decomposition into left cosets.  Since the map $H\bs G\to G/H$ given
by $gH\mapsto Hg^{-1}$ is a bijection, we also know that $HgH$ has a
finite decomposition into the same number of right cosets, i.e. we can
write
$$
HgH=\bigsqcup_{j=1}^s x_jgH=\bigsqcup_{j=1}^s Hgy_j,
$$
where $x_j,y_j\in H$ for all $j$.

The \emph{Hecke algebra}\index{Hecke algebra} of the Hecke pair $(G,H)$ is by definition $\mathcal{H}(G,H):=\zed[H\bs G/H]$, with the multiplication described above.  From now on we set $G=\GL_2(\cue_\ell)$, $H=\GL_2(\zed_\ell)$ and $\mathcal{H}_\ell=\mathcal{H}(\GL_2(\cue_\ell),\GL_2(\zed_\ell))$.  We call $\mathcal{H}_\ell$ the \emph{local Hecke algebra}\index{Hecke algebra!local} at $\ell$.

Take an isogeny $\phi:E_1\to E_2$ whose degree is a power of $\ell$.  It induces an injective $\zed_\ell$-linear map $T_{\ell}\phi:T_{\ell} E_1\to T_{\ell} E_2$ which gives an element $g\in G=\GL_2(\cue_\ell)$.  Since $g$ is defined only up to changes of bases for $T_{\ell} E_1$ and $T_{\ell} E_2$, $\phi$ actually defines a double coset $HgH$, where $H=\GL_2(\zed_{\ell})$.  In this situation we say that $\phi$ is \emph{of type}\index{type!of an isogeny} $HgH$.  We say that a finite subgroup $C$ of an elliptic curve $E$ is \emph{of type}\index{type!of a finite subgroup} $HgH$ if the quotient isogeny $E\to E/C$ is of type $HgH$.

If $HgH\in\mathcal{H}_\ell$, we denote by $\det(HgH)$ the $\ell$-part of the determinant of any representative of $HgH$.  The action of $\mathcal{H}_\ell$ on $\Sigma^0$ is defined as follows.  If $\det(HgH)>1$, set
$$
T_{HgH}([\phi:E_0\to E,\alpha,\eta]):=\sum_{\text{$C\subset E$ of type $HgH$}}
[E_0\xrightarrow{\phi} E\xrightarrow{\psi_C} E/C, \alpha_C, \eta_C],
$$
where $\eta_C:=M(\psi_C')^{-1}(\eta)$, and $\alpha_C$ is
defined by the diagram
\begin{equation}
\label{diag:def_alpha,eta}
\UseTips
\xymatrix{
E[N] \ar[r]^<<<<<<<{\psi_C}_<<<<<<<{\sim} \ar[d]_{\alpha} & (E/C)[N]
\ar[d]_{\alpha_C}\\
(\zed/N\zed)^2 \ar@{=}[r] & (\zed/N\zed)^2.
}
\end{equation}
Note that these definitions make sense because $(\deg\psi_C,pN)=1$.

Now suppose $\det(HgH)<1$.  Given a finite subgroup $C$ of $E$ of type $Hg^{-1}H$, let $\hat{\psi}_C:(E/C)\to E$ be the dual isogeny to the quotient $\psi_C:E\to E/C$.  The action is defined by
$$
T_{HgH}([\phi:E_0\to E,\alpha,\eta]):=\sum_{\text{$C\subset E$ of type $Hg^{-1}H$}}
[E_0\xrightarrow{\phi} E\xleftarrow{\hat{\psi}_C} E/C, \alpha_C, \eta_C],
$$
where $\eta_C=M(\hat{\psi}_C')(\eta)$, and $\alpha_C$ is defined by the diagram
\begin{equation}
\label{diag:def_alpha,eta_quasi}
\UseTips
\xymatrix{
E[N] \ar[d]_{\alpha} & \ar[l]_<<<<<{\hat{\psi}_C}^<<<<<{\sim} (E/C)[N] \ar[d]_{\alpha_C}\\
(\zed/N\zed)^2 \ar@{=}[r] & (\zed/N\zed)^2.
}
\end{equation}

The algebra $\mathcal{H}_\ell$ acts on $H\bs G$ as follows: let $HgH=\coprod_i Hg_i$, let $Hx\in H\bs G$ and choose a representative $x\in Hx$.  Then there exist representatives $g_i\in Hg_i$ such that
$$
T_{HgH}(Hx)=\sum_i Hg_ix.
$$
The algebra $\mathcal{H}_\ell$ acts on $\Omega$ by acting on the component $Hx_l$ of $[x]\in\Omega$.

\begin{lem}
The bijection $\gamma:\Sigma^0\to\Omega$ is compatible with the action of the local Hecke algebra $\mathcal{H}_\ell$, i.e. for all $HgH\in\mathcal{H}_{\ell}$ and $[\phi,\alpha,\eta]$ we have
$$
\gamma\left(T_{HgH}([\phi,\alpha,\eta])\right)=
T_{HgH}(\gamma([\phi,\alpha,\eta])).
$$
\end{lem}
\begin{proof}
Let $HgH\in\mathcal{H}_\ell$, let $[\phi:E_0\to E,\alpha,\eta]\in\Sigma^0$ and let $[x]=\gamma([\phi,\alpha,\eta])$.  

Suppose at first that $\det(HgH)>1$ and let $C$ be a subgroup of $E$ of type $HgH$.  Let $[x_C]:=\gamma([\psi_C\circ\phi,\alpha_C,\eta_C])$.  If $(\ell',p\ell)=1$, we have a diagram
$$
\UseTips
\xymatrix{
T_{\ell'}E_0 \ar[r]^{T_{\ell'}\phi} \ar[d]_{x_{\ell'}} & T_{\ell'}E
\ar[r]^<<<<<{T_{\ell'}\psi_C}_<<<<<{\sim} & T_{\ell'}(E/C).\\
T_{\ell'}E_0 \ar[ur]_{\sim}^{k_{\ell'}}
}
$$
Since $(T_{\ell'}\psi_C)\circ k_{\ell'}:T_{\ell'}E_0\to T_{\ell'}(E/C)$ is an isomorphism restricting to $\alpha_C^{-1}\circ\alpha_0$ (see diagram~(\ref{diag:def_alpha,eta})), we get that $[x_{C,{\ell'}}]=[x_{\ell'}]$.

A similar argument, based on the following diagram, shows that $[x_{C,p}]=[x_p]$:
$$
\UseTips
\xymatrix{
M_C \ar[r]^{M(\psi_C')}_{\sim} & M \ar[r]^{M(\phi')}
\ar[dr]_{\sim}^{k_p} & M_0\\
& & M_0. \ar[u]_{x_p}
}
$$

Let's figure out what happens at $\ell$.  Fix $x_\ell\in Hx_\ell$, then the isomorphism $k_\ell:T_\ell E_0\to T_\ell E$ is fixed and allows us to identify these two $\zed_\ell$-modules.  Choose a $\zed_\ell$-linear isomorphism $k_C:T_\ell E\to T_\ell(E/C)$ and set $y_C=k_C^{-1}\circ T_\ell \psi_C$.  Via the identification $k_\ell$, $y_C$ induces a map $z_C:T_\ell E_0\to T_\ell E_0$.  We have a diagram
$$
\UseTips
\xymatrix{
T_{\ell} E_0 \ar[r]^{T_{\ell}\phi} \ar[d]_{x_{\ell}} & T_{\ell} E \ar[d]^{y_C} \ar[r]^<<<<<{T_{\ell}\psi_C} &
T_{\ell}(E/C).\\
T_{\ell} E_0 \ar[ur]_{k_{\ell}}^{\sim} \ar[d]_{z_C} & T_{\ell} E\ar[ur]_{k_C}^{\sim}\\
T_{\ell} E_0 \ar[ur]_{k_{\ell}}^{\sim}
}
$$
Since $k_C\circ k_\ell$ is an isomorphism $T_\ell E_0\to T_\ell (E/C)$ and $z_C\circ x_\ell$ satisfies all the properties $x_{C,\ell}$ should, we conclude that
$Hx_{C,\ell}=Hz_Cx_{\ell}$.  The assumption that $C$ is of type $HgH$ implies that $Hz_C\subset HgH$.

It remains to show that the map $C\mapsto Hz_C$ gives a bijection between the set of subgroups $C$ of $E$ of type $HgH$ and the set of right cosets $Hz$ contained in $HgH$.  We start by constructing an inverse map.  Let $Hz\subset HgH$ and pick a representative $z$.  This corresponds to a map $z:T_\ell E_0\to T_\ell E_0$, and hence induces via $k_\ell$ a map $y:T_\ell E\to T_\ell E$.  We use the same construction as in the definition of the map $\delta$ in \S{}\ref{sect:elliptic_bijection} (pages~\pageref{page:inverse} and following) to get a subgroup $C$ of $E$ which is canonically isomorphic to the cokernel of $y$.  $C$ will be of type $HgH$ because $Hz\subset HgH$.  The proof of the bijectivity of $C\mapsto z_C$ is now the same as the proof of Lemma~\ref{lem:elliptic_bijection}.

It remains to deal with the case $\det(HgH)<1$.  This works essentially the same, except that various arrows are reversed.  We illustrate the point by indicating how to obtain the equivalent of the map $C\mapsto Hz_C$ in this setting.  Let $C$ be a subgroup of $E$ of type $Hg^{-1}H$.  This defines a new element of $\Sigma^0$ which we denote by $[\hat{\psi}_C^{-1}\circ\phi,\alpha_C,\eta_C]$ (by a slight abuse of notation since $\hat{\psi}_C$ is not invertible as an isogeny).  Let $[x_C]:=\gamma([\hat{\psi}_C^{-1}\circ\phi,\alpha_C,\eta_C])$.  If $(\ell',p\ell)=1$, we have a diagram
$$
\UseTips
\xymatrix{
T_{\ell'}E_0 \ar[r]^{T_{\ell'}\phi} \ar[d]_{x_{\ell'}} & T_{\ell'}E &
\ar[l]_<<<<{T_{\ell'}\hat{\psi}_C}^<<<<{\sim} T_{\ell'}(E/C).\\
T_{\ell'}E_0 \ar[ur]_{\sim}^{k_{\ell'}}
}
$$
Since $(T_{\ell'}\hat{\psi}_C)^{-1}\circ k_{\ell'}:T_{\ell'}E_0\to T_{\ell'}(E/C)$ is an isomorphism restricting to $\alpha_C^{-1}\circ\alpha_0$ (see diagram~(\ref{diag:def_alpha,eta_quasi})), we get that $[x_{C,{\ell'}}]=[x_{\ell'}]$.  The situation at $p$ is similar and we have $[x_{C,p}]=[x_p]$.

What about $\ell$?  As before, we fix $x_\ell\in Hx_\ell$ and with it the isomorphism $k_\ell:T_\ell E_0\to T_\ell E$.  Choose an isomorphism $k_C:T_\ell E\to T_\ell (E/C)$ and set $y_C:=T_\ell\hat{\psi}_C\circ k_C$.  Via the identification $k_\ell$, $y_C$ induces a map $z_C:T_\ell E_0\to T_\ell E_0$.  We have a diagram
$$
\UseTips
\xymatrix{
T_{\ell} E_0 \ar[r]^{T_{\ell}\phi} \ar[d]_{x_{\ell}} & T_{\ell} E &\ar[l]_<<<<{T_{\ell}\hat{\psi}_C}
T_{\ell}(E/C).\\
T_{\ell} E_0 \ar[ur]_{k_{\ell}}^{\sim} & T_{\ell} E\ar[ur]_{k_C}^{\sim} \ar[u]_{y_C}\\
T_{\ell} E_0 \ar[ur]_{k_{\ell}}^{\sim} \ar[u]^{z_C}
}
$$
It is now clear that $z_C\circ x_{C,\ell}=x_\ell$.  $z$ is only defined up to right multiplication by elements of $H$ (because of the choice of $k_C$), so we get the formula $Hx_{C,\ell}=Hz_C^{-1}x_\ell$.  The assumption that $C$ is of type $Hg^{-1}H$ guarantees that $Hz_C^{-1}\subset HgH$.  The rest of the proof proceeds similarly to the case $\det(HgH)>1$.
\end{proof}

\subsubsection{Action of $\GL_2(\zed/N\zed)$}
Within this section we'll write $G$ to denote $\GL_2(\zed/N\zed)$.  The group $G$ acts on $\Sigma^0$ as follows: 
$$
g\cdot [\phi,\alpha,\eta]:=[\phi,g\circ\alpha,\eta].
$$
The action on $\Omega$ is more delicate.  It is easy to see that since $\oh_\ell=\End(T_\ell E_0)$, we have $\oh_\ell^\times(N)\bs \oh_\ell=\Aut(E_0[\ell^n])$, where $\ell^n\|N$.  Our fixed isomorphism $\alpha_0:E_0[N]\to(\zed/N\zed)^2$ identifies $G$ with $\Aut(E_0[N])$ via $g\mapsto \alpha_0^{-1}\circ g\circ\alpha_0$.  Therefore we get an identification
\begin{eqnarray*}
G &\xrightarrow{\sim}& \prod_\ell \oh_\ell^\times(N)\bs\oh_\ell^\times\\
g &\mapsto& \prod_\ell \oh_\ell^\times(N)(\alpha_0^{-1}\circ g\circ\alpha_0),
\end{eqnarray*}
where the product is finite since the terms with $\ell\nmid N$ are $1$.  The action of $G$ on $\Omega$ is then given by
$$
g\cdot \left[\prod_\ell \oh_\ell^\times(N)x_\ell\right]:=\left[\prod_\ell \oh_\ell^\times(N)(\alpha_0^{-1}\circ g\circ\alpha) x_\ell\right].
$$
\begin{lem}
The bijection $\gamma:\Sigma^0\to\Omega$ is compatible with the action of $G=\GL_2(\zed/N\zed)$.
\end{lem}
\begin{proof}
Let $\left[\prod \oh_\ell^\times(N)x_\ell\right]:=\gamma([\phi,\alpha,\eta])$ and $\left[\prod \oh_\ell^\times(N)x'_\ell\right]:=\gamma(g\cdot[\phi,\alpha,\eta])=\gamma([\phi,g\circ\alpha,\eta])$.  Pick some $\ell\neq p$ and set $H:=\oh_\ell^\times(N)$; we claim that $Hx'_\ell=H(\alpha_0^{-1}\circ g\circ\alpha)x_\ell$.  Recall that $x_\ell=k_\ell^{-1}\circ T_\ell\phi$, where $k_\ell:T_\ell E_0\to T_\ell E$ is some isomorphism extending $\alpha^{-1}\circ\alpha_0$.  Therefore $k_\ell':=k_\ell\circ (\alpha_0^{-1}\circ g\circ\alpha_0)$ is an isomorphism extending $\alpha^{-1}\circ g\circ \alpha_0$ and is thus precisely what we need in order to define $x_\ell'=(k'_\ell)^{-1}\circ T_\ell\phi$.  By the definition of $k'_\ell$ we have
$$
x'_\ell=(\alpha_0^{-1}\circ g^{-1}\circ\alpha) \circ k_\ell^{-1}\circ T_\ell\phi=(\alpha_0^{-1}\circ g^{-1}\circ\alpha) \circ x_\ell,
$$
which is what we wanted to show.
\end{proof}

\subsubsection{Raising the level}
Suppose $N'=dN$ for some positive integer $d$.  A level $N'$ structure $\alpha':E[N']\to(\zed/N'\zed)^2$ on an elliptic curve $E$ induces a level $N$ structure on $E$ in the following way.  Multiplication by $d$ on $E[N']$ gives a surjection $d:E[N']\to E[N]$, and there is a natural surjection $\pi:(\zed/N'\zed)^2\to (\zed/N\zed)^2$ given by reduction mod $N$.  We want to define an isomorphism $\alpha:E[N]\to (\zed/N\zed)^2$ that completes the following square
\begin{equation*}
\UseTips
\xymatrix{
E[N'] \ar[r]_<<<<<<{\sim}^<<<<<<{\alpha'} \ar@{->>}[d]_{d} & (\zed/N'\zed)^2 \ar@{->>}[d]^{\pi} \\
E[N] \ar@{.>}[r]^<<<<<<{\alpha} & (\zed/N\zed)^2.
}
\end{equation*}
This is straightforward: let $P\in E[N]$ and take some preimage $Q$ of it in $E[N']$.  Set $\alpha(P):=\pi(\alpha'(Q))$.  This is easily seen to be well-defined and a bijection.  We conclude that $[\phi,\alpha',\eta]\mapsto[\phi,\alpha,\eta]$ gives a map
$$
\Sigma^0(N')\to\Sigma^0(N).
$$

There is a similar map on the $\Omega$'s.  We only need to consider primes $\ell|N'$.  Here we have $\oh_\ell^\times(N')\subset\oh_\ell^\times(N)$ so we get maps $\oh_\ell^\times(N')\bs B_\ell^\times\to\oh_\ell^\times(N)\bs B_\ell$, which can be put together to form
$$
\Omega(N')\to\Omega(N).
$$

We want to show that the bijection $\gamma$ commutes with these maps.  This is clear at primes $\ell\nmid N'$, so suppose $\ell$ is a prime divisor of $N'$; say $\ell^m\|N$ and $\ell^n\|N'$.  Let $[\phi,\alpha',\eta]\in\Sigma^0(N')$, $[x']=\gamma([\phi,\alpha',\eta])$ and $[x]=\gamma([\phi,\alpha,\eta])$.  By definition, $x'_\ell=(k'_\ell)^{-1}\circ\phi$ where $k'_\ell:T_\ell E_0\to T_\ell E$ is an isomorphism restricting to
$$
\UseTips
\xymatrix{
E_0[\ell^n] \ar^{k'_\ell}_{\sim}[r] \ar_{\alpha'_0}^{\sim}[d] & E[\ell^n]
\ar_{\alpha'}^{\sim}[d]\\ 
(\zed/\ell^n\zed)^2 \ar@{=}[r] & (\zed/\ell^n\zed)^2.
}
$$
This defines the local component $\oh_\ell^\times(N')x'_\ell$.  We can restrict $k'_\ell$ even further to the $\ell^m$-torsion, and then by the definition of $\alpha$ we have
$$
\UseTips
\xymatrix{
E_0[\ell^m] \ar^{k'_\ell}_{\sim}[r] \ar_{\alpha'_0}^{\sim}[d] & E[\ell^m]
\ar_{\alpha}^{\sim}[d]\\ 
(\zed/\ell^m\zed)^2 \ar@{=}[r] & (\zed/\ell^m\zed)^2.
}
$$
But this means that $k_\ell'$ plays the role of the $k_\ell$ in the definition of $x_\ell$, so 
$$
\oh_\ell^\times(N)x'_\ell=\oh_\ell^\times(N)x_\ell.
$$  
This is precisely what the map $\Omega(N')\to\Omega(N)$ looks like at $\ell$, so we're done.

\section{Siegel modular forms}\label{chap:abelian}

Following a suggestion of Gross, we generalize the results of the previous section to dimension~$g\geq 2$.  More precisely, we link Siegel modular forms (mod~$p$) to reductions modulo~$p$ of modular forms on the algebraic group~$\GU_g(B)$, as defined in~\cite{gross2}.

Fix a prime~$p$ and an integer~$N>1$ prime to~$p$.

\subsection{The geometric theory of Siegel modular forms}
We review the basic definitions and results from~\cite{chai1}.  

All the schemes we consider are locally noetherian.  A~$g$-dimensional \emph{abelian scheme}\index{abelian scheme}~$A$ over a scheme~$S$ is a proper smooth group scheme
$$
\UseTips
\xymatrix{
A \ar[d]_{\pi}\\
S, \ar@/_.8pc/[u]_{0}
}
$$
whose (geometric) fibers are connected of dimension~$g$.

A \emph{polarization}\index{polarization} of~$A$ is an~$S$-homomorphism~$\lambda:A\to A^t=\Pic^0(A/S)$ such that for any geometric point~$s$ of~$S$, the homomorphism~$\lambda_s:A_s\to A_s^t$ is of the form~$\lambda_s(a)=t_a^*\mathcal{L}_s\otimes\mathcal{L}_s^{-1}$ for some ample invertible sheaf~$\mathcal{L}_s$ on~$A_s$.  Such~$\lambda$ is necessarily an isogeny.  In this case,~$\lambda_*\oh_A$ is a locally free~$\oh_{A^t}$-module whose rank is constant over each connected component of~$S$.  This rank is called the \emph{degree} of~$\lambda$; if this degree is~$1$ (so~$\lambda$ is an isomorphism) then~$\lambda$ is said to be \emph{principal}\index{polarization!principal}.  Any polarization is symmetric: $\lambda^t=\lambda$ via the canonical isomorphism~$A\cong A^{tt}$.

Let~$\phi:A\to B$ be an isogeny of abelian schemes over~$S$.  Cartier duality (Theorem~III.19.1 in~\cite{oort1}) states that~$\ker\phi$ is canonically dual to~$\ker\phi^t$.  There is a canonical non-degenerate pairing
$$
\ker\phi\times\ker\phi^t\to\gee_m.
$$
An important example is~$\phi=[N]$ for an integer~$N$.  The kernel~$A[N]$ of multiplication by~$N$ on~$A$ is a finite flat group scheme of rank~$N^{2g}$ over~$S$; it is \'etale over~$S$ if and only if~$S$ is a scheme over~$\zed[\frac{1}{N}]$.  We get the \emph{Weil pairing}\index{Weil pairing}
$$
A[N]\times A^t[N]\to\gee_m.
$$
A principal polarization~$\lambda$ on~$A$ induces a canonical non-degenerate skew-symmetric pairing
$$
A[N]\times A[N]\to\mu_N,
$$
which is also called the \emph{Weil pairing}.

For our purposes, a \emph{level $N$ structure}\index{level structure} on~$(A,\lambda)$ will be a symplectic similitude from~$A[N]$ with the Weil pairing to~$(\zed/N\zed)^{2g}$ with the standard symplectic pairing, i.e. an isomorphism of group schemes~$\alpha:A[N]\to(\zed/N\zed)^{2g}$ such that the following diagram commutes:
$$
\UseTips
\xymatrix{
A[N]\times A[N] \ar[r]^<<<<<{(f,f)} \ar[d]_{\text{Weil}} &
(\zed/N\zed)^{2g}\times(\zed/N\zed)^{2g} \ar[d]_{\text{std}}\\
\mu_N \ar[r]^<<<<<<<<<<<<<<<<{\sim} & {\zed/N\zed}
}
$$
for some isomorphism~$\mu_N\cong\zed/N\zed$.

If~$N\geq 3$, the functor ``isomorphism classes of principally polarized~$g$-dimensional abelian varieties with level~$N$ structure'' is representable by a scheme~$\mathcal{A}_{g,N}$ which is faithfully flat over~$\zed$, smooth and quasi-projective over~$\zed[\frac{1}{N}]$.  Let
$$
\UseTips
\xymatrix{
Y \ar[d]_{\pi}\\
{\mathcal{A}_{g,N}} \ar@/_.8pc/[u]_{0}
}
$$
be the corresponding universal abelian variety.  Let~$\mathbb{E}=0^*(\Omega_{Y/\mathcal{A}_{g,N}})$; this is called the \emph{Hodge bundle}\index{Hodge bundle}.

\subsubsection{Twisting the sheaf of differentials}\label{sect:twisting}
Let~$X$ be a scheme and let~$\mathcal{F}$ be a locally free~$\oh_X$-module whose rank is the same integer~$n$ on all connected components of~$X$.  Let~$\{U_i:i\in I\}$ be an open cover of~$X$ that trivializes~$\mathcal{F}$, then we have~$\mathcal{F}|_{U_i}\cong (\oh_X|_{U_i})^n$, and for all~$i$ and~$j$ we have isomorphisms~$\mathcal{F}|_{U_i\cap U_j}\cong\mathcal{F}|_{U_j\cap U_i}$ given by~$g_{ij}\in\GL_n(\oh_X|_{U_i\cap U_j})$ satisfying the usual cocycle identities.

Now suppose we are given a rational linear representation~$\rho:\GL_n\to\GL_m$.  We construct a new locally free~$\oh_X$-module~$\mathcal{F}_\rho$ as follows:  set~$(\mathcal{F}_\rho)_i=(\oh_X|_{U_i})^m$, and for any~$i,j$ define an isomorphism~$(\mathcal{F}_\rho)_i|_{U_i\cap U_j}\to(\mathcal{F}_\rho)_j|_{U_i\cap U_j}$ by~$\rho(g_{ij})\in\GL_m(\oh_X|_{U_i\cap U_j})$.  Since the transition functions~$\rho(g_{ij})$ satisfy the required properties, we can glue the~$(\mathcal{F}_\rho)_i$ together to get the locally free~$\oh_X$-module~$\mathcal{F}_\rho$.  We say that it was obtained by \emph{twisting~$\mathcal{F}$ by~$\rho$}.  It is obvious that $\mathcal{F}=\mathcal{F}_\text{std}$, where~$\text{std}:\GL_n\to\GL_n$ is the standard representation.

The correspondence~$\rho\mapsto\mathcal{F}_\rho$ is a covariant functor from the category of rational linear representations of~$\GL_n$ to the category of locally free~$\oh_X$-modules.  This functor is exact and it commutes with tensor products.

The scheme~$\mathcal{A}_{g,N}$ can be compactified in various ways (see~\cite{faltings1}).  Pick some compactification~$\mathcal{A}^*_{g,N}$ and let~$X$ denote its base-change to~$\fpbar$.  According to the Koecher principle, the Hodge bundle~$\mathbb{E}$ extends uniquely to a locally free sheaf~$\mathbb{E}$ on~$X$.

Given a rational representation~$\rho:\GL_g\to\GL_m$, the global sections of~$\mathbb{E}_\rho$ are called \emph{Siegel modular forms}\index{modular form!Siegel} (mod~$p$) of weight~$\rho$ and level~$N$ and they can be written
\begin{multline*}
M_\rho(N)=\H^0(X,\mathbb{E}_\rho)\\
=\left\{f:\{[A,\lambda,\alpha,\eta]\}\to\fpbar^m|f(A,\lambda,\alpha,M\eta)=\rho(M)^{-1}f(A,\lambda,\alpha,\eta),\forall M\in\GL_g(\fpbar)\right\},
\end{multline*}
where~$\eta$ is a basis of invariant differentials on~$A$.

\subsubsection{The Kodaira-Spencer isomorphism}\label{sect:kodaira-spencer}
We recall the properties of the Kodaira-Spencer isomorphism.  For a detailed account see \S{}III.9 and \S{}VI.4 in~\cite{faltings1}.

If~$\pi:A\to S$ is projective and smooth, there is a \emph{Kodaira-Spencer map}\index{Kodaira-Spencer map}
$$
\kappa:\mathcal{T}_S\to \R^1\pi_*(\mathcal{T}_{A/S}).
$$
If
$$
\UseTips
\xymatrix{
A \ar[d]_{\pi}\\
S, \ar@/_.8pc/[u]_{0}
}
$$
is an abelian scheme, set~$\mathbb{E}_{A/S}:=0^*(\Omega^1_{A/S})$.  Then
$$
\mathcal{T}_{A/S}=\pi^*(0^*(\mathcal{T}_{A/S}))=\pi^*(\mathbb{E}_{A/S}^\vee).
$$
The projection formula gives
$$
\R^1\pi_*(\pi^*(\mathbb{E}_{A/S}^\vee))=(\R^1\pi_*\oh_A)\otimes_{\oh_S}\mathbb{E}_{A/S}^\vee.
$$
Let~$\pi^t:A^t\to S$ be the dual abelian scheme, then
$$
\R^1\pi_*\oh_A=0^*(\mathcal{T}_{A^t/S})=(\mathbb{E}_{A^t/S})^\vee.
$$
So the Kodaira-Spencer map can be written as follows:
$$
\kappa:\mathcal{T}_S\to(\mathbb{E}_{A^t/S})^\vee\otimes_{\oh_S}\mathbb{E}_{A/S}^\vee,
$$
which after dualizing gives
$$
\kappa^\vee:\mathbb{E}_{A^t/S}\otimes_{\oh_S}\mathbb{E}_{A/S}\to\Omega^1_S.
$$
Now suppose that~$\lambda:A/S\to A^t/S$ is a principal polarization, i.e. an isomorphism.  Then the pullback map~$\lambda^*:\mathbb{E}_{A^t/S}\to\mathbb{E}_{A/S}$ is an isomorphism and we get a map
$$
\mathbb{E}_{A/S}^{\otimes 2}\to\Omega^1_S.
$$
This map factors through the projection map to~$\Sym^2(\mathbb{E}_{A/S})$, and the resulting map
$$
\Sym^2(\mathbb{E}_{A/S})\to\Omega^1_S
$$ 
is an isomorphism.  In particular, in the notation of \S{}\ref{sect:twisting} we have a Hecke isomorphism~$\mathbb{E}_{\Sym^2\text{std}}\cong\Omega^1_X$, where~$X=\mathcal{A}_{g,N}^*$.

\subsection{Superspecial abelian varieties}
For a commutative group scheme $X$ over a field $K$ we define the
\emph{$a$-number} of $X$ by 
$$
a(X)=\dim_K\Hom(\alpha_p,X).
$$
If $K\subset L$ then
$\dim_K\Hom(\alpha_p,X)=\dim_L\Hom(\alpha_p,X\otimes L)$ so $a(X)$
does not depend on the base field.

An abelian variety $A$ over $K$ of dimension $g\geq 2$ is said to be
\emph{superspecial}\index{abelian variety!superspecial} if $a(A)=g$.
Let $k$ be an algebraic closure of $K$.  By Theorem~2 of~\cite{oort1},
$a(A)=g$ if and only if 
$$
A\otimes k\cong E_1\times\ldots\times E_g,
$$
where the $E_i$ are supersingular elliptic curves over $k$.  On the
other hand, for any $g\geq 2$ and any supersingular elliptic curves
$E_1,\ldots,E_{2g}$ over $k$ we have (see Theorem~3.5 in~\cite{shioda1}) 
$$
E_1\times\ldots\times E_g\cong E_{g+1}\times\ldots\times E_{2g}.
$$
We conclude that $A$ is superspecial if and only if $A\otimes k\cong
E^g$ for some (and therefore any) supersingular elliptic curve $E$
over $k$.

Any abelian subvariety of a superspecial abelian variety $A$ is also
superspecial.  If $A$ is superspecial and $G\subset A$ is a finite
\'etale subgroup scheme, then $A/G$ is also superspecial. 

\begin{lem}\label{lem:canonical_superspec}
Let $A$ be a superspecial abelian variety over $\fpbar$.  Then $A$ has
a canonical $\fptwo$-structure $A'$, namely the one whose geometric
Frobenius is $[-p]$.  The correspondence $A\mapsto A'$ is functorial.
\end{lem}
\begin{proof}
Let $E$ be a supersingular elliptic curve over $\fpbar$, then $A\cong E^g$.  By Lemma~\ref{lem:canonical_supersing} we know that $E$ has an $\fptwo$-structure $E'$ with $\pi_{E'}=[-p]_{E'}$, therefore $A':=(E')^g$ is an $\fptwo$-structure for $A$ such that 
$$
\pi_{A'}=\pi_{E'}\times\pi_{E'}\times\ldots\times\pi_{E'}=[-p]_{E'}\times[-p]_{E'}\times\ldots\times[-p]_{E'}=[-p]_{A'}.
$$

The functoriality statement follows from the corresponding statement in Lemma~\ref{lem:canonical_supersing}.  Since any superspecial abelian variety over $\fpbar$ is isomorphic to $E^g$, it suffices to consider a morphism $f:E^g\to E^g$.  This is built of a bunch of morphisms $E\to E$, which by Lemma~\ref{lem:canonical_supersing} come from morphisms $E'\to E'$.  These piece together to give a morphism $f':(E')^g\to (E')^g$ over $\fptwo$, which is just $f$ after tensoring with $\fpbar$.
\end{proof}
An easy consequence of the functoriality is that if $\lambda$ is a principal polarization on $A$, there exists a principal polarization $\lambda'$ of the canonical $\fptwo$-structure $A'$ of $A$ such that $\lambda'\otimes\fpbar=\lambda$.  We say that $(A',\lambda')$ is the canonical $\fptwo$-structure of $(A,\lambda)$.

\subsubsection{Isogenies}\label{subsect:isogenies}
We need to define what it means for two principally polarized abelian
varieties $(A_1,\lambda_1)$ and $(A_2,\lambda_2)$ to be isogenous\index{isogeny!of polarized abelian varieties}.
The natural tendency is to consider isogenies $\phi:A_1\to A_2$ such
that the following diagram commutes:
$$
\UseTips
\xymatrix{
A_1 \ar[r]^{\phi} \ar[d]_{\lambda_1}^{\sim} & A_2
\ar[d]_{\lambda_2}^{\sim}\\
A_1^t & A_2^t \ar[l]_{\phi^t},
}
$$
i.e.
$$
\phi^t\circ\lambda_2\circ\phi=\lambda_1.
$$
But then $\deg\phi=1$ so the only isogenies that satisfy this
condition are isomorphisms.  We therefore relax the condition by
requiring $\phi$ to satisfy 
$$
\phi^t\circ\lambda_2\circ\phi=m\lambda_1,
$$
where $m\in\en$.  By computing degrees we get $(\deg\phi)^2=m^g$. 

We now consider the local data given by the presence of a
principal polarization.  Let $(A,\lambda)$ be a $g$-dimensional
principally polarized abelian variety defined over
$\fpbar$.  Let ${\ell}$ be a prime different from $p$ and set as
usual $\zed_{\ell}(1)=\invlim \mu_{{\ell}^n}$.  We have
the canonical Weil pairing (see \S{}12 of~\cite{milne2} or \S{}16 of~\cite{milne1})
$$
e_{\ell}:T_{\ell} A\times T_{\ell} A^t\to\zed_{\ell}(1),
$$
which is a non-degenerate $\zed_{\ell}$-bilinear map.  When combined with a
homomorphism of the form $\alpha:A\to A^t$ it gives 
\begin{eqnarray*}
e_{\ell}^\alpha:T_{\ell} A\times T_{\ell} A &\to& \zed_{\ell}(1)\\
(a,a') &\mapsto& e_{\ell}(a,\alpha a').
\end{eqnarray*}

If $\alpha$ is a polarization then $e_{\ell}^\alpha$ is an alternating
(also called symplectic) form, i.e.
$$
e_{\ell}^\alpha(a',a)=e_{\ell}^\alpha(a,a')^{-1}
$$
for all $a,a'\in T_{\ell} A$.  If $f:A\to B$ is a homomorphism, then
$$
e_{\ell}^{f^t\circ\alpha\circ f}(a,a')=e_{\ell}^\alpha(f(a),f(a'))
$$
for all $a,a'\in T_{\ell} A$, $\alpha:B\to B^t$.

Let $\phi:(A_1,\lambda_1)\to(A_2,\lambda_2)$ be an isogeny of
principally polarized abelian varieties.  $\phi$ induces an injective
$\zed_{\ell}$-linear map on Tate modules $T_{\ell}\phi:T_{\ell} A_1\to T_{\ell} A_2$, with
finite cokernel $T_{\ell} A_2/(T_{\ell} \phi)(T_{\ell} A_1)$ isomorphic to the
${\ell}$-primary part $(\ker\phi)_{\ell}$ of $\ker\phi$.  Since
$\phi^t\circ\lambda_2\circ\phi=m\lambda_1$, we have 
\begin{multline*}
e_{\ell}^{\lambda_2}((T_{\ell}\phi)a,(T_{\ell}\phi)a')=
e_{\ell}^{\phi^t\circ\lambda_2\circ\phi}(a,a')=
e_{\ell}^{m\lambda_1}(a,a') \\
=e_{\ell}(a,m\lambda_1 a')= e_{\ell}(a,\lambda_1a')^m=
e_{\ell}^{\lambda_1}(a,a')^m.
\end{multline*}
We say that $T_{\ell}\phi$ is a symplectic similitude between the
symplectic modules $(T_{\ell} A_1,e_{\ell}^{\lambda_1})$ and $(T_{\ell}
A_2,e_{\ell}^{\lambda_2})$.

What happens at $p$?  Let $W=W(k)$ for $k$ a perfect field of
characteristic $p$ and let $M$ be a free $W$-module with semi-linear
maps $F$ and $V$ satisfying
$$
FV=VF=p,\quad\quad Fx=x^pF,\quad\quad Vx=x^{1/p}V.
$$
A \emph{principal quasi-polarization}\index{quasi-polarization!principal} on $M$ is an alternating form
$$
e:M\times M\to W
$$
which is a perfect pairing over $W$, such that $F$ and $V$ are
adjoints:
$$
e(Fx,y)=e(x,Vy)^p.
$$

A principal polarization on $A$ defines a principal quasi-polarization
$e_p$ on the Dieudonn\'e module $M$ of $A$ (see Jeff Achter's thesis~\cite{achter3} for references).  If $A$ is defined over
$\fptwo$ then $e_p$ defines a hermitian form on $M/FM$ as follows:
\begin{eqnarray*}
M/FM\times M/FM &\to&\fptwo\\
(x,y) &\mapsto&\langle x,y\rangle:=e_p(x,Fy)^p\pmod{p}.
\end{eqnarray*}

An isogeny $\phi:(A_1,\lambda_1)\to(A_2,\lambda_2)$ induces a
symplectic similitude $\phi^*:M_2\to M_1$ of principally
quasi-polarized Dieudonn\'e modules.

\subsubsection{Quaternion hermitian forms}
Let $B$ be a quaternion algebra over a field $F$.  Let
$\overline{\cdot}$ denote the canonical involution of $B$ (i.e.
conjugation) and let $N$ denote the norm map.  Let $V$ be a left
$B$-module which is free of dimension $g$.  A \emph{quaternion hermitian
form}\index{quaternion hermitian form} on $V$ is an $F$-bilinear map
$$
f:V\times V\to B
$$
such that
$$
f(bx,y)=bf(x,y),\quad\quad \overline{f(x,y)}=f(y,x)
$$
for all $b\in B$, $x,y\in V$.  We say $f$ is non-degenerate if
$f(x,V)=0$ implies $x=0$.

The following result says that any quaternion hermitian form is
diagonalizable (see \S{}2.2 in~\cite{shimura1}):
\begin{prop}\label{prop:hermitian_diagonalize}
For every quaternion hermitian form $f(x,y)$ on $V$, there exists a
basis $\{x_1,\ldots,x_g\}$ of $V$ over $B$ such that
$$
f(x_i,x_j)=\alpha_i\delta_{ij}
$$
for $1\leq i,j\leq n$, where $\alpha_i\in F$.  Moreover if $f$ is
non-degenerate and the norm $N:B\to F$ is surjective, then there
exists a basis $\{y_1,\ldots,y_g\}$ of $V$ over $B$ such that
$$
f(y_i,y_j)=\delta_{ij}.
$$
\end{prop}

Furthermore, we have the following result (see \S{}3.4 of~\cite{vigneras1}):
\begin{thm}[The norm theorem]\label{thm:norm}
Let $B$ be a quaternion algebra over a field $F$, and let $F_B$ be the
set of elements of $F$ which are positive at all the real places of
$F$ which ramify in $B$.  Then the image of reduced norm map $n:B\to
F$ is precisely $F_B$.
\end{thm}
We conclude that if $B$ is the quaternion algebra over $\cue$ ramified
at $p$ and $\infty$, then $n(B)=\cue_{>0}$.

\subsubsection{The similitude groups}
Let $B$ be a quaternion algebra over a field $F$.  We define the group
of unitary $g\times g$ matrices and its similitude group by\index{unitary similitudes}
\begin{eqnarray*}
\U_g(B) &=& \{M\in\GL_g(B):M^*M=I\},\\
\GU_g(B) &=& \{M\in\GL_g(B):M^*M=\gamma(M)I,\gamma(M)\in F^\times\}.
\end{eqnarray*}
These are algebraic groups over $F$: let $(f_{ij})=M^*M$, then
$\U_g(B)$ is defined by the equations $f_{ij}=0$ ($i\neq j$),
$f_{ii}=1$, and $\GU_g(B)$ is defined by the equations 
$$
f_{ij}=0\text{ for }i\neq j,\quad f_{11}=f_{22}=\ldots=f_{gg}
$$
(these are automatically in $F$ because they are sums of norms of elements of $B$).

We define the group of symplectic $2g\times 2g$ matrices and its
similitude group as follows:\index{symplectic similitudes} 
\begin{eqnarray*}
\Sp_{2g}(F) &=& \{M\in\GL_{2g}(F):M^tJ_{2g}M=J_{2g}\},\\
\GSp_{2g}(F) &=&
\{M\in\GL_{2g}(F):M^tJ_{2g}M=\gamma(M)J_{2g},\gamma(M)\in F^\times\}, 
\end{eqnarray*}
where $J_{2g}=\left(\begin{smallmatrix} 0 & I_g\\ -I_g & 0\end{smallmatrix}\right)$.  

\begin{lem}
Let $K$ be a field.  The subgroups $\GU_g(M_2(K))$ and $\GSp_{2g}(K)$
are conjugate inside $\GL_{2g}(K)$.  In particular, they are
isomorphic and the $F$-algebraic group $\GU_g(B)$ is an $F$-form of
$\GSp_{2g}$.
\end{lem}
\begin{proof}
If $A=\left(\begin{smallmatrix}a&b\\c&d\end{smallmatrix}\right)\in M_2(K)$, then the
conjugate of $A$ is
$$
\bar{A}=
\left(\begin{smallmatrix}d&-b\\-c&a\end{smallmatrix}\right),
$$
therefore the adjoint of $A$ is
$$
A^*=\left(\begin{smallmatrix}d&-c\\-b&a\end{smallmatrix}\right)=J_2^{-1}AJ_2.
$$

Set $\tilde{J}_{2g}=\diag(J_2,\ldots,J_2)$ and let $M=(A_{ij})_{1\leq
  i,j\leq g}\in M_g(M_2(K))$.  We have 
$$
M^*=\tilde{J}_{2g}^{-1}M^t\tilde{J}_{2g}
$$
therefore
$$
M^*M=\tilde{J}_{2g}^{-1}M^t\tilde{J}_{2g}M.
$$

I claim that there exists a permutation matrix $P$ such that
$P^t\tilde{J}_{2g}P=J_{2g}$.  If $g$ is odd (resp. even), let
$P_1$ be the matrix corresponding to the product of transpositions 
$$
(2 (2g-1))(4(2g-3))\ldots((g-1)g), \text{ (resp. }(2 (2g-1))(4
(2g-3))\ldots(g(g+1))\text{)},
$$  
then
$$
\setcounter{MaxMatrixCols}{12}
P_1^t\tilde{J}_{2g}P_1=\left(\begin{smallmatrix}
& & & & & & & & & & 1 & 0\\
& & & & & & & & & & 0 & 1\\
& & & & & & & & & \;\;\;\;\text{\begin{rotate}{90}$\ddots$\end{rotate}} \\
& & & & & & & 1 & 0 \\
& & & & & & & 0 & 1 \\
& & & & & & 1 \\
& & & & & -1 \\
& & & -1 & 0 \\
& & & 0 & -1 \\
& & \;\;\;\;\text{\begin{rotate}{90}$\ddots$\end{rotate}} \\
-1 & 0 \\
0 & -1
\end{smallmatrix}\right),
$$
respectively
$$
P_1^t\tilde{J}_{2g}P_1=\left(\begin{smallmatrix}
& & & & & & & & 1 & 0\\
& & & & & & & & 0 & 1\\
& & & & & & & \;\;\;\;\text{\begin{rotate}{90}$\ddots$\end{rotate}} \\
& & & & & 1 & 0 \\
& & & & & 0 & 1 \\
& & & -1 & 0 \\
& & & 0 & -1 \\
& & \;\;\;\;\text{\begin{rotate}{90}$\ddots$\end{rotate}} \\
-1 & 0 \\
0 & -1
\end{smallmatrix}\right).
$$
In both cases it is clear what permutation matrix $P_2$ will finish
the job; set $P=P_2P_1$ and the claim is proved.

We have $J_{2g}=P^t\tilde{J}_{2g}P$, so $\tilde{J}_{2g}=PJ_{2g}P^t$
and
$$
M^*M=PJ_{2g}^{-1}P^tM^tPJ_{2g}P^tM.
$$
Now if $M\in\GU_g(M_2(K))$, then $M^*M=\gamma I$ for some $\gamma\in
K^\times$ and a little manipulation gives
$$
(P^tMP)^tJ_{2g}(P^tMP)=\gamma J_{2g},
$$
i.e. $P^tMP\in\GSp_{2g}(K)$.  Conversely, if $P^tMP\in\GSp_{2g}(K)$
then
$$
M^*M=PJ_{2g}^{-1}(P^tMP)^tJ_{2g}P^tM=PJ_{2g}^{-1}\gamma
J_{2g}P^t=\gamma I
$$
so $M\in\GU_g(M_2(K))$.  Therefore
$P^{-1}\GU_g(M_2(K))P=\GSp_{2g}(K)$, as desired.

Since $B\otimes\bar{F}\cong M_2(\bar{F})$, we conclude that
$\GU_g(B)\otimes\bar{F}\cong\GSp(\bar{F})$.
\end{proof}

\subsubsection{Algebraic modular forms (mod $p$)}\label{sect:algebraic_modular}
We give the definition of algebraic modular forms (mod $p$) on the group $G=\GU_g(B)$, where $B$ is the quaternion algebra over $\cue$ ramified at $p$ and $\infty$.  See~\cite{gross2} and~\cite{gross4} for more details.

The definition given by Gross requires $G$ to be a reductive algebraic group over $\cue$ which satisfies a technical condition for which it sufficient to know that $G_0(\arr)$ is a compact Lie group.  Our $G$ is reductive, being a form of the reductive group $\GSp_{2g}$.  We also know that $G_0(\arr)$ is compact, since it is a subgroup of the orthogonal group $\text{O}(4g)$.

Let $\oh_p$ be the maximal order of $B\otimes\cue_p$.  We define $U_p$ to be the kernel of the reduction modulo a uniformizer $\pi$ of $\oh_p$, i.e.
$$
1\to U_p\to G(\oh_p)\xrightarrow{\text{mod }\pi} \GU_g(\fptwo)\to 1.
$$
For $\ell\neq p$, we set 
$$
U_\ell(N):=\{x\in G(\zed_\ell):x\equiv 1\pmod{\ell^n},\ell^n\|N\}.
$$
The product
$$
U:=U_p\times\prod_{\ell\neq p} U_\ell(N)
$$
is an open compact subgroup of $G(\hat{\cue})$, called the level ($\hat{\cue}$ is the ring of finite ad\`eles).  Set $\Omega(N):=U\bs G(\hat{\cue})/G(\cue)$.  By Proposition 4.3 in~\cite{gross2}, the double coset space $\Omega(N)$ is finite.

Now let $\rho:\GU_g(\fptwo)\to\GL(W)$ be an irreducible representation, where $W$ is a finite-dimensional $\fp$-vector space.  We define the space of \emph{algebraic modular forms}\index{modular form!algebraic} (mod $p$) of weight $\rho$ and level $U$ on $G$ as follows:
$$
M(\rho,U):=\{f:\Omega(N)\to W:f(\lambda g)=\rho(\lambda)^{-1}f(g)\text{ for all }\lambda\in\GU_g(\fptwo)\}.
$$
Since $\Omega(N)$ is a finite set and $W$ is finite-dimensional, $M(\rho,U)$ is a finite-dimensional $\fp$-vector space.

\subsubsection{Differentials defined over $\fptwo$}
We know from Lemma~\ref{lem:canonical_superspec} that a principally polarized superspecial abelian variety $(A,\lambda)$ has a canonical $\fptwo$-structure $(A',\lambda')$.  We are therefore lead to consider invariant differentials on $A$ defined over $\fptwo$.  
\begin{lem}
Let $A$ be a superspecial abelian
variety over $\fpbar$, let $A'$ be its canonical
$\fptwo$-structure and $M=M(A'[p^\infty])$.  Then giving a basis of
invariant differentials on $A$ defined over $\fptwo$ is
equivalent to giving a basis of $M/FM$ over $\fptwo$.
\end{lem}
\begin{proof}
The invariant differentials on $A$ defined over $\fptwo$ are identified with $\omega(A')$, the invariant differentials on $A'$.  We have $\omega(A')\cong\omega(E'^g)\cong\omega(E')^g$.
By Lemma~\ref{lem:ec_diff} we know that
$$
\omega(E')\cong M(E'[p^\infty])/FM(E'[p^\infty]),
$$ 
and since $M(A'[p^\infty])\cong M(E'[p^\infty])^g$ we conclude that
$$
\omega(A')\cong M/FM.
$$
\end{proof}

Note that as we've seen in \S{}\ref{subsect:isogenies}, the presence of a principal polarization $\lambda'$ on an $\fptwo$-abelian variety $A'$ induces a hermitian form on the $g$-dimensional $\fptwo$-vector space $M/FM$.  We say that an invariant differential on $A$ defined over $\fptwo$ is an invariant differential on $(A,\lambda)$ if it respects this hermitian structure.  We can therefore conclude that
\begin{cor}
Let $(A,\lambda)$ be a principally polarized superspecial abelian
variety over $\fpbar$, let $(A',\lambda')$ be its canonical
$\fptwo$-structure and $M=M(A'[p^\infty])$.  Then giving a basis of
invariant differentials\index{invariant differential!on polarized superspecial abelian variety} on $(A,\lambda)$ defined over $\fptwo$ is
equivalent to giving a hermitian basis of $M/FM$ over $\fptwo$.
\end{cor}

\subsubsection{Dieudonn\'e module of a superspecial abelian variety}
We use the notation introduced in \S{}\ref{subsect:dieudonne_elliptic}.

We want to describe the structure of the principally quasi-polarized Dieudonn\'e module $M=M(E'[p^\infty]^g)$, where $E'$ is a supersingular elliptic curve defined over $\fptwo$.  What we need is a simple consequence of the following result (Proposition 6.1 from~\cite{oort5}):
\begin{prop}
Let $K$ be a perfect field containing $\fptwo$, and suppose $\{M,\langle,\rangle\}$ is a quasi-polarized superspecial Dieudonn\'e module of genus $g$ over $W=W(K)$ such that $M\cong A_{1,1}^g$.  Then one can decompose
$$
M\cong M_1\oplus M_2\oplus\ldots\oplus M_d\quad\quad\quad (\langle M_i,M_j\rangle=0\text{ if }i\neq j),
$$
where each $M_i$ is of either of the following types:
\begin{enumerate}
\renewcommand{\labelenumi}{(\roman{enumi})}
\item a genus $1$ quasi-polarized superspecial Dieudonn\'e module over $W$ generated by some $x$ such that $\langle x,Fx\rangle=p^r\epsilon$ for some $r\in\zed$ and $\epsilon\in W\setminus pW$ with $\epsilon^\sigma=-\epsilon$; or
\item a genus $2$ quasi-polarized superspecial Dieudonn\'e module over $W$ generated by some $x$, $y$ such that $\langle x,y\rangle=p^r$ for some $r\in\zed$, and $\langle x,Fx\rangle=\langle y,Fy\rangle=\langle x,Fy\rangle=\langle y,Fx\rangle=0$.
\end{enumerate}
\end{prop}

\begin{cor}\label{cor:dieudonne_superspecial}
$M(E'[p^\infty]^g)\cong A_{1,1}^g$ as principally quasi-polarized Dieudonn\'e modules, where $A_{1,1}^g$ is endowed with the product quasi-polarization. 
\end{cor}
\begin{proof}
In the direct sum decomposition of the proposition, the degree of the quasi-polariza\-tion on $M$ is the product of the degrees of the quasi-polarizations of each of the summands.  Since our $M$ is principally quasi-polarized we conclude that each summand is also principally quasi-polarized, i.e. the bilinear form $\langle,\rangle$ is a perfect pairing on each summand.  

Let $M_0$ be such a summand and suppose $M_0$ is of type (ii) from the proposition.  This gives a $W$-basis for $M_0$ consisting of $x$, $Fx$, $y$ and $Fy$.  The quasi-polarization $\langle, \rangle$ defines a map $M_0\to M_0^t$ given by $z\mapsto f_z$, where $f_z(v):=\langle z,v\rangle$.  Let $x^t$, $(Fx)^t$, $y^t$ and $(Fy)^t$ be the dual basis to $x$, $Fx$, $y$ and $Fy$.  It is an easy computation to see that $f_x=p^ry^t$, $f_{Fx}=p^{r+1}(Fy)^t$, $f_y=-p^rx^t$ and $f_{Fy}=-p^{r+1}(Fx)^t$.  For instance
$$
f_{Fy}(Fx)=\langle Fy,Fx\rangle=\langle y,VFx\rangle^\sigma=\langle y,px\rangle^\sigma=-p\langle x,y\rangle^\sigma=-p^{r+1}.
$$
But the map $M_0\to M_0^t$ given by $z\mapsto f_z$ is an isomorphism, hence $p^r=p^{r+1}=1$, contradiction.

So $M$ has only summands of type (i).  A similar (but even simpler) computation shows that each summand must have $\langle x,Fx\rangle=1$.
\end{proof}

\begin{cor}\label{cor:symplectic_superspecial}
There exists an isomorphism $\End(M,e_0)^\times\cong\GU_g(\oh_p)$, such that the subgroup of symplectic automorphisms which lift the identity map on $(M/FM,e_0)$ is identified with $U_p$ defined by the short exact sequence
$$
1\to U_p\to \GU_g(\oh_p)\to \GU_g(\fptwo)\to 1,
$$
where the surjective map is reduction modulo the uniformizer $\pi$ of $\oh_p$.
\end{cor}
\begin{proof}
Recall the identification $\End(A_{1,1})\cong \oh_p$ from the proof of Corollary~\ref{cor:isom_p}:
\begin{eqnarray*}
\varphi:\End(A_{1,1}) &\to & \oh_p\\
\left(\begin{smallmatrix} x & y\\-py^p & x^p\end{smallmatrix}\right) &\mapsto& x+\pi y.
\end{eqnarray*}
On the other hand, any $T\in\End(M)=\End(A_{1,1}^g)$ is a $2g\times 2g$ matrix made of $2\times 2$ blocks of the form
$$
T_{ij}=\left(\begin{smallmatrix}x_{ij} & y_{ij}\\-py_{ij}^p & x_{ij}^p\end{smallmatrix}\right).
$$
Therefore we have an isomorphism
\begin{eqnarray*}
\varphi:\End(M)^\times &\to& \GL_g(\oh_p)\\
T=(T_{ij})_{i,j} &\mapsto& (x_{ij}+\pi y_{ij})_{i,j}.
\end{eqnarray*}
We want to prove that under this isomorphism, $\End(M,e_0)^\times$ corresponds to $\GU_g(\oh_p)$.  For this we use Corollary~\ref{cor:dieudonne_superspecial}, which says that the bilinear form $e_0$ is given by the block-diagonal matrix
$$
E_0=\left(\begin{smallmatrix}0 & 1\\-1 & 0\\&&\ddots\\&&&0 & 1\\&&&-1 & 0\end{smallmatrix}\right).
$$
Therefore we have
$$
\End(M,e_0)^\times=\{T\in\End(M)^\times: T^tE_0T=\gamma E_0, \gamma\in\zed_p\}.
$$
Note that for the $2\times 2$ block $T_{ij}$ we have
$$
\left(\begin{smallmatrix}0 & 1\\-1 & 0\end{smallmatrix}\right)^{-1}T_{ij}^t \left(\begin{smallmatrix}0 & 1\\-1 & 0\end{smallmatrix}\right)=\left(\begin{smallmatrix}x_{ij}^p & -y_{ij}\\py_{ij}^p & x_{ij}\end{smallmatrix}\right),
$$
which maps under $\varphi$ to $x_{ij}^p-\pi y_{ij}=\overline{x_{ij}+\pi y_{ij}}=\overline{\varphi(T_{ij})}$, where $\bar{\cdot}$ denotes the conjugation in the quaternion algebra $B_p=\oh_p\otimes\cue_p$.  This means that $E_0^{-1} T^t E_0$ maps to $\varphi(T)^*$, where we write $U^*=\overline{U^t}$.  Putting it all together we conclude that for any $T\in\End(M)^\times$ we have
$$
T\in\End(M,e_0)^\times\iff E_0^{-1}T^tE_0T=\gamma\iff \varphi(T)^*\varphi(T)=\gamma\iff \varphi(T)\in\GU_g(\oh_p),
$$
which is precisely what we wanted to show.

For the second part of the statement note that 
$$
M/FM=\{(0,a_1,0,a_2,\ldots,0,a_g)^t+FM:a_i\in\fptwo\}.
$$  
Let $T=(T_{ij})\in\End(M,e_0)^\times$, then its induced map on $M/FM$ is
$$
T((0,a_1,0,a_2,\ldots,0,a_g)^t+FM)=\left(0,\sum_j a_j \bar{x}_{1j}^p,\ldots,0,\sum_j a_j \bar{x}_{gj}^p\right)+FM,
$$
where $\bar{x}_{ij}$ denotes the reduction modulo $\pi$ of $x_{ij}$.  Therefore $T$ induces the identity map on $M/FM$ if and only if
$$
\left(\begin{smallmatrix}\bar{x}_{11} &\bar{x}_{12} & \ldots & \bar{x}_{1g}\\ \bar{x}_{21} & \bar{x}_{22} & \ldots & \bar{x}_{2g}\\ \vdots & \vdots & \ddots & \vdots \\ \bar{x}_{g1} & \bar{x}_{g2} & \ldots & \bar{x}_{gg}\end{smallmatrix}\right)=1.
$$
But the matrix above is precisely the matrix of the reduction of $\varphi(T)$ modulo $\pi$, so $T$ induces the identity on $M/FM$ if and only if $\varphi(T)\in U_p$.
\end{proof}

\subsection{Construction of the bijection}\label{sect:abelian_bijection}
Let $A$ be a superspecial abelian variety of dimension $g$ over
$\fpbar$.  Let $A'\cong E'^g$ be its canonical $\fptwo$-structure,
then $A\cong E^g$ for $E:=E'\otimes\fpbar$.  Until further notice, I
will write $A$ to mean $E^g$ and $A'$ to mean $E'^g$.  Let
$\lambda_0'$ be the principal polarization on $A'$ defined by the
$g\times g$ identity matrix, let $\lambda_0:=\lambda_0'\otimes\fpbar$,
let $\alpha_0:A[N]\to(\zed/N\zed)^{2g}$ be a level $N$ structure on
$A$, and let $\eta_0$ be a basis of invariant
differentials on $(A,\lambda_0)$ defined over $\fptwo$ (i.e. a hermitian basis of $M/FM$), where
$M=M(A'[p^\infty])$.  The various Weil pairings induced by $\lambda_0$, resp. $\lambda'_0$ will be denoted $e_0$, resp. $e'_0$.

Let $\Sigma$ denote the finite set of isomorphism classes of pairs $(\lambda,\alpha)$, where $\lambda$ is a principal polarization on $A$ and $\alpha$ is a level $N$ structure.  $\Sigma$ is a subscheme of $X$.  We also define $\tilde{\Sigma}$ to be the set of isomorphism classes of triples $(\lambda,\alpha,\eta)$ with $\lambda$ and $\alpha$ as above and $\eta$ a basis of invariant differentials on $(A,\lambda)$ defined over $\fptwo$.  Isomorphism is given by the condition $f'(\eta_2)=\eta_1$ and the commutativity of the diagrams
\begin{equation}\label{diag:abelian_iso}
\UseTips
\xymatrix{
A \ar[r]^{f}_{\sim} \ar[d]_{\lambda_1}^{\sim} & A
\ar[d]_{\lambda_2}^{\sim}\\
A^t & A^t \ar[l]_{f^t}^{\sim}
}
\quad\quad
\xymatrix{
(A[N],e_1) \ar[r]^f_{\sim} \ar[d]_{\alpha_1}^{\sim} &
(A[N],e_2) \ar[d]_{\alpha_2}^{\sim} \\
((\zed/N\zed)^{2g},\text{std}) \ar@{=}[r] &
((\zed/N\zed)^{2g},\text{std}),
}
\end{equation}
where $\text{std}$ denotes the standard symplectic pairing on the various modules.

Let $\oh:=\End(E)$ and $B:=\oh\otimes\cue$.  Let $G:=\GU_g(B)$, and recall the notation of \S{}\ref{sect:algebraic_modular}.  The purpose of this section is to construct a bijection between the finite sets $\Sigma$ and $\Omega=\Omega(N)$, generalizing the result of \S{}\ref{sect:elliptic_bijection}.

\begin{lem}\label{lem:abelian_isog}
Given any principal polarization $\lambda$ on $A$, there exists an
isogeny of principally polarized abelian varieties
$\phi:(A,\lambda_0)\to(A,\lambda)$.
\end{lem}
\begin{proof}
We want an isogeny $\phi:A\to A$ such that 
$$
\phi^t\circ\lambda\circ\phi=m\lambda_0
$$
for some $m\in\en$.

There is an obvious bijective correspondence associating to a
homomorphism $\psi:A\to A$ a matrix $\Psi\in M_g(\oh)$.  Under
this bijection, $\psi^t:A^t\to A^t$ corresponds to the adjoint
$\Psi^*$.  If $\phi:A\to A$ is an isogeny, then $\Phi\in\GL_g(B)$.
If $\lambda:A\to A^t$ is a polarization, then
$\lambda^t=\lambda$ so $\Lambda^*=\Lambda$.  Also $\Lambda$ is
positive-definite.  If $\lambda$ is a principal polarization, then
$\Lambda\in\GL_g(\oh)$ defines a positive-definite quaternion
hermitian form $f$.  By Proposition~\ref{prop:hermitian_diagonalize}
we know that $\Lambda$ can be diagonalized, i.e. there exists
$M\in\GL_g(B)$ such that 
$$
M^{-1}\Lambda M=\diag(\alpha_1,\ldots,\alpha_g),
$$ 
with $\alpha_i\in\cue$.  The form $f$ is positive-definite so
$\alpha_i\in\cue_{>0}$.  But the norm theorem~\ref{thm:norm} says that
the norm map is surjective onto $\cue_{>0}$, so by the last part of
Proposition~\ref{prop:hermitian_diagonalize} there exists
$M'\in\GL_g(B)$ such that
$$
(M')^{-1}\Lambda M'=I.
$$

So there is a basis of $B^g$ such that the quaternion hermitian form
$f$ is represented by the matrix $I$.  But the matrices representing
$f$ are all of the form $Q^*\Lambda Q$ for $Q\in\GL_g(B)$.  Now
$B=\oh\otimes\cue$ so there exists a positive integer $n$ such that
$nQ$ has coefficients in $\oh$.  Let $\Phi=nQ$ and let $\phi:A\to A$
be the homomorphism corresponding to $\Phi$.  Since 
$\Phi\in\GL_g(B)$ and the fixed principal polarization $\lambda_0$
corresponds to the identity matrix, we conclude that $\phi$ is an
isogeny and 
$$
\phi^t\circ\lambda\circ\phi=n^2.
$$
\end{proof}

We can now mimic our approach from \S~\ref{chap:elliptic}.
Lemma~\ref{lem:abelian_isog} allows us to identify $\tilde{\Sigma}$ with
the set $\tilde{\Sigma}^0$ consisting of isomorphism classes of triples
$$
\left((A,\lambda_0)\xrightarrow{\phi}(A,\lambda),
\alpha:A[N]\to(\zed/N\zed)^{2g}, \eta\right),
$$
where $(A,\lambda_0)\xrightarrow{\phi} (A,\lambda)$ is an isogeny of principally polarized abelian varieties and isomorphism is defined by the diagrams~(\ref{diag:abelian_iso}).

\begin{prop}
An isogeny $\phi_1:(A,\lambda_0)\to(A,\lambda_1)$ defines for any prime
${\ell}\neq p$ an element $[x_{\ell}]\in U_{\ell}(N)\bs G_{\ell}$.  If ${\ell}\nmid\deg\phi_1$
then $[x_{\ell}]=1$.
\end{prop}
\begin{proof}
Pick a prime ${\ell}\neq p$ and let $n$ satisfy ${\ell}^n\|N$.  As we've seen in \S{}\ref{subsect:isogenies}, $\phi$ induces an injective symplectic similitude $T_{\ell}\phi_1:(T_{\ell} A,e_{\ell}^{\lambda_0})\to(T_{\ell} A,e_{\ell}^{\lambda_1})$, with finite cokernel isomorphic to $(\ker\phi_1)_{\ell}$.  To ease notation, we'll just write $e_0$ for $e_{\ell}^{\lambda_0}$ and $e_1$ for $e_{\ell}^{\lambda_1}$ (and we use the same letters for the corresponding Weil pairings on $A[{\ell}^n]$).

Let $k_{{\ell},1}:(T_{\ell} A,e_0)\to (T_{\ell} A,e_1)$ be a symplectic isomorphism whose restriction gives a commutative diagram
$$
\UseTips
\xymatrix{
(A[{\ell}^n],e_0) \ar^{k_{{\ell},1}}_{\sim}[r] \ar_{\alpha_0}^{\sim}[d] & (A[{\ell}^n],e_1) \ar_{\alpha_1}^{\sim}[d]\\ 
(\zed/{\ell}^n\zed)^{2g} \ar@{=}[r] & (\zed/{\ell}^n\zed)^{2g}.
}
$$

Let $x_{\ell}=k_{{\ell},1}^{-1}\circ T_{\ell}\phi_1$, then $x_{\ell}:(T_{\ell} A,e_0)\to (T_{\ell} A,e_0)$ is a symplectic similitude and sits in the commutative diagram 
\begin{equation}\label{diag:glab}
\UseTips
\xymatrix{
(T_{\ell} A,e_0) \ar[r]^{T_{\ell}\phi_1} \ar[d]_{x_{\ell}} & (T_{\ell} A,e_1)\\
(T_{\ell} A,e_0). \ar[ur]^{k_{{\ell},1}}_{\sim}
}
\end{equation}

The map $x_{\ell}$ is not necessarily invertible, but since it's injective with finite cokernel it defines a symplectic automorphism of $(V_{\ell} A,e_0)$, i.e. $x_{\ell}\in\GSp_{2g}(\cue_{\ell})=G_{\ell}$.  If ${\ell}\nmid\deg\phi$ then $T_{\ell}\phi$ is a symplectic isomorphism so we can take $x_{\ell}=1$.

How does this depend on the particular choice of $k_{{\ell},1}$?  Let $\tilde{k}_{{\ell},1}:(T_{\ell} A,e_0)\xrightarrow{\sim} (T_{\ell} A,e_1)$ be some other symplectic isomorphism that restricts to $\alpha_1^{-1}\circ\alpha_0$.  Let 
$$
u=(\tilde{k}_{{\ell},1})^{-1}\circ k_{{\ell},1}\in\GSp_{2g}(\zed_{\ell})=U_{\ell}.
$$  
Note that $u$ restricts to the identity on $A[{\ell}^n]$ so actually $u\in U_{\ell}(N)$.  Conversely, if $u\in U_{\ell}(N)$ then $k_{{\ell},1}\circ u^{-1}:(T_{\ell} A,e_0)\to (T_{\ell} A,e_1)$ is a symplectic isomorphism restricting to $\alpha_1^{-1}\circ\alpha$.  Therefore $\phi_1$ gives us a well-defined element $[x_{\ell}]\in U_{\ell}(N)\bs G_{\ell}$.
\end{proof}

What happens at $p$?  The isogeny $\phi_1$ induces an injective symplectic similitude 
$$
M(\phi'_1):(M,e_1)\to (M,e_0)
$$ 
with finite cokernel.  Let $k_{p,1}:(M,e_1)\to (M,e_0)$ be a symplectic isomorphism whose reduction $(M/FM,e_1)\to (M/FM,e_0)$ maps $\eta_1$ to $\eta_0$.  Set $x_p:=M(\phi'_1)\circ k_{p,1}^{-1}$, then the map $x_p:(M,e_0)\to (M,e_0)$ is an injective symplectic similitude with finite cokernel.  Hence $x_p$ induces a symplectic isomorphism of $(M\otimes\cue_p,e_0)$, so by Corollary~\ref{cor:symplectic_superspecial}, $x_p$ gives an element of $\GU_g(B_p)$.  Since $k_{p,1}$ is well-defined up to multiplication by $U_p$, we have that $\phi_1$ defines a element $[x_p]\in U_p\bs\GU_g(B_p)$.

\begin{lem}\label{lem:av_isogenies}
Any two isogenies $\phi_1,\tilde{\phi}_1:(A,\lambda_0)\to (A,\lambda_1)$ are related by $\tilde{\phi}_1=\phi_1\circ u$, where $u$ corresponds to a matrix $U\in\GU_g(B)$. 
\end{lem}
\begin{proof}
Suppose $\phi_1$, $\tilde{\phi}_1$ satisfy
\begin{eqnarray*}
\phi_1^t\circ\lambda_1\circ\phi_1&=&m\lambda_0,\\
\tilde{\phi}_1^t\circ\lambda_1\circ\tilde{\phi}_1&=&\tilde{m}\lambda_0.
\end{eqnarray*}

We treat $\phi_1$, $\tilde{\phi}_1$ as quasi-isogenies, i.e. elements of $\End(A)\otimes\cue$.  Let $n=\deg\phi_1$, then we have that as quasi-isogenies:
$$
\left(\hat{\phi_1}\otimes\frac{1}{n}\right)\circ\phi_1 =
n\otimes\frac{1}{n}=1=\phi_1\circ\left(\hat{\phi_1}\otimes\frac{1}{n}\right).
$$
We can therefore write $\phi_1^{-1}=\hat{\phi_1}\otimes\frac{1}{n}$ and we've shown that any isogeny has an inverse quasi-isogeny -- actually a trivial modification of the argument shows that any quasi-isogeny is invertible.  Set
$$
u=\phi_1^{-1}\circ\tilde{\phi}_1\in\left(\End(A)\otimes\cue\right)^\times.
$$

Denote by capital letters the matrices corresponding to the various maps.  We have
$$
U^*U = \tilde{\Phi}_1^*\left(\Phi_1^{-1}\right)^*\Phi_1^{-1}\tilde{\Phi}_1 =
\tilde{\Phi}_1^*\left(\frac{1}{m}\Lambda_1\right)\tilde{\Phi}_1 =\frac{\tilde{m}}{m}I
$$
so $U\in\GU_g(B)$.
\end{proof}

The next lemma says that we have indeed constructed a map 
$$
\gamma:\tilde{\Sigma}^0\to\Omega=U\bs G(\hat{\cue})/G(\cue).
$$

\begin{lem}
The map $\gamma$ is well-defined.
\end{lem}
\begin{proof}
We need to show that $\gamma$ only depends on the isomorphism class $[\phi_1,\alpha_1,\eta_1]$.  Suppose $f:(\phi_1,\alpha_1,\eta_1)\to(\phi_2,\alpha_2,\eta_2)$ is an isomorphism of triples.  By Lemma~\ref{lem:av_isogenies} we can assume without loss of generality that $\phi_2=f\circ\phi_1$.  For ${\ell}\neq p$, we get the following diagrams
$$
\UseTips
\xymatrix{
(T_{\ell} A,e_0) \ar[d]_{x_{\ell}} \ar[r]^{T_{\ell}\phi_1} \ar@(ur,ul)[rr]^{T_{\ell}\phi_2}
& (T_{\ell} A,e_1) \ar[r]^{T_{\ell} f}_{\sim} & (T_{\ell} A,e_2)\\
(T_{\ell} A,e_0) \ar@{=}[r] & (T_{\ell} A,e_0) \ar[u]_{\sim}^{k_{{\ell},1}} \ar@{=}[r] & (T_{\ell}
A,e_0) \ar[u]_{\sim}^{k_{{\ell},2}} 
}
\quad\quad
\xymatrix{
(A[{\ell}^n],e_0) \ar[d]_{\alpha_0}^{\sim} \ar[r]^{k_{{\ell},1}}_{\sim}
\ar@(ur,ul)[rr]^{k_{{\ell},2}}_{\sim} & (A[{\ell}^n],e_1) \ar[r]^{T_{\ell} f}_{\sim}
\ar[d]_{\alpha_1}^{\sim} & (A[{\ell}^n],e_2) \ar[d]_{\alpha_2}^{\sim}\\
(\zed/{\ell}^n\zed)^{2g} \ar@{=}[r] & (\zed/{\ell}^n\zed)^{2g} \ar@{=}[r] &
(\zed/{\ell}^n\zed)^{2g},
}
$$
where $k_{{\ell},2}:=T_{\ell} f\circ k_{{\ell},1}$.  It is now clear that we end up with the same $x_{\ell}\in\oh_{\ell}^\times(N)\bs B_{\ell}^\times$ as the one obtained from $\phi_1$.  The exact same thing happens at the prime $p$.
\end{proof}

We need to construct an inverse.\label{page:inverse_abelian}  Let $[x]\in\Omega$ and pick a representative $x=(x_v)\in G(\hat{\cue})$.  Let ${\ell}\neq p$.  We have $x_{\ell}\in G(\cue_\ell)=\GSp_{2g}(\cue_\ell)=\Aut(V_\ell,e_0)$.  Let $n_{\ell}\in\zed$ be the smallest integer such that $y_\ell:={\ell}^{n_{\ell}}x_{\ell}\in\GSp_{2g}(\zed_\ell)=\End(T_\ell A,e_0)$.  The endomorphism $y_\ell$ is injective with finite cokernel $C_{\ell}$.  Let ${\ell}^k$ be the order of $C_{\ell}$.  Let $K_{\ell}$ be the kernel of the map induced by $y_{\ell}$ on $A[{\ell}^k]$: 
$$
0\to K_{\ell}\to A[{\ell}^k]\xrightarrow{y_{\ell}} A[{\ell}^k]\to C_{\ell}\to 0.
$$

For ${\ell}=p$ we have $x_p\in\GU_g(B_p)=(\End(M,e_0)\otimes\cue_p)^\times$.  Write $x_p=a+\pi b$, where $a,b\in M_g(L_p)$ and $\pi^2=-p$.  We have
\begin{eqnarray*}
a=\sum_i a_i\otimes\frac{1}{p^i},\\
b=\sum_j b_j\otimes\frac{1}{p^j}.
\end{eqnarray*}
with $a_i,b_j\in\End(M,e_0)$.  Let $n_p\in\zed$ be the smallest integer such that 
$$
p^{n_p}x_p=(a'\otimes 1)+\pi(b'\otimes 1)
$$ 
and set $y_p=a'+\pi b'\in\End(M,e_0)$.  This $y_p$ is an endomorphism of the Dieudonn\'e module $M$ which induces an automorphism of $M\otimes\cue_p$, therefore this endomorphism must be injective with finite cokernel $C_p$.  Let $p^k$ be the order of $C_p$, then $y_p$ induces a map  
$$
M(A[p^k])\xrightarrow{y_p} M(A[p^k])\to C_p\to 0.
$$
Then $C_p$ is the Dieudonn\'e module of a subgroup scheme $K_p$ of $A$ of rank $p^k$.

Since $x\in G(\hat{\cue})$, $n_{\ell}=0$ for all but finitely many ${\ell}$.  Therefore it makes sense to set $q:=\prod {\ell}^{n_{\ell}}\in\cue^\times$ and $y:=xq$; the ${\ell}$-th component of $y$ is precisely the $y_{\ell}$ above, and clearly $[x]=[y]$.  Now set $K:=\bigoplus K_{\ell}$, then $K$ is a finite subgroup of $A$.  So to the given $[x]\in\Omega$ we can associate the quotient isogeny $A\to A/K$.  After picking an isomorphism $A/K\cong A$ we get an isogeny $\phi:A\to A$, and this induces a principal polarization $\lambda$ on $A$ such that $\phi$ is an isogeny of polarized abelian varieties.  For ${\ell}\neq p$, our construction gives for any positive integer $m$
$$
\UseTips
\xymatrix{
0\ar[r] & \ker \ar[r] \ar@{=}[d] & (A[{\ell}^m],e_0) \ar[r]^{\phi} \ar@{=}[d] &
(A[{\ell}^m],e)\\
0\ar[r] & \ker \ar[r] & (A[{\ell}^m],e_0) \ar[r]^{y_{\ell}} & (A[{\ell}^m],e_0).
}
$$
Due to the structure of ${\ell}^m$-torsion, it is not hard to see that one
can construct a symplectic isomorphism (actually, there exist many of them)
$(A[{\ell}^m],e_0)\cong (A[{\ell}^m],e)$ which makes the above diagram commute.
On the level of Tate modules, we get
$$
\UseTips\label{diag:hl_abelian}
\xymatrix{
0\ar[r] & (T_{\ell} A,e_0) \ar[r]^{T_{\ell}\phi}\ar@{=}[d] & (T_{\ell}A,e)\\
0\ar[r] & (T_{\ell} A,e_0) \ar[r]^{y_{\ell}} & (T_{\ell} A,e_0)\ar[u]^{k_{\ell}}_{\sim}.
}
$$
In particular, we can set $\alpha:=\alpha_0\circ k_{\ell}^{-1}$, then the symplectic isomorphisms 
$$
\alpha:(A[{\ell}^n],e)\xrightarrow{\sim}((\zed/{\ell}^n\zed)^{2g},\text{std})
$$
for ${\ell}|N$ piece together to give a level $N$ structure on $(A,\lambda)$.

For ${\ell}=p$ we have similarly 
$$
\UseTips
\xymatrix{
0 \ar[r] & (M,e) \ar[d]^{\sim}_{k_p} \ar[r]^{M(\phi)} & (M,e_0) \ar@{=}[d]
\ar[r] & \coker M(\phi) \ar[d]^{\sim}\ar[r] & 0\\ 
0 \ar[r] & (M,e_0) \ar[r]^{y_p} & (M,e_0) \ar[r] & C_p \ar[r] & 0,
}
$$
and $\eta:=k_p^{-1}(\eta_0)$ gives a nonzero invariant
differential on $(A,\lambda)$.

The next result tells us that we have indeed constructed a map
$$
\delta:\Omega\to\tilde{\Sigma}^0.
$$
\begin{prop}
The map $\delta$ is well-defined.
\end{prop}
\begin{proof}
First suppose that $\bar{x}=xu$, where $u\in\End(A,\lambda_0)$ is not
divisible by any rational prime.  Let ${\ell}\neq p$, then
$\bar{x}_{\ell}=x_{\ell}u$, so $\bar{y}_{\ell}=y_{\ell}u$:
$$
\UseTips
\xymatrix{
0\ar[r] & (T_{\ell} A,e_0)\ar[r]^{y_{\ell}} & (T_{\ell} A,e_0)\ar[r]\ar@{=}[d] & C_{\ell}\ar[r] &
0\\
0\ar[r] & (T_{\ell} A,e_0)\ar[r]^{\bar{y}_{\ell}}\ar[u]^{u} & (T_{\ell} A,e_0)\ar[r] &
\bar{C}_{\ell}\ar[r]\ar[u]^{v_{\ell}} & 0.
}
$$
The snake lemma gives
$$
\coker v_{\ell}=0,\quad\ker v_{\ell}\cong\coker u.
$$
Let ${\ell}^k$ be the order of $\bar{C}_{\ell}$, then we can restrict the above
diagram to the ${\ell}^k$-torsion and get
$$
\UseTips
\xymatrix{
0\ar[r] & K_{\ell}\ar[r] & (A[{\ell}^k],e_0)\ar[r]^{y_{\ell}} & (A[{\ell}^k],e_0)\ar@{=}[d]\ar[r]
& C_{\ell}\ar[r] & 0\\
0\ar[r] & \bar{K}_{\ell}\ar[u]^{g_{\ell}}\ar[r] &
(A[{\ell}^k],e_0)\ar[u]^{u_{\ell}}\ar[r]^{\bar{y}_{\ell}} & 
(A[{\ell}^k],e_0)\ar[r] & \bar{C}_{\ell}\ar[r]\ar[u]^{v_{\ell}} & 0,
}
$$
where $u_{\ell}$ is the restriction of $u$ to $A[{\ell}^k]$ and $g_{\ell}$ is the
restriction of $u$ to $\bar{K}_{\ell}$.  Note that $\coker(u_{\ell}:T_{\ell} A\to T_{\ell}
A)=\coker(u:A[{\ell}^k]\to A[{\ell}^k])$.  Since there's no snake lemma
for diagrams of long exact sequences, we split the above diagram in
two:
\begin{equation}\label{eqn:snake1}
\UseTips
\xymatrix{
0\ar[r] & K_{\ell}\ar[r] & (A[{\ell}^k],e_0)\ar[r] & (A[{\ell}^k],e_0)/\ker y_{\ell}\ar[r] & 0\\ 
0\ar[r] & \bar{K}_{\ell}\ar[r]\ar[u]^{g_{\ell}} & (A[{\ell}^k],e_0)\ar[r]\ar[u]^{u_{\ell}} &
(A[{\ell}^k],e_0)/\ker \bar{y}_{\ell}\ar[r]\ar[u]^{h_{\ell}} & 0,
}
\end{equation}
\begin{equation}\label{eqn:snake2}
\UseTips
\xymatrix{
0\ar[r] & \im y_{\ell}\ar[r] & (A[{\ell}^k],e_0)\ar[r]\ar@{=}[d] & C_{\ell}\ar[r] & 0\\
0\ar[r] & \im \bar{y}_{\ell}\ar[r]\ar[u]^{h_{\ell}} & (A[{\ell}^k],e_0)\ar[r] &
\bar{C}_{\ell}\ar[r]\ar[u]^{v_{\ell}} & 0,
}
\end{equation}
where I have taken the liberty of using the same label $h_{\ell}$ for two
maps which are canonically isomorphic.  We first apply the snake lemma
to diagram~\ref{eqn:snake2} and get
$$
\ker h_{\ell}=0,\quad\coker h_{\ell}\cong\ker v_{\ell}.
$$
Using this information together with the snake lemma in diagram~\ref{eqn:snake1} gives
$$
\ker g_{\ell}\cong\ker u_{\ell},\quad 0\to\coker g_{\ell}\to\coker u_{\ell}\to\coker
h_{\ell}\to 0.
$$
But we already have $\coker u_{\ell}=\coker u\cong\ker v_{\ell}\cong\coker h_{\ell}$
so the short exact sequence above becomes $0\to\coker g_{\ell}\to 0$, i.e.
$\coker g_{\ell}=0$.

Let $g=\bigoplus g_{\ell}:\bar{K}\to K$ and let $f:(A,\bar{\lambda})\to (A,\lambda)$ be defined by the diagram 
$$
\UseTips
\xymatrix{
0\ar[r] & K\ar[r] & (A,\lambda_0)\ar[r]^>>>>>{\phi} & (A,\lambda)\ar[r] & 0\\
0\ar[r] & \bar{K}\ar[u]^g\ar[r] & (A,\lambda_0)\ar[u]^u\ar[r]^>>>>{\bar{\phi}} &(A,\bar{\lambda})\ar[u]^f\ar[r] & 0,
}
$$
where we use some isomorphism $A/\bar{K}\cong A$ to define the isogeny $\bar{\phi}$ and the principal polarization $\bar{\lambda}$.  We apply the snake lemma and get an exact sequence
$$
0\to\ker g\to\ker u\to\ker f\to\coker g=0\to\coker u=0\to\coker f\to 0.
$$
But the map $\ker g\to\ker u$ is the sum of the isomorphisms $\ker
g_{\ell}\cong\ker u_{\ell}$, so $\ker u\to\ker f$ is the zero map; therefore
$\ker f=0$.  Clearly $\coker f=0$, so $f$ is an isomorphism.

We check that this isomorphism preserves level $N$ structures.  We
have a diagram 
$$
\UseTips
\xymatrix{
(T_{\ell}A, e_0)\ar[rr]^<<<<<<<<<<<<<<<<<{T_{\ell}\phi}\ar[dr]_{y_{\ell}} & & (T_{\ell}A,e)\\
& (T_{\ell} A,e_0)\ar[ur]_{k_{\ell}}^{\sim}\ar[dr]^{\bar{k}_{\ell}}_{\sim}\\
(T_{\ell} A,e_0)\ar[rr]^<<<<<<<<<<<<<<<<<{T_{\ell}\bar{\phi}}\ar[uu]^{T_{\ell}
u=u_{\ell}}\ar[ur]^{\bar{y}_{\ell}} & & 
(T_{\ell}A,\bar{e})\ar[uu]_{T_{\ell} f}^{\sim},
}
$$
where we know that the outer square commutes, and that the triangles
situated over, to the left, and under the central $(T_{\ell}A, e_0)$ commute.
Therefore the triangle to the right of the central $(T_{\ell}A, e_0)$ also
commutes, i.e. $k_{\ell}=T_{\ell} f\circ \bar{k}_{\ell}$.  The level $N$ structures on $(A,\lambda)$ and $(A,\bar{\lambda})$ are defined in such a way that the inner squares
in the following diagram commute:
$$
\UseTips
\xymatrix{
(A[{\ell}^n],e)\ar[r]^>>>>>>>>>{k_{\ell}^{-1}}_>>>>>>>>>{\sim}\ar[d]_{\alpha}^{\sim}
\ar@(ur,ul)[rr]^{f}_{\sim} & (A[{\ell}^n],e_0)\ar[r]^<<<<<<<<<<{\bar{k}_{\ell}}_<<<<<<<<<<{\sim}\ar[d]_{\alpha_0}^{\sim} & (A[{\ell}^n],\bar{e})\ar[d]_{\bar{\alpha}}^{\sim}\\
((\zed/{\ell}^n\zed)^{2g},\text{std})\ar@{=}[r] & ((\zed/{\ell}^n\zed)^{2g},\text{std})\ar@{=}[r] & ((\zed/{\ell}^n\zed)^{2g},\text{std}),
}
$$
therefore the outer rectangle also commutes, i.e. $f$ preserves the
level $N$ structures.

The same argument with reversed arrows shows that $f$ preserves
differentials.

Now suppose $\bar{x}=x{\ell}$, ${\ell}\neq p$ (the case ${\ell}=p$ is analogous, even
easier).  If ${\ell}'\nmid {\ell}p$, then $\bar{x}_{{\ell}'}=x_{{\ell}'}{\ell}$ and
$\bar{y}_{{\ell}'}=y_{{\ell}'}{\ell}$.  Multiplication by ${\ell}$ is an isomorphism of 
$(T_{{\ell}'}A,e_0)$, so it induces an isomorphism $\bar{K}_{{\ell}'}\cong K_{{\ell}'}$ by
applying the same argument as before on the diagram:
$$
\UseTips
\xymatrix{
0\ar[r] & K_{{\ell}'}\ar[r] & (A[{\ell}'^k],e_0)\ar[r]^{y_{{\ell}'}} &
(A[{\ell}'^k],e_0)\ar@{=}[d]\ar[r] & C_{{\ell}'}\ar[r] & 0\\
0\ar[r] & \bar{K}_{{\ell}'}\ar[u]^{g_{{\ell}'}}_{\sim}\ar[r] &
(A[{\ell}'^k],e_0)\ar[u]^{{\ell}}_{\sim}\ar[r]^{\bar{y}_{{\ell}'}} & (A[{\ell}'^k],e_0)\ar[r] &
\bar{C}_{{\ell}'}\ar[r]\ar[u]^{v_{{\ell}'}}_{\sim} & 0.
}
$$
Something similar occurs at $p$.  If ${\ell}'={\ell}$, we get $\bar{x}_{\ell}=x_{\ell}{\ell}$ and
$\bar{y}_{\ell}=y_{\ell}$ so $\bar{K}_{\ell}=K_{\ell}$.  We have an isomorphism
$\bar{K}\cong K$ so $(A,\bar{\lambda})\cong(A,\lambda)$.  We need to check that this isomorphism is
compatible with the level structures and the differentials.  Let
${\ell}'\nmid {\ell}p$, then we have a diagram
$$
\UseTips
\xymatrix{
(T_{\ell'}A,e) \ar[d]_{\alpha} \ar@(ur,ul)[rr]^{\ell}_{\sim} & (T_{\ell'}A,e_0)
\ar[l]_<<<<<<<<<<{k_{{\ell}'}} \ar[d]_{\alpha_0} \ar[r]^>>>>>>>>>>>{\bar{k}_{{\ell}'}} &
(T_{{\ell}'}A,\bar{e}) \ar[d]_{\bar{\alpha}}\\
((\zed/{\ell}'^n\zed)^{2g},\text{std}) \ar@{=}[r] & ((\zed/{\ell}'^n\zed)^{2g},\text{std}) \ar@{=}[r] &
((\zed/{\ell}'^n\zed)^{2g},\text{std}).
}
$$
Since the top ``triangle'' commutes, we see that the level structures
commute with the isomorphism.  The same thing happens at $p$.  When
${\ell}'={\ell}$, then $\bar{K}_{\ell}=K_{\ell}$ so we get the same diagram as above, except
that the top isomorphism is actually the identity map.

It remains to check the local choices.  The group $C_{\ell}$ (therefore $K_{\ell}$) depends
on the chosen isomorphism $(T_{\ell} A,e_0)\cong(\zed_{\ell}^{2g},\text{std})$, and this can change
$y_{\ell}$ by right multiplication by an element of $U_\ell(N)$.  Suppose
we have another such candidate $\bar{y}_{\ell}=u_{\ell}y_{\ell}$, then we would get a
commutative diagram 
$$
\UseTips
\xymatrix{
0\ar[r] & (T_{\ell} A,e_0)\ar@{=}[d]\ar[r]^{y_{\ell}} & (T_{\ell} A,e_0)\ar[r] & C_{\ell}\ar[r] & 0\\
0\ar[r] & (T_{\ell} A,e_0)\ar[r]^{\bar{y}_{\ell}} & (T_{\ell} A,e_0)\ar[r]\ar[u]^{u_{\ell}}_{\sim} &
\bar{C}_{\ell}\ar[u]^{v_{\ell}}_{\sim}\ar[r] & 0,
}
$$
from which we conclude as before that $\bar{K}_{\ell}\cong K_{\ell}$ and
$(A,\bar{\lambda})\cong(A,\lambda)$.  For the level $N$ structure, we have the
diagram 
$$
\UseTips
\xymatrix{
(A[{\ell}^n],e) \ar@(ur,ul)[rr]^{\alpha} \ar[r]^<<<<<<{k_{\ell}^{-1}} \ar[d]^{\sim}
& (A[{\ell}^n],e_0) \ar@{=}[d] \ar[r]^<<<<<{\alpha_0} & ((\zed/{\ell}^n\zed)^{2g},\text{std})
\ar@{=}[d]\\ 
(A[{\ell}^n],\bar{e}) \ar@(dr,dl)[rr]_{\bar{\alpha}}
\ar[r]^<<<<<{(\bar{k}_{\ell})^{-1}} 
& (A[{\ell}^n],e_0)
\ar[r]^<<<<<{\alpha_0} & ((\zed/{\ell}^n\zed)^{2g},\text{std})
}
$$
and a similar argument holds for the $\eta$ and $\bar{\eta}$.
\end{proof}

\begin{lem}\label{lem:abelian_bijection}
The map $\gamma$ is bijective with inverse $\delta$.
\end{lem}
\begin{proof}
Suppose we started with $[x]\in\Omega$ and got $[(A,\lambda_0)\xrightarrow{\phi} 
(A,\lambda),\alpha,\eta]$.  For ${\ell}\neq p$ we get the exact sequence
$$
0\to (T_{\ell} A,e_0)\xrightarrow{T_{\ell}\phi} (T_{\ell}A,e)\to\coker T_{\ell}\phi\to 0.
$$
We see from diagram~(\ref{diag:hl_abelian}) that $y_{\ell}=k_{\ell}^{-1}\circ T_{\ell}\phi$,
where $k_{\ell}$ is an isomorphism that restricts to
$\alpha^{-1}\circ\alpha_0$.  Therefore $[y_{\ell}]$ is exactly the local
element that's obtained in the computation of
$\gamma([\phi,\alpha,\eta])$.  The same thing happens at $p$, so
indeed $\gamma\circ\delta=1$.

Conversely, suppose we start with a triple $((A,\lambda_0)\xrightarrow{\phi}
(A,\lambda),\alpha,\eta)$.  We get local elements $x_{\ell}$ forming an ad\`ele
$x$.  We have $\ker\phi=\prod_{\ell}\coker x_{\ell}$.  Now when we apply $\delta$
we already have $x_{\ell}\in\GSp_{2g}(\zed_\ell)$ so $y_{\ell}=x_{\ell}$ and $K=\bigoplus\coker
x_{\ell}=\ker\phi$.  We get an isogeny $(A,\lambda_0)\to (A,\bar{\lambda})$ which has the same
kernel as $\phi$, therefore $(A,\bar{\lambda})\cong(A,\lambda)$.  It is clear from the
construction of $\delta$ that the level $N$ structure and the
invariant differential will stay the same.  
\end{proof}

We have just proved
\begin{thm}
There is a canonical bijection $\tilde{\Sigma}^0\to\Omega$.
\end{thm}

\subsection{Compatibilities}
We now turn to the proof of the following result:
\begin{thm}\label{thm:main_bijection}
The canonical bijection $\gamma:\tilde{\Sigma}^0(N)\to\Omega(N)$ is compatible with the action of the Hecke algebra, with the action of $\GSp_{2g}(\zed/N\zed)$, and with the operation of raising the level.
\end{thm}

\subsubsection{Hecke action}
In this section $\ell$ will denote a fixed prime not dividing $pN$.

An isogeny of polarized abelian varieties $\phi:(A_1,\lambda_1)\to (A_2,\lambda_2)$ is said to be an \emph{${\ell}$-isogeny}\index{l-isogeny@${\ell}$-isogeny} if its degree is a power of ${\ell}$.  Such $\phi$ induces a symplectic similitude 
$$
T_{\ell}\phi:(T_{\ell} A_1,e_1)\to (T_{\ell} A_2,e_2)
$$ 
which gives an element $g\in\GSp_{2g}(\cue_{\ell})$.  Since $g$ is defined only up to changes of symplectic bases for $T_{\ell} A_1$ and $T_{\ell} A_2$, $\phi$ actually defines a double coset $HgH$, where $H=\GSp_{2g}(\zed_{\ell})$.  We say that $\phi$ is \emph{of type}\index{type!of an isogeny} $HgH$.  If $C$ is a finite $\ell$-subgroup of an abelian variety $A$ and $\lambda$ is a principally polarization on $A$, there exists a principal polarization $\lambda_C$ on $A/C$ that makes the quotient isogeny $A\to A/C$ into an isogeny of principally polarized abelian varieties.  We say that $C$ is \emph{of type}\index{type!of a finite subgroup} $HgH$ if the quotient isogeny $(A,\lambda)\to (A/C,\lambda_C)$ is of type $HgH$.

Since $(\GSp_{2g}(\cue_\ell),\GSp_{2g}(\zed_\ell))$ is a Hecke pair (see \S{}3.3.1 of~\cite{andrianov2}), we can talk about the local Hecke algebra $\mathcal{H}_\ell:=\mathcal{H}(G,H)$, where $G=\GSp_{2g}(\cue_\ell)$ and $H=\GSp_{2g}(\zed_\ell)$.

If $HgH\in\mathcal{H}_\ell$, we denote by $\det(HgH)$ the $\ell$-part of the determinant of any representative of $HgH$.  The action of $\mathcal{H}_\ell$ on $\tilde{\Sigma}^0$ is defined as follows.  If $\det(HgH)>1$, let $C$ be a subgroup of $A$ of type $HgH$ and let $[(A,\lambda_0)\xrightarrow{\phi}(A,\lambda),\alpha,\eta]\in\tilde{\Sigma}^0$.  The abelian variety $A/C$ is also superspecial, so it can be identified with $A$.  We denote by $\psi_C$ the composition $A\to A/C\cong A$, and we denote by $\lambda_C$ the principal polarization induced on the image $A$.  We set
$$
T_{HgH}([(A,\lambda_0)\xrightarrow{\phi}(A,\lambda),\alpha,\eta]):=\sum_{\text{$C$ of type $HgH$}}
[(A,\lambda_0)\xrightarrow{\phi}(A,\lambda)\xrightarrow{\psi_C}(A, \lambda_C), \alpha_C, \eta_C],
$$
where $\eta_C:=M(\psi_C')^{-1}(\eta)$, and $\alpha_C$ is
defined by the diagram
\begin{equation}
\label{diag:def_alpha,eta_abelian}
\UseTips
\xymatrix{
(A[N],e) \ar[r]^<<<<<<<<<<{\psi_C}_<<<<<<<<<<{\sim} \ar[d]_{\alpha} & (A[N],e_C)
\ar[d]_{\alpha_C}\\
((\zed/N\zed)^{2g},\text{std}) \ar@{=}[r] & ((\zed/N\zed)^{2g},\text{std}).
}
\end{equation}
Note that these definitions make sense because $(\deg\psi_C,pN)=1$.

Now suppose $\det(HgH)<1$.  Given $C$ a subgroup of $A$ of type $Hg^{-1}H$, let $\psi_C$ be the composition $A\to A/C\cong A$ and let $\hat{\psi}_C:A\to A$ be the dual isogeny to $\psi_C$.  Given a principal polarization $\lambda$ on $A$, there is a principal polarization $\lambda_C$ on $A$ such that the following diagram commutes:
$$
\UseTips
\xymatrix{
A \ar[d]_{\lambda} & A \ar[d]^{\lambda_C} \ar[l]_{\hat{\psi}_C}\\
A^t \ar[r]^{(\hat{\psi}_C)^t}& A^t.
}
$$

The action is defined by
$$
T_{HgH}([(A,\lambda_0)\xrightarrow{\phi}(A,\lambda),\alpha,\eta]):=\sum_{\text{$C$ of type $Hg^{-1}H$}}
[(A,\lambda_0)\xrightarrow{\phi} (A,\lambda)\xleftarrow{\hat{\psi}_C} (A,\lambda_C), \lambda_C, \alpha_C, \eta_C],
$$
where $\eta_C=M(\hat{\psi}_C')(\eta)$, and $\alpha_C$ is defined by the diagram
\begin{equation}
\label{diag:def_alpha,eta_quasi_abelian}
\UseTips
\xymatrix{
(A[N],e) \ar[d]_{\alpha} & \ar[l]_<<<<<<<<<<{\hat{\psi}_C}^<<<<<<<<<<{\sim} (A[N],e_C) \ar[d]_{\alpha_C}\\
((\zed/N\zed)^{2g},\text{std}) \ar@{=}[r] & ((\zed/N\zed)^{2g},\text{std}).
}
\end{equation}

The algebra $\mathcal{H}_\ell$ acts on $H\bs G$ as follows: let $HgH=\coprod_i Hg_i$, let $Hx\in H\bs G$ and choose a representative $x\in Hx$.  Then there exist representatives $g_i\in Hg_i$ such that
$$
T_{HgH}(Hx)=\sum_i Hg_ix.
$$
The algebra $\mathcal{H}_\ell$ acts on $\Omega$ by acting on the component $Hx_l$ of $[x]\in\Omega$.

\begin{lem}
The bijection $\gamma:\tilde{\Sigma}^0\to\Omega$ is compatible with the action of the local Hecke algebra $\mathcal{H}_\ell$, i.e. for all $HgH\in\mathcal{H}_{\ell}$ and $[\phi,\alpha,\eta]$ we have
$$
\gamma\left(T_{HgH}([\phi,\alpha,\eta])\right)=
T_{HgH}(\gamma([\phi,\alpha,\eta])).
$$
\end{lem}
\begin{proof}
Let $HgH\in\mathcal{H}_\ell$, let $[(A,\lambda_0)\xrightarrow{\phi}(A,\lambda),\alpha,\eta]\in\Sigma^0$ and let $[x]=\gamma([\phi,\alpha,\eta])$.  

Suppose at first that $\det(HgH)>1$ and let $C$ be a subgroup of $A$ of type $HgH$.  Let $[x_C]:=\gamma([\psi_C\circ\phi,\alpha_C,\eta_C])$.  If $(\ell',p\ell)=1$, we have a diagram
$$
\UseTips
\xymatrix{
(T_{\ell'}A,e_0) \ar[r]^{T_{\ell'}\phi} \ar[d]_{x_{\ell'}} & (T_{\ell'}A,e)
\ar[r]^<<<<<{T_{\ell'}\psi_C}_<<<<<{\sim} & (T_{\ell'}A,e_C).\\
(T_{\ell'}A,e_0) \ar[ur]_{\sim}^{k_{\ell'}}
}
$$
Since $(T_{\ell'}\psi_C)\circ k_{\ell'}:(T_{\ell'}A,e_0)\to (T_{\ell'}A,e_C)$ is a symplectic isomorphism restricting to $\alpha_C^{-1}\circ\alpha_0$ (see diagram~(\ref{diag:def_alpha,eta_abelian})), we get that $[x_{C,{\ell'}}]=[x_{\ell'}]$.

A similar argument, based on the following diagram, shows that $[x_{C,p}]=[x_p]$:
$$
\UseTips
\xymatrix{
(M,e_C) \ar[r]^<<<<<{M(\psi_C')}_<<<<<{\sim} & (M,e) \ar[r]^{M(\phi')}
\ar[dr]_{\sim}^{k_p} & (M,e_0)\\
& & (M,e_0) \ar[u]_{x_p}.
}
$$

Let's figure out what happens at $\ell$.  Fix $x_\ell\in Hx_\ell$, then the symplectic isomorphism $k_\ell:(T_\ell A,e_0)\to (T_\ell A,e)$ is fixed and allows us to identify these two symplectic $\zed_\ell$-modules.  Choose a symplectic isomorphism $k_C:(T_\ell A,e)\to (T_\ell A,e_C)$ and set $y_C=k_C^{-1}\circ T_\ell \psi_C$.  Via the identification $k_\ell$, $y_C$ induces a map $z_C:(T_\ell A,e_0)\to (T_\ell A,e_0)$.  We have a diagram
$$
\UseTips
\xymatrix{
(T_{\ell} A,e_0) \ar[r]^{T_{\ell}\phi} \ar[d]_{x_{\ell}} & (T_{\ell} A,e) \ar[d]^{y_C} \ar[r]^<<<<<{T_{\ell}\psi_C} &
(T_{\ell}A,e_C).\\
(T_{\ell} A,e_0) \ar[ur]_{k_{\ell}}^{\sim} \ar[d]_{z_C} & (T_{\ell} A,e)\ar[ur]_{k_C}^{\sim}\\
(T_{\ell} A,e_0) \ar[ur]_{k_{\ell}}^{\sim}
}
$$
Since $k_C\circ k_\ell$ is a symplectic isomorphism $(T_\ell A,e_0)\to (T_\ell A,e_C)$ and $z_C\circ x_\ell$ satisfies all the properties $x_{C,\ell}$ should, we conclude that
$Hx_{C,\ell}=Hz_Cx_{\ell}$.  The assumption that $C$ is of type $HgH$ implies that $Hz_C\subset HgH$.

It remains to show that the map $C\mapsto Hz_C$ gives a bijection between the set of subgroups $C$ of $A$ of type $HgH$ and the set of right cosets $Hz$ contained in $HgH$.  We start by constructing an inverse map.  Let $Hz\subset HgH$ and pick a representative $z$.  This corresponds to a map $z:(T_\ell A,e_0)\to (T_\ell A,e_0)$, and hence induces via $k_\ell$ a map $y:(T_\ell A,e)\to (T_\ell A,e)$.  We use the same construction as in the definition of the inverse map $\delta$ in \S{}\ref{sect:abelian_bijection} (pages~\pageref{page:inverse_abelian} and following) to get a subgroup $C$ of $A$ which is canonically isomorphic to the cokernel of $y$.  This $C$ will be of type $HgH$ because $Hz\subset HgH$.  The proof of the bijectivity of $C\mapsto z_C$ is now the same as the proof of Lemma~\ref{lem:abelian_bijection}.

It remains to deal with the case $\det(HgH)<1$.  This works essentially the same, except that various arrows are reversed.  We illustrate the point by indicating how to obtain the equivalent of the map $C\mapsto Hz_C$ in this setting.  Let $C$ be a subgroup of $A$ of type $Hg^{-1}H$.  This defines a new element of $\tilde{\Sigma}^0$ which we denote by $[\hat{\psi}_C^{-1}\circ\phi,\alpha_C,\eta_C]$ (by a slight abuse of notation since $\hat{\psi}_C$ is not invertible as an isogeny).  Let $[x_C]:=\gamma([\hat{\psi}_C^{-1}\circ\phi,\alpha_C,\eta_C])$.  If $(\ell',p\ell)=1$, we have a diagram
$$
\UseTips
\xymatrix{
(T_{\ell'}A,e_0) \ar[r]^{T_{\ell'}\phi} \ar[d]_{x_{\ell'}} & (T_{\ell'}A,e) &
\ar[l]_{T_{\ell'}\hat{\psi}_C}^{\sim} (T_{\ell'}A,e_C).\\
(T_{\ell'}A,e_0) \ar[ur]_{\sim}^{k_{\ell'}}
}
$$
Since $(T_{\ell'}\hat{\psi}_C)^{-1}\circ k_{\ell'}:(T_{\ell'}A,e_0)\to (T_{\ell'}A,e_C)$ is a symplectic isomorphism restricting to $\alpha_C^{-1}\circ\alpha_0$ (see diagram~(\ref{diag:def_alpha,eta_quasi_abelian})), we get that $[x_{C,{\ell'}}]=[x_{\ell'}]$.  The situation at $p$ is similar and we have $[x_{C,p}]=[x_p]$.

What about $\ell$?  As before, we fix $x_\ell\in Hx_\ell$ and with it the symplectic isomorphism $k_\ell:(T_\ell A,e_0)\to (T_\ell A,e)$.  Choose a symplectic isomorphism $k_C:(T_\ell A,e)\to (T_\ell A,e_C)$ and set $y_C:=T_\ell\hat{\psi}_C\circ k_C$.  Via the identification $k_\ell$, $y_C$ induces a map $z_C:(T_\ell A,e_0)\to (T_\ell A,e_0)$.  We have a diagram
$$
\UseTips
\xymatrix{
(T_{\ell} A,e_0) \ar[r]^{T_{\ell}\phi} \ar[d]_{x_{\ell}} & (T_{\ell} A,e) &\ar[l]_{T_{\ell}\hat{\psi}_C}
(T_{\ell}A,e_C).\\
(T_{\ell} A,e_0) \ar[ur]_{k_{\ell}}^{\sim} & (T_{\ell} A,e)\ar[ur]_{k_C}^{\sim} \ar[u]_{y_C}\\
(T_{\ell} A,e_0) \ar[ur]_{k_{\ell}}^{\sim} \ar[u]^{z_C}
}
$$
It is now clear that $z_C\circ x_{C,\ell}=x_\ell$.  $z$ is only defined up to right multiplication by elements of $H$ (because of the choice of $k_C$), so we get the formula $Hx_{C,\ell}=Hz_C^{-1}x_\ell$.  The assumption that $C$ is of type $Hg^{-1}H$ guarantees that $Hz_C^{-1}\subset HgH$.  The rest of the proof proceeds similarly to the case $\det(HgH)>1$.
\end{proof}

\subsubsection{Action of $\GSp_{2g}(\zed/N\zed)$}
Within this section we'll write $G$ to denote $\GSp_{2g}(\zed/N\zed)$.  The group $G$ acts on $\tilde{\Sigma}^0$ as follows: 
$$
g\cdot [\phi,\lambda,\alpha,\eta]:=[\phi,\lambda,g\circ\alpha,\eta].
$$
The action on $\Omega$ is more delicate.  It is easy to see that since
$U_\ell=\Aut(T_\ell A,e_0)$, we have $U_\ell(N)\bs
U_\ell=\Aut(A[\ell^n],e_0)$, where $\ell^n\|N$.  Our fixed symplectic
isomorphism 
$$
\alpha_0:(A[N],e_0)\to((\zed/N\zed)^{2g},\text{std})
$$ 
identifies $G$ with $\Aut(A[N],e_0)$ via $g\mapsto \alpha_0^{-1}\circ g\circ\alpha_0$.  Therefore we get an identification
\begin{eqnarray*}
G &\xrightarrow{\sim}& \prod_\ell U_\ell(N)\bs U_\ell\\
g &\mapsto& \prod_\ell U_\ell(N)(\alpha_0^{-1}\circ g\circ\alpha_0),
\end{eqnarray*}
where the product is finite since the terms with $\ell\nmid N$ are $1$.  The action of $G$ on $\Omega$ is then given by
$$
g\cdot \left[\prod_\ell U_\ell(N)x_\ell\right]:=\left[\prod_\ell U_\ell(N)(\alpha_0^{-1}\circ g\circ\alpha) x_\ell\right].
$$
\begin{lem}
The bijection $\gamma:\tilde{\Sigma}^0\to\Omega$ is compatible with the action of $\GSp_{2g}(\zed/N\zed)$.
\end{lem}
\begin{proof}
Let $\left[\prod
  U_\ell(N)x_\ell\right]:=\gamma([\phi,\lambda,\alpha,\eta])$ and
$$
\left[\prod
  U_\ell(N)x'_\ell\right]:=\gamma(g\cdot[\phi,\lambda,\alpha,\eta])=\gamma([\phi,\lambda,g\circ\alpha,\eta]).
$$
  Pick some $\ell\neq p$ and set $H:=U_\ell(N)$; we claim that
  $Hx'_\ell=H(\alpha_0^{-1}\circ g\circ\alpha)x_\ell$.  Recall that
  $x_\ell=k_\ell^{-1}\circ T_\ell\phi$, where $k_\ell:(T_\ell
  A,e_0)\to (T_\ell A,e)$ is some symplectic isomorphism extending
  $\alpha^{-1}\circ\alpha_0$.  Therefore $k_\ell':=k_\ell\circ
  (\alpha_0^{-1}\circ g\circ\alpha_0)$ is a symplectic isomorphism extending $\alpha^{-1}\circ g\circ \alpha_0$ and is thus precisely what we need in order to define $x_\ell'=(k'_\ell)^{-1}\circ T_\ell\phi$.  By the definition of $k'_\ell$ we have
$$
x'_\ell=(\alpha_0^{-1}\circ g^{-1}\circ\alpha) \circ k_\ell^{-1}\circ T_\ell\phi=(\alpha_0^{-1}\circ g^{-1}\circ\alpha) \circ x_\ell,
$$
which is what we wanted to show.
\end{proof}

\subsubsection{Raising the level}
Suppose $N'=dN$ for some positive integer $d$.  A level $N'$ structure 
$$
\alpha':(A[N'],e)\to((\zed/N'\zed)^{2g},\text{std})
$$ 
on the principally polarized abelian variety $(A,\lambda)$ induces a level $N$ structure on $(A,\lambda)$ in the following way.  Multiplication by $d$ on $A[N']$ gives a surjection $d:A[N']\to A[N]$, and there is a natural surjection $\pi:(\zed/N'\zed)^{2g}\to (\zed/N\zed)^{2g}$ given by reduction mod $N$.  We want to define a map $\alpha:A[N]\to (\zed/N\zed)^{2g}$ that completes the following square
\begin{equation*}
\UseTips
\xymatrix{
A[N'] \ar[r]_<<<<<<{\sim}^<<<<<<{\alpha'} \ar@{->>}[d]_{d} & (\zed/N'\zed)^{2g} \ar@{->>}[d]^{\pi} \\
A[N] \ar@{.>}[r]^<<<<<<{\alpha} & (\zed/N\zed)^{2g}
}
\end{equation*}
This is straightforward: let $P\in A[N]$ and take some preimage $Q$ of it in $A[N']$.  Set $\alpha(P):=\pi(\alpha'(Q))$.  This is easily seen to be well-defined and a bijection.  Since both surjections $d$ and $\pi$ respect the symplectic structure, $\alpha$ is a symplectic isomorphism.  We conclude that $[\phi,\lambda,\alpha',\eta]\mapsto[\phi,\lambda,\alpha,\eta]$ gives a map
$$
\tilde{\Sigma}^0(N')\to\tilde{\Sigma}^0(N).
$$

There is a similar map on the $\Omega$'s.  We only need to consider primes $\ell|N'$.  Here we have $U_\ell(N')\subset U_\ell(N)$ so we get maps $U_\ell(N')\bs G_\ell\to U_\ell(N)\bs G_\ell$, which can be put together to form
$$
\Omega(N')\to\Omega(N).
$$

We want to show that the bijection $\gamma$ commutes with these maps.  This is clear at primes $\ell\nmid N'$, so suppose $\ell$ is a prime divisor of $N'$; say $\ell^m\|N$ and $\ell^n\|N'$.  Choose elements $[\phi,\lambda,\alpha',\eta]\in\tilde{\Sigma}^0(N')$, $[x']:=\gamma([\phi,\lambda,\alpha',\eta])$ and $[x]:=\gamma([\phi,\lambda,\alpha,\eta])$.  By definition, we have $x'_\ell=(k'_\ell)^{-1}\circ\phi$ where $k'_\ell:(T_\ell A,e_0)\to (T_\ell A,e)$ is a symplectic isomorphism restricting to
$$
\UseTips
\xymatrix{
(A[\ell^n],e_0) \ar^{k'_\ell}_{\sim}[r] \ar_{\alpha'_0}^{\sim}[d] & (A[\ell^n],e)
\ar_{\alpha'}^{\sim}[d]\\ 
((\zed/\ell^n\zed)^{2g},\text{std}) \ar@{=}[r] & ((\zed/\ell^n\zed)^{2g},\text{std}).
}
$$
This defines the local component $U_\ell(N')x'_\ell$.  We can restrict $k'_\ell$ even further to the $\ell^m$-torsion, and then by the definition of $\alpha$ we have
$$
\UseTips
\xymatrix{
(A[\ell^m],e_0) \ar^{k'_\ell}_{\sim}[r] \ar_{\alpha'_0}^{\sim}[d] & (A[\ell^m],e)
\ar_{\alpha}^{\sim}[d]\\ 
((\zed/\ell^m\zed)^{2g},\text{std}) \ar@{=}[r] & ((\zed/\ell^m\zed)^{2g},\text{std}).
}
$$
But this means that $k_\ell'$ plays the role of the $k_\ell$ in the definition of $x_\ell$, so $U_\ell(N)x'_\ell=U_\ell(N)x_\ell$.  This is precisely what the map $\Omega(N')\to\Omega(N)$ looks like at $\ell$, so we're done.

\subsection{Restriction to the superspecial locus}
Let $V$ be an $\fpbar$-vector space and let $\rho:\GU_g(\fptwo)\to\GL(V)$ be a representation.  A \emph{superspecial modular form}\index{modular form!superspecial}\footnote{This is not standard terminology.} of weight $\rho$ and level $N$ is a function $f:\Sigma\to V$ satisfying
$$
f([A,\lambda,\alpha,M\eta])=\rho(M)^{-1}f([A,\lambda,\alpha,\eta]),\quad\text{for all }M\in\GU_g(\fptwo).
$$
The space of all such forms will be denoted $S_\rho$.  If $\tau$ is a subrepresentation of $\rho$, then $S_\tau\subset S_\rho$.  If $\rho$ and $\tau$ are representations, then
$$
S_{\rho\otimes\tau}=S_\rho\otimes S_\tau.
$$

Let $\mathcal{I}$ denote the ideal sheaf of $i:\Sigma\into
X$, i.e. the kernel in:
$$
0\to\mathcal{I}\to\oh_X\to i_*\oh_\Sigma\to 0.
$$
The sheaf $\mathcal{I}$ is coherent (by Proposition II.5.9 of~\cite{hartshorne1}).  Given one of our sheaves $\mathbb{E}_\rho$, we obtain after tensoring and taking cohomology
$$
0\to\H^0(X,\mathcal{I}\otimes\mathbb{E}_\rho)\to\H^0(X,\mathbb{E}_\rho)\to\H^0(X,i_*\oh_\Sigma\otimes\mathbb{E}_\rho)=\H^0(\Sigma,i^*\mathbb{E}_\rho)\to\H^1(X,\mathcal{I}\otimes\mathbb{E}_\rho).
$$
We rewrite the part that interests us in a more familiar notation:
$$
0\to\H^0(X,\mathcal{I}\otimes\mathbb{E}_\rho)\to M_\rho(N)\xrightarrow{r} S_{\Res\rho} \to H^1(X,\mathcal{I}\otimes\mathbb{E}_\rho),
$$
where $\Res$ restricts representations on $\GL_g$ to the finite subgroup $\GU_g(\fptwo)$.

We'll use the following result to determine when the map $r$ is surjective:
\begin{thm}[Serre, see Theorem III.5.2 in~\cite{hartshorne1}]
\label{thm:very_ample_vanish}
Let $X$ be a projective scheme over a noetherian ring $A$, and let
$\mathcal{F}$ be a very ample invertible sheaf on $X$ over $\Spec A$.
Let $\mathcal{G}$ be a coherent sheaf on $X$.  Then
\begin{enumerate}
\renewcommand{\labelenumi}{(\alph{enumi})}
\item for each $i\geq 0$, $\H^1(X,\mathcal{G})$ is a finitely
generated $A$-module;
\item there is an integer $n_0$, depending on $\mathcal{G}$, such that
for each $i>0$ and each $n\geq n_0$
$$
\H^i(X,\mathcal{G}\otimes\mathcal{F}^n)=0.
$$
\end{enumerate}
\end{thm}

Let $\omega=\Lambda^g\mathbb{E}=\mathbb{E}_{\det}$; it is an ample invertible sheaf (see Theorem V.2.5 in~\cite{faltings1}).  Then for $n\gg 0$ we have
$$
\H^1(X,(\mathcal{I}\otimes\mathbb{E}_\rho)\otimes\omega^n)=0,
$$
so for $n\gg 0$, $r$ is a surjective map
$$
r:M_{\rho\otimes\det^n}(N)\to S_{\Res(\rho\otimes\det^n)}.
$$

\subsubsection{Lifting weights}
If $H$ is a subgroup of a group $G$, we say that a representation $\rho$ of $H$ \emph{lifts}\index{lift!of a representation} to $G$ if there exists a representation $\bar{\rho}$ of $G$ such that $\rho=\Res\bar{\rho}$.  It is clear that if $\rho$ lifts to $\bar{\rho}$ and $\tau$ lifts to $\bar{\tau}$, then $\rho\oplus\tau$ lifts to $\bar{\rho}\oplus\bar{\tau}$.  

Let $q$ be some power of $p$.  The following is a direct consequence of Theorems 6.1 and 7.4 from~\cite{steinberg1}:
\begin{prop}\label{prop:lift_sl}
Every irreducible representation of $\SL_g(\fq)$ lifts to a unique irreducible rational representation of $\SL_g(\fpbar)$.
\end{prop}

We now extend this to
\begin{prop}
Every irreducible representation of $\GL_g(\fq)$ lifts to an irreducible rational representation of $\GL_g(\fpbar)$.
\end{prop}
\begin{proof}
It suffices to prove that every irreducible representation lifts to a completely reducible one.  Let $\rho:\GL_g(\fq)\to\GL(V)$ be irreducible.

Via the canonical embeddings $\SL_g(\fq)\subset\GL_g(\fq)$ and $\gee_m(\fq)\subset\GL_g(\fq)$, $\rho$ induces representations $\rho_s:\SL_g(\fq)\to\GL(V)$ and $\rho_m:\gee_m(\fq)\to\GL(V)$, such that $\im\rho_s$ commutes with $\im\rho_m$.  Since $\GL_g(\fq)=\SL_g(\fq)\cdot\gee_m(\fq)$ and $\SL_g(\fq)\cap\gee_m(\fq)=\boldsymbol{\mu}_g(\fq)$, we also have that $\rho_s(\zeta)=\rho_m(\zeta)$ for all $\zeta\in\boldsymbol{\mu}_g(\fq)$.

Any representation of $\gee_m(\fq)$ is of the form
\begin{eqnarray*}
\gee_m(\fq) &\to& \GL(V)\\
\lambda &\mapsto& \left(\begin{smallmatrix}\lambda^{a_1}\\&\ddots\\&&\lambda^{a_n}\end{smallmatrix}\right)
\end{eqnarray*}
with $a_i\in\zed/(q-1)\zed$.  I claim that in our case $\gee_m(\fq)$ acts by scalars on $V$.  Suppose this is false, then there exists $\lambda\in\gee_m(\fq)$ such that at least two of the diagonal entries of $\rho_m(\lambda)$ are distinct.  By changing the basis of $V$ we can assume $\rho_m(\lambda)$ is in Jordan canonical form.  Let $A\in\SL_g(\fptwo)$, then the fact that $\rho_s(A)$ commutes with $\rho_m(\lambda)$ forces $A$ to have the same shape as $\rho_m(\lambda)$ (i.e. it is block-diagonal with blocks of the same dimensions as $\rho_m(\lambda)$).  Since this holds for all $A\in\SL_g(\fq)$, we conclude that as an $\SL_g(\fq)$-module, $V$ has a direct sum decomposition
$$
V=V_1\oplus\ldots\oplus V_j
$$
corresponding to the shape of $\rho_m(\lambda)$ (in the chosen basis for $V$, $V_1$ is the span of the first $k$ vectors, where $k$ is the size of the first Jordan block of $\rho_m(\lambda)$, etc.).  But this means that $V_1$ is a proper subspace of $V$ which invariant under both $\SL_g(\fq)$ and $\gee_m(\fq)$, contradicting the hypothesis that $V$ is an irreducible representation of $\GL_g(\fq)$.  So $\gee_m(\fq)$ acts by scalars on $V$, say $\rho_m(\lambda)v=\lambda^av$ for some $a\in\zed/(q-1)\zed$.

From this it is clear that $\rho_m$ is completely reducible and that any choice of $\bar{a}\in\zed$ with $\bar{a}\equiv a\pmod{q-1}$ yields a completely reducible lift $\bar{\rho}_m:\gee_m(\fpbar)\to\GL(V)$ given simply by $\lambda\mapsto\lambda^{\bar{a}}$.  Note that $\bar{\rho}_m$ is a rational representation.  Later on we'll need to choose a lift of $a$ to $\bar{a}\in\zed$ that suits us better.

It is also pretty clear that $\rho_s$ is irreducible: if $W$ is an irreducible $\SL_g(\fq)$-submodule, then $W$ is also $\gee_m(\fq)$-invariant so it is $\GL_g(\fq)$-invariant, hence either $W=0$ or $W=V$.

By Proposition~\ref{prop:lift_sl}, $\rho_s$ lifts to an irreducible rational $\bar{\rho}_s:\SL_g(\fpbar)\to\GL(V)$.  Since $\gee_m$ acts by scalars, $\im\bar{\rho}_m$ commutes with $\im\bar{\rho}_s$.  I claim that the maps $\bar{\rho}_m$ and $\bar{\rho}_s$ agree on $\boldsymbol{\mu}_g(\fpbar)=\SL_g(\fpbar)\cap\gee_m(\fpbar)$.  Assuming this is true, we can construct a rational representation 
\begin{eqnarray*}
\bar{\rho}:\GL_g(\fpbar)&\to& \GL(V)\\
M &\mapsto& \bar{\rho}_m(\det M)\cdot\bar{\rho}_s\left((\det M)^{-1}M\right).
\end{eqnarray*}
Since the restriction of $\bar{\rho}$ to $\SL_g(\fpbar)$ is $\bar{\rho}_s$ and in particular irreducible, we conclude that $\bar{\rho}$ is irreducible.

It remains to prove that $\bar{\rho}_m$ and $\bar{\rho}_s$ agree on the $g$-th roots of unity.  It suffices to do this for a primitive $g$-th root $\zeta$.  Write $g=p^sg'$ with $(p,g')=1$.  We have $(\zeta^{g'})^{p^s}=\zeta^g=1$, so $\zeta^{g'}=1$ since the only $p^s$-th root of unity in characteristic $p$ is $1$.  Therefore $\zeta$ is a $g'$-th root of unity, so without loss of generality we may assume that $(p,g)=1$.

Consider the linear transformation $\bar{\rho}_s(\zeta)$.  It is diagonalizable if and only if its minimal polynomial has distinct roots.  But the transformation satisfies $X^g-1=0$, which has distinct roots, and hence the minimal polynomial will also have distinct roots.  So we can choose a basis for $V$ such that $\bar{\rho}_s(\zeta)$ is diagonal.  If it has at least two distinct diagonal entries, we can apply the same argument as before to conclude that since it commutes with all of $\bar{\rho}_s(\SL_g(\fpbar))$ the representation $\bar{\rho}_s$ is reducible, which is a contradiction.  So 
$$
\bar{\rho}_s(\zeta)=\zeta^b,\quad\quad\text{for some }b\in\zed/g\zed.
$$
We want to show that $\bar{\rho}_m(\zeta)=\bar{\rho}_s(\zeta)$, i.e. that we can choose $\bar{a}\in\zed$ such that $\bar{a}\equiv b\pmod{g}$.  Let $d=(g,q-1)$ and write $g=dm$, $q-1=dn$.  We have $(\zeta^m)^d=\zeta^g=1$ so $(\zeta^m)^{q-1}=(\zeta^{md})^n=1$ so $\zeta^m\in\fq$.  Therefore $\zeta^m\in\boldsymbol{\mu}_g(\fq)$ and hence
$$
(\zeta^b)^m=\bar{\rho}_s(\zeta^m)=\bar{\rho}_m(\zeta^m)=(\zeta^m)^{\bar{a}}.
$$
This implies that $m\bar{a}\equiv mb\pmod{g}$, i.e. $\bar{a}\equiv b\pmod{d}$.  Since $d=(g,q-1)$ and $d|(\bar{a}-b)$ there exist integers $u,v$ such that $\bar{a}-b=ug+v(q-1)$ and therefore
$$
(\bar{a}-v(q-1))\equiv b\pmod g,
$$
which is what we wanted.
\end{proof}

Note that in contrast with Proposition~\ref{prop:lift_sl} the lift of $\rho$ to $\GL_g(\fpbar)$ is not unique.  Fix some lift $\bar{\rho}$, then any lift can be written in the form $\det^m\otimes\bar{\rho}$, where $m$ is a common multiple of $g$ and $q-1$.

\begin{cor}\label{cor:lift_weight}
Given an irreducible representation $\tau:\GU_g(\fptwo)\to\GL(W)$, there exists an irreducible rational representation $\bar{\rho}:\GL_g(\fpbar)\to\GL(V)$ such that $\tau\subset\Res\bar{\rho}$.
\end{cor}
\begin{proof}
Consider the induced representation from $\GU_g(\fptwo)$ to $\GL_g(\fptwo)$.  This has an irreducible subrepresentation $\rho:\GL_g(\fptwo)\to\GL(V)$ with the property that $\tau\subset\Res\rho$.  The result now follows from the previous proposition.
\end{proof}

\subsubsection{Proof of the main result}
We have come to the main result of the section.  Recall the notation $U_\ell(N)=\GSp_{2g}(\zed_\ell)(N)$ for $\ell\neq p$, $U_p=\ker(\GU_g(\oh_p)\to\GU_g(\fptwo))$ and
$$
U=U_p\times\prod_{\ell\neq p} U_\ell(N).
$$
\begin{thm}
Fix a dimension $g>1$, a level $N\geq 3$ and a prime $p$ not dividing $N$.  The systems of Hecke eigenvalues coming from Siegel modular forms (mod $p$) of dimension $g$, level $N$ and any weight $\rho$, are the same as the systems of Hecke eigenvalues coming from algebraic modular forms (mod $p$) of level $U$ and any weight $\rho_\Sigma$ on the group $\GU_g(B)$.
\end{thm}
\begin{proof}
Let $f$ be a Siegel modular form of weight $\rho:\GL_g\to\GL_m$ which is a Hecke eigenform.  If $r(f)=0$, then $f\in\H^0(X,\mathcal{I}\otimes\mathbb{E}_\rho)$.  The quotient map of $\oh_X$-modules $\mathcal{I}\to\mathcal{I}/\mathcal{I}^2$ induces (after tensoring with $\mathbb{E}_\rho$ and taking global sections) a map 
$$
\bar{\cdot}:\H^0(X,\mathcal{I}\otimes\mathbb{E}_\rho)\to\H^0(X,\mathcal{I}/\mathcal{I}^2\otimes\mathbb{E}_\rho).
$$
Consider $\bar{f}\in\H^0(X,\mathcal{I}/\mathcal{I}^2\otimes\mathbb{E}_\rho)$.  We have an exact sequence
$$
0\to\mathcal{I}\otimes\mathcal{I}/\mathcal{I}^2\otimes\mathbb{E}_\rho\to\mathcal{I}/\mathcal{I}^2\otimes\mathbb{E}_\rho\to i_*\oh_\Sigma\otimes\mathcal{I}/\mathcal{I}^2\otimes\mathbb{E}_\rho\to 0
$$
which gives us a long exact sequence that starts with
$$
0\to\H^0(X,\mathcal{I}^2/\mathcal{I}^3\otimes\mathbb{E}_\rho)\to\H^0(X,\mathcal{I}/\mathcal{I}^2\otimes\mathbb{E}_\rho)\xrightarrow{r_1}\H^0(\Sigma,i^*(\mathcal{I}/\mathcal{I}^2\otimes\mathbb{E}_\rho)).
$$
If $r_1(\bar{f})=0$ then $\bar{f}\in\H^0(X,\mathcal{I}^2/\mathcal{I}^3\otimes\mathbb{E}_\rho)$ and we can similarly consider $r_2(\bar{f})$, $r_3(\bar{f})$ etc.  There exists some $n$ such that $r_n(\bar{f})\neq 0$.  Let $f_S=r_n(\bar{f})\in\H^0(\Sigma,i^*(\mathcal{I}^n/\mathcal{I}^{n+1}\otimes\mathbb{E}_\rho))$.  Note that $\mathcal{I}^n/\mathcal{I}^{n+1}=\Sym^n(\mathcal{I}/\mathcal{I}^2)$ and that $i^*(\mathcal{I}/\mathcal{I}^2)=i^*(\Omega^1_X)$.  Recall from \S{}\ref{sect:kodaira-spencer} the Kodaira-Spencer isomorphism $\Omega^1_X\cong\mathbb{E}_{\Sym^2\std}$.  We conclude that $f_S\in S_{\Res((\Sym^{2n}\std)\otimes\rho)}$.  So our process associates to a Siegel modular form $f$ of weight $\rho$ a superspecial modular form $f_S$ of weight $\Res((\Sym^{2n}\std)\otimes\rho)$ for some integer $n$ depending on $f$.  Moreover, since the restrictions $r_i$ and the Kodaira-Spencer isomorphism are Hecke maps, we conclude that $f_S$ is a Hecke eigenform with the same eigenvalues as $f$.

Now let $f_S$ be a superspecial Siegel modular form of weight $\rho_S:\GU_g(\fptwo)\to\GL_m(\fpbar)$.  By applying Corollary~\ref{cor:lift_weight} we get a rational representation $\bar{\rho}:\GL_g\to\GL_m$ such that $\rho_S\subset\Res\bar{\rho}$.  By functoriality we get $S_{\rho_S}\subset S_{\Res\bar{\rho}}$.  We know that the map
$$
r:M_{\bar{\rho}\otimes\det^n}(N)\to S_{\Res(\bar{\rho}\otimes\det^n)}
$$
is surjective for $n\gg 0$, and therefore there exists an integer $k$ such that
$$
r:M_{\bar{\rho}\otimes\det^{k(p^2-1)}}(N)\to S_{\Res(\bar{\rho}\otimes\det^{k(p^2-1)})}=S_{\Res\bar{\rho}}\supset S_{\rho_S}
$$
is surjective.  Since this map is also Hecke-invariant, we conclude from Proposition 1.2.2 of~\cite{ash1} that any system of Hecke eigenvalues that occurs in $S_{\rho_S}$ also occurs in $M_{\bar{\rho}\otimes\det^{k(p^2-1)}}$.

So far we showed that the systems of Hecke eigenvalues given by Siegel modular forms (mod $p$) of all weights are the same as the systems of Hecke eigenvalues given by superspecial modular forms $S_{\rho_S}$ of all weights.  By Theorem~\ref{thm:main_bijection} we know that $S_{\rho_S}$ is isomorphic as a Hecke module to the space of algebraic modular forms (mod $p$) of weight $\rho_S$, and we're done.
\end{proof}

\subsubsection{Agreement with the definition of Gross}
In this section we'll write $G=\GU_g(\fptwo)$.

Recall from \S{}\ref{sect:algebraic_modular} that Gross defines algebraic modular forms (mod $p$) as follows: let $\rho:G\to\GL(V)$ be an irreducible representation where $V$ is a finite-dimensional vector space \emph{over $\fp$}, then set
$$
M(\rho):=\{f:\Omega\to V|f(\lambda x)=\rho(\lambda)^{-1}f(x)\text{ for all }\lambda\in G\}.
$$

For comparison, our spaces of modular forms on $\Omega$ are defined as
$$
M(\tau):=\{f:\Omega\to W|f(\lambda x)=\rho(\lambda)^{-1}f(x)\text{ for all }\lambda\in G\},
$$
where $\tau:G\to\GL(W)$ is an irreducible representation and $W$ is a finite-dimensional vector space \emph{over $\fpbar$}.

The purpose of this section is to show that the spaces $M(\rho)$ and $M(\tau)$ for varying $\rho$ and $\tau$ give the same systems of Hecke eigenvalues.

First suppose that $(a_T:T)$ is a system of Hecke eigenvalues coming from $M(\rho)$.  Then there exists $f\in M(\rho)\otimes\fpbar$ such that $T(f)=a_Tf$ for all $T$.  Let $\rho\otimes\fpbar$ denote the composition $G\xrightarrow{\rho}\GL(V)\into\GL(V\otimes\fpbar)$.  The map
\begin{eqnarray*}
M(\rho)\otimes\fpbar & \to & M(\rho\otimes\fpbar) \\
m\otimes\alpha & \mapsto & \alpha m
\end{eqnarray*}
is an isomorphism compatible with the action of the Hecke operators, so the image of $f$ in $M(\rho\otimes\fpbar)$ is an eigenform with the same eigenvalues as $f$.  Therefore the system $(a_T)$ also comes from $M(\rho\otimes\fpbar)$.

Conversely, suppose that $(a_T:T)$ is a system of Hecke eigenvalues coming from $M(\tau)$ for some $\tau:G\to\GL(W)$, $W$ a finite-dimensional $\fpbar$-vector space.  Then there exists $f\in M(\tau)$ such that $T(f)=a_T f$ for all $T$.  Since $G$ is a finite group there exist $q=p^a$, a finite-dimensional $\fq$-vector space $W'$ and a representation $\tau':G\to\GL(W')$ such that $\tau'\otimes\fpbar=\tau$.  Similarly, $\Omega$ is a finite set and $f$ is a map $\Omega\to W$ so by enlarging $q$ if necessary, there exists $f'\in M(\tau')$ such that $f$ is the image of $f'\otimes 1$ under the isomorphism $M(\tau')\otimes\fpbar\cong M(\tau)$.  Clearly $T(f')=a_Tf'$ for all $T$; in particular $a_T\in\fq$ for all $T$.

We now use the following
\begin{prop}\label{prop:galois_eigen}
Suppose $L/K$ is a finite Galois extension with Galois group $G$ and $V$ is a finite-dimensional vector space over $L$.  Let $\mathcal{T}$ be a collection of commuting diagonalizable linear operators on $V$ and let $V_K$ be the space $V$ viewed as a vector space over $K$.  If a $\mathcal{T}$-eigenvector $v$ has system of eigenvalues $\{a_T:T\in\mathcal{T}\}$, then for every $\sigma\in G$ there exists an eigenvector $v_\sigma\in V_K$ with system of eigenvalues $\{\sigma(a_T):T\in\mathcal{T}\}$.
\end{prop}

Let's first see how this concludes our argument.  We apply the proposition to the finite Galois extension $\fq/\fp$, the vector space $M(\tau')$, the Hecke operators $T$, the eigenvector $f'$ and the identity Galois element $\sigma=1$.  We conclude that if we consider $M(\tau')$ as a vector space over $\fp$, there exists an eigenvector $f''$ with the same system of eigenvalues as $f'$.  This is precisely what we needed to show.

\begin{proof}[Proof of Proposition~\ref{prop:galois_eigen}]
The isomorphism $\varphi$ of the next lemma induces an isomorphism of $L$-vector spaces
\begin{eqnarray*}
\varphi:L\otimes_K V &\to& \bigoplus_{\sigma\in G} Ve_\sigma\\
\alpha\otimes w &\mapsto& \sum_{\sigma\in G} \sigma(\alpha)w e_\sigma.
\end{eqnarray*}
Let $v_\sigma:=\varphi^{-1}(ve_{\sigma^{-1}})$.  We have
$$
Tv_\sigma=\varphi^{-1}((Tv)e_{\sigma^{-1}})=\varphi^{-1}((a_Tv)e_{\sigma^{-1}})=\sigma(a_T)\varphi^{-1}(ve_{\sigma^{-1}})=\sigma(a_T)v_\sigma,
$$
so $v_\sigma$ is an eigenvector of $T$ with eigenvalue $\sigma(a_T)$, and this holds for all $T\in\mathcal{T}$.
\end{proof}

\begin{lem}
Suppose $L/K$ is a finite Galois extension with Galois group $G$.  The map
$$
\varphi:L\otimes_K L\to\bigoplus_{\sigma\in G} Le_\sigma
$$
defined by $\alpha\otimes\beta\mapsto \sum_{\sigma\in G} \sigma(\alpha)\beta e_\sigma$ is an isomorphism of $L$-algebras.
\end{lem}
\begin{proof}
It is pretty clear that $\varphi$ is an $L$-algebra homomorphism.  Since the dimensions of the domain and of the range are equal (and equal to $[L:K]$), it suffices to prove that $\varphi$ is injective.

Let $\{\alpha_1,\ldots,\alpha_n\}$ be a basis of $L$ as a $K$-vector space.  Then $\{\alpha_i\otimes\alpha_j:1\leq i,j\leq n\}$ is a basis of $L\otimes_K L$ as a $K$-vector space.  Suppose $\varphi(\sum c_{ij}\alpha_i\otimes\alpha_j)=0$.  If we write $G=\{\sigma_1,\ldots,\sigma_n\}$, then we have
\begin{equation}\label{eq:indep_char}
\sum_{i,j}c_{ij}\sigma_k(\alpha_i)\alpha_j=0\quad\text{for all }k.
\end{equation}
Let $A$ be the $n\times n$ matrix whose $(i,j)$-th entry is $\sigma_i(\alpha_j)$, and let $c$ be the column vector whose $i$-th entry is $\sum_j c_{ij}\alpha_j$.  Then the system~(\ref{eq:indep_char}) can be written as $Ac=0$.  But it is an easy consequence of independence of characters (see Corollary VI.5.4 in~\cite{lang1}) that $A\in\GL_n(L)$, therefore we must have $c=0$, i.e.
$$
\sum_j c_{ij}\alpha_j=0\quad\text{for all }i.
$$
Since the $\alpha_j$ are linearly independent we conclude that $c_{ij}=0$ for all $i$ and $j$, hence $\varphi$ is injective.
\end{proof}

\bibliographystyle{alpha}
\bibliography{master}
\printindex
\end{document}